\documentclass[11pt]{article}
\usepackage[T1]{fontenc}
\usepackage[utf8]{inputenc}
\usepackage{lmodern}
\usepackage[margin=1in]{geometry}
\usepackage{amsfonts,amsmath,amssymb,amsthm,mathtools,mathrsfs,times,color}
\usepackage{microtype}
\usepackage{enumitem,booktabs,array}
\usepackage[colorlinks,linkcolor=blue,anchorcolor=blue,citecolor=blue]{hyperref}
\hypersetup{pdftitle={Calabi-Energy Sublevels, K-Stability, and Calabi Flow},
pdfauthor={Haozhao Li and Bing Wang}}
\numberwithin{equation}{section}
\newtheorem{theorem}{Theorem}[section]
\newtheorem{proposition}[theorem]{Proposition}
\newtheorem{lemma}[theorem]{Lemma}
\newtheorem{corollary}[theorem]{Corollary}
\theoremstyle{definition}

\newtheorem{remark}[theorem]{Remark}
\newcommand{\Ca}{\operatorname{Ca}}
\newcommand{\AM}{\operatorname{AM}}
\newcommand{\Ric}{\operatorname{Ric}}
\newcommand{\Rm}{\operatorname{Rm}}
\newcommand{\tr}{\operatorname{tr}}
\newcommand{\Vol}{\operatorname{Vol}}
\newcommand{\Aut}{\operatorname{Aut}}
\newcommand{\DF}{\operatorname{DF}}
\newcommand{\ddc}{\sqrt{-1}\partial\bar\partial}
\newcommand{\Mab}{\mathcal M}
\newcommand{\HH}{\mathcal H}
\newcommand{\EE}{\mathcal E}
\newcommand{\CS}{C_{\mathrm S}}
\newcommand{\CE}{C_{\mathrm E}}
\newcommand{\norm}[1]{\left\|#1\right\|}
\allowdisplaybreaks[1]
\title{Calabi-energy sublevels, K-stability and Calabi flow}
\author{Haozhao Li\footnote{Supported by NSFC grant No. 12471058 and a
research fund from University of Science and Technology of China.}
\quad and \quad Bing Wang\footnote{Supported by NSFC grant No. 12431003 and a
research fund from University of Science and Technology of China.}}
\date{}
\begin{document}
\maketitle
\begin{abstract} In this paper, we prove path connectedness of every nonempty strict Calabi-energy sublevel and identify
the infimum of the square root of the Calabi energy over smooth metrics
with the infimum of the descending Mabuchi slope in the metric completion.
In particular, there exist smooth metrics of arbitrarily small Calabi
energy in every K-semistable real ample class. On compact K\"ahler surfaces,
we prove the existence of constant scalar curvature K\"ahler metrics in the
K-polystable classes of the controlled cone. In any class containing such a metric, the
Calabi flow exists smoothly for all time with a uniform Sobolev bound whenever the
initial metric is smooth and has sufficiently small Calabi energy.
\end{abstract}
\begingroup
\small
\setcounter{tocdepth}{2}
\tableofcontents
\endgroup
\section{Introduction}\label{sec:intro}

A basic problem in complex differential geometry is to find a canonical metric in a
fixed K\"ahler class. Calabi~\cite{CalabiExtremal}
formulated a variational approach through the squared $L^2$ norm of
scalar curvature. Constant scalar curvature K\"ahler metrics are
absolute minimizers whenever they exist. Even when it is unknown whether the infimum is
attained, we can still study its value and the geometry of low-energy metrics.

The relation with stability provides another source of questions.
For polarized manifolds, Donaldson~\cite[Theorem 1]{DonaldsonCalabi} obtained
lower bounds for the Calabi functional from test configurations.
This links scalar-curvature error to degenerations of the complex
manifold. On the analytic side, low-energy metrics can be studied through the Mabuchi functional and its negative gradient flow, the Calabi flow.
We ask how this geometry controls the topology of energy sublevels,
the smooth energy infimum, and the evolution of small-energy metrics.

J.~Streets~\cite[Lemma 4.9, equation (4.4)]{Streets} derived the
resolvent equation by computing the first
variation of the penalized Mabuchi functional at a smooth potential.
This is the Euler--Lagrange equation for one minimizing step of the
Calabi flow. We construct classical solutions of this equation in the
fixed K\"ahler class and use them to prove that the strict sublevels are
path connected and that the smooth Calabi infimum, in square-root form,
equals the weak Mabuchi-slope infimum. K-semistability gives arbitrarily
small smooth Calabi energy in real ample classes. On compact K\"ahler
surfaces, K-polystability in the controlled cone gives constant scalar
curvature K\"ahler metrics, abbreviated as cscK metrics. If the prescribed class on a
compact K\"ahler surface contains a cscK metric, we also obtain a Sobolev bound and global smooth
Calabi flow when the initial energy is sufficiently small.

Let $(X^n,J,\omega_*)$ be a compact connected K\"ahler manifold,
$n\ge1$. We use the conventions
\[
 V=\int_X\omega_*^n,\qquad \omega_u=\omega_*+\ddc u,
 \qquad d\mu_u=V^{-1}\omega_u^n.
\]
Let $S_u$ denote complex scalar curvature, $\underline S$ its class average,
and $\Delta_u$ the nonpositive complex Laplacian. The centered, volume-normalized Calabi energy is defined by
\begin{equation}\label{eq:energy-definition}
 \Ca(u)=\int_X(S_u-\underline S)^2d\mu_u.
\end{equation}
Thus $\Ca(u)=0$ exactly when $\omega_u$ is cscK. In a fixed class,
centering subtracts a constant from the uncentered Calabi functional,
and probability normalization changes it by the fixed volume factor.
Write $\HH=\{u\in C^\infty(X,\mathbb R):\omega_u>0\}$ and normalize
additive constants by
\[
 \AM(u)=\frac1{(n+1)V}\sum_{j=0}^n
              \int_Xu\,\omega_u^j\wedge\omega_*^{n-j}=0.
\]
Put $\HH_0=\HH\cap\AM^{-1}(0)$ and, for $0<\alpha<1$, put
\[
 \HH^{4,\alpha}_0
 =\{u\in C^{4,\alpha}(X,\mathbb R):\omega_u>0,\ \AM(u)=0\}.
\]

\subsection{Sublevels, energy infima, and stability}
\label{sec:intro-resolvents}

Let $d_2$ be the Mabuchi distance for the tangent norm
$\|\xi\|_{2,u}^2=\int_X\xi^2d\mu_u$.
Darvas~\cite[Theorem 1]{DarvasCompletion} identified the completion of
$(\HH,d_2)$ with the space $\EE^2$ of finite-energy potentials.
Mabuchi~\cite{MabuchiKenergy} introduced the K-energy $\Mab$, whose
first variation at a smooth potential is
\[
 d\Mab_u(\xi)=-\int_X\xi(S_u-\underline S)d\mu_u.
\]
Thus the smooth critical points are precisely the cscK potentials, and
the Calabi flow is the negative gradient flow of $\Mab$.
Streets~\cite[Theorem 1.3]{Streets} constructed global minimizing
movements in the metric completion. Berman, Darvas, and
Lu~\cite[Theorem 4.7 and Section 6]{BDLflow} constructed the lower
semicontinuous extension of $\Mab$ to $\EE^2$ and identified its
gradient flow with these minimizing movements.
At smooth potentials, Xia~\cite[Proposition 3.6]{Xia} showed that
the descending metric slope $|\partial\Mab|$ equals $\sqrt{\Ca}$.

The first question concerns the geometry of the sublevels themselves.
Affine interpolation joins any two smooth potentials, but the scalar curvature may
vary along the resulting path. The following theorem gives a path along which the
energy stays below the prescribed level.

\begin{theorem}[Path connectedness]\label{thm:connected}
Let $(X^n,J,\omega_*)$ be a compact connected K\"ahler manifold. Fix $0<\alpha<1$ and $c\in\mathbb R$. Consider the strict sublevel
\[
 \mathcal U_c^\alpha
 =\{u\in C^{4,\alpha}(X,\mathbb R):\ \omega_u>0,\ \AM(u)=0,\ \Ca(u)<c\}.
\]
 If $\mathcal U_c^\alpha$ is nonempty,
it is path connected in the $C^{4,\alpha}$ topology. Any two smooth
potentials in $\mathcal U_c^\alpha$ can be joined within this sublevel
by a path of smooth potentials continuous in the $C^\infty$ topology.
\end{theorem}

Any two metrics in the same strict sublevel can therefore be deformed
into one another while retaining a uniform bound below its threshold.
This gives a way to transfer continuous geometric inequalities along
low-energy paths. In Section~\ref{sec:sobolev-application}, we use it
to transfer a Sobolev bound from a cscK reference metric.

The next question is how small the energy can be. At a smooth potential, the descending Mabuchi
slope equals $\sqrt{\Ca}$, so its infimum in the
completion gives a lower bound for the smooth infimum.

\begin{theorem}[Smooth and weak Calabi infima]
\label{thm:smooth-weak-infima}
On every compact connected K\"ahler manifold $(X^n,J,\omega_*)$, we have
\begin{equation}\label{eq:intro-smooth-weak-infima}
 \inf_{u\in\HH}\sqrt{\Ca(u)}
 =\inf_{w\in\EE^2}|\partial\Mab|(w).
\end{equation}
In particular, if the right-hand side vanishes, then for every
$\varepsilon>0$ there is a smooth positive K\"ahler metric
$\omega_u\in[\omega_*]$ with $\Ca(u)<\varepsilon$.
\end{theorem}

The equality identifies the least scalar-curvature error among smooth
metrics with the least Mabuchi slope in the metric completion.
In particular, arbitrarily small weak slopes give smooth minimizing
sequences in the original class. Xia~\cite[Theorem 1.1(1)]{Xia}
identifies the weak infimum with the maximum of the normalized
negative Mabuchi slopes over nonconstant finite-energy geodesic rays.
Xia~\cite[Remark 4.2]{Xia} obtains the smooth infimum identity
under global smooth existence of Calabi flow. Here the classical
resolvent construction gives this identity in every compact
K\"ahler class.

This infimum comparison connects the metric theory to stability.
For real classes, we use cohomological K-stability as in
Sj\"ostr\"om Dyrefelt~\cite{SDSemi} and
Dervan and Ross~\cite{DRKahler}. K-semistability requires a nonnegative
Donaldson--Futaki invariant for every normal relatively K\"ahler test
configuration whose class restricts to $\Omega$ on the general fibre.
K-polystability requires every zero-invariant configuration to be a
product, allowing a nontrivial $\mathbb C^*$-action on $X$.
For a polarization, we test normal relatively ample algebraic test
configurations of every positive exponent. In the anticanonical class
of a Fano manifold, Li~\cite[Theorem 3]{LiSemiFano} proved that
K-semistability is equivalent to the vanishing of the smooth Calabi
infimum. The next corollary gives the implication from semistability
for any real ample class.

\begin{corollary}[Small smooth Calabi energy from semistability]
\label{cor:semistable-small-energy}
Let $X$ be a smooth connected projective manifold of any complex
dimension, and let $\Omega\in\operatorname{NS}(X)_{\mathbb R}$
be a real ample class. If $(X,\Omega)$ is cohomologically
K-semistable, then
\begin{equation}\label{eq:intro-semistable-small-energy}
 \inf_{\omega\in2\pi\Omega,\ \omega\ \text{smooth K\"ahler}}
       \Ca(\omega)=0.
\end{equation}
Thus the original class contains a smooth metric with energy below any prescribed
positive threshold. For $\Omega=c_1(L)$, ordinary algebraic
K-semistability of $(X,L)$ suffices.
\end{corollary}

The corollary produces smooth minimizing sequences in each fixed real
ample class satisfying semistability. This lets us choose an initial
metric below any positive energy threshold arising in a geometric
existence criterion.

On a compact K\"ahler surface, the controlled cone provides such a
positive threshold.
By the flow criterion of Li, Wang, and Zheng~\cite[Corollary 1.3]{LWZSobolev}, the
resulting Calabi flow exists smoothly for all time with uniformly bounded Sobolev
constants. Write $c_1=c_1(X)$, $[\omega]=2\pi\Omega$, and set
\begin{equation}\label{eq:controlled-cone}
 \mathcal U=\left\{\Omega\in\mathcal K_X:
 c_1\cdot\Omega>0,\quad
 \delta(\Omega):=c_1^2-
 \frac{2(c_1\cdot\Omega)^2}{3\Omega^2}>0\right\}.
\end{equation}
Within this cone we can ask when the flow yields a cscK metric on
the original complex manifold $(X,J)$. The next theorem gives a sufficient
condition in terms of K-polystability.

\begin{theorem}[Stability and cscK existence in the controlled cone]
\label{thm:stability}
Let $X$ be a smooth compact connected K\"ahler surface and let
$\Omega\in\mathcal U$ be a real K\"ahler class.
If $(X,\Omega)$ is cohomologically K-polystable, then
$2\pi\Omega$ contains a smooth cscK metric on $(X,J)$.
When $\Aut^0(X,J)=\{1\}$, existence is equivalent to
cohomological K-stability.
For $\Omega=c_1(L)$, existence is equivalent to ordinary
algebraic K-polystability of $(X,L)$.
\end{theorem}

Here semistability supplies the small-energy initial metric, while
polystability identifies the geometric limit with $(X,J)$ in the
prescribed K\"ahler class. The cone inequalities provide the analytic
control required by the degeneration theorem of
Li, Wang, and Zheng~\cite{ControlledCompanion}.
Thus the resulting metric lies in the prescribed real class.
Remark~\ref{rem:controlled-cone-role} details the role of these inequalities.

\subsection{Sobolev control and global Calabi flow}
\label{sec:intro-sobolev}

Let
$X$ be a compact K\"ahler surface with cscK reference metric $\omega_*$.
Chen and He~\cite[Theorem 1.2]{ChenHeStability} proved stability of
Calabi flow for potentials sufficiently close to a cscK metric in
$C^{3,\alpha}$. To study an energy sublevel in the whole class, we
seek estimates depending on the initial Calabi energy.
On a compact K\"ahler surface, a uniform Sobolev bound suffices to extend a smooth Calabi flow past any finite
time by the theorem of Li, Wang, and Zheng~\cite[Theorem 1.2]{LWZSobolev}. Set
\[
 D_u(f)=\int_X|\partial f|_u^2d\mu_u
       =-\int_Xf\Delta_uf\,d\mu_u,
\]
and define the optimal squared Sobolev constant by
\begin{equation}\label{eq:sobolev-definition}
 \CS(\omega_u)=\sup_{0\ne f\in W^{1,2}(X,\omega_u)}
 \frac{\left(\int_X|f|^4d\mu_u\right)^{1/2}}
 {D_u(f)+\int_Xf^2d\mu_u}.
\end{equation}
The lower-order coefficient is fixed, and $\CS\ge1$.

\begin{theorem}[Small energy controls the Sobolev constant]
\label{thm:sobolev-small}
Let $(X,J)$ be a compact connected K\"ahler surface with cscK metric $\omega_*$. For every
$A>\CS(\omega_*)$ there is
$\varepsilon_A=\varepsilon_A(X,J,\omega_*,A)>0$ such that every
smooth K\"ahler metric $\omega_u\in[\omega_*]$ satisfies
\[
 \Ca(u)<\varepsilon_A\quad\Longrightarrow\quad\CS(\omega_u)<A.
\]
In particular, sufficiently small Calabi energy gives
$\CS(\omega_u)<2\CS(\omega_*)$.
\end{theorem}

Since the theorem holds for every $A>\CS(\omega_*)$, it follows that
\[
 \Ca(u_j)\to0\quad\Longrightarrow\quad
 \limsup_{j\to\infty}\CS(\omega_{u_j})\le\CS(\omega_*).
\]
Tosatti and Weinkove~\cite[Theorems 1 and 2]{TosattiWeinkoveSmallEnergy}
proved long-time existence and smooth convergence of the Calabi flow
for sufficiently small initial Calabi energy under additional curvature
bounds. When $c_1(X)=0$, they assumed two-sided bounds on the initial
Ricci form. When $c_1(X)<0$, they treated classes near the canonical
class and assumed an upper bound on the initial Ricci form and a lower
bound on the initial scalar curvature. On a K\"ahler surface in a class
containing a cscK metric, Theorem~\ref{thm:sobolev-small} gives the
Sobolev bound from the initial energy alone. Since the energy decreases
along the smooth flow, this bound persists and the extension theorem
of Li, Wang, and Zheng~\cite[Theorem 1.2]{LWZSobolev} gives the
following corollary.

\begin{corollary}[Global Calabi flow]\label{cor:calabi-global}
Let $(X,J)$ be a compact connected K\"ahler surface with cscK metric $\omega_*$.
For every $A>\CS(\omega_*)$, there is
$\varepsilon_A=\varepsilon_A(X,J,\omega_*,A)>0$ such that every
smooth K\"ahler metric $\omega_{u_0}\in[\omega_*]$ with
$\Ca(u_0)<\varepsilon_A$ generates a smooth Calabi flow
\[
 \partial_tu=S_u-\underline S,\qquad u(0)=u_0,
\]
for all $t\ge0$. Along this flow,
\[
 \Ca(u(t))\le\Ca(u_0)<\varepsilon_A\quad(t\ge0),\qquad
 \sup_{t\ge0}\CS(\omega_{u(t)})\le A.
\]
\end{corollary}

Every initial metric in this strict energy sublevel therefore has a
global smooth evolution with the same uniform Sobolev bound.
The proof uses two properties of the sublevel. Paths in the sublevel carry the Sobolev
bound from the reference metric to the initial metric, and energy dissipation keeps
the flow in that sublevel.

\subsection{Outline of the proofs}\label{sec:intro-proofs}

To prove Theorem~\ref{thm:connected} for smooth endpoints in the strict
sublevel, we start with the normalized affine path between them.
For a single sufficiently large parameter, the resolvent maps the entire path below the
energy threshold. At each end we let the parameter decrease to zero to recover
the original endpoint. Positivity and the strict energy bound are open
conditions in $C^4$. We can therefore approximate the resulting path
by affine segments with smooth vertices and normalize each segment.

For Theorem~\ref{thm:smooth-weak-infima}, we instead compare the energy
of a resolvent with the Mabuchi slope at a finite-energy competitor.
The upper bound for the smooth infimum follows by approximating the resolvent by smooth
potentials and then letting the parameter tend to infinity. The reverse bound follows from the slope identity at smooth
potentials.

This comparison reduces Corollary~\ref{cor:semistable-small-energy} to
showing that the weak slope infimum vanishes. Entropy approximation in
the fixed real ample class extends the semistability inequality from
model metrics to finite-energy non-Archimedean potentials. Combining the ray comparison
of Mesquita-Piccione~\cite[Corollary 6.2.3 and Proposition 6.3.1]{MP}
with the duality theorem of Xia~\cite[Theorem 1.1(1)]{Xia}, we obtain
the desired vanishing.

For Theorem~\ref{thm:stability}, the preceding corollary supplies an initial metric
below the energy threshold for Calabi flow in the controlled cone. By the
degeneration theorem of Li, Wang, and Zheng~\cite[Proposition 4.4]{ControlledCompanion},
the cscK limit occurs in a test configuration with zero invariant.
K-polystability forces this configuration to be a product, so the original
class contains a cscK metric.

Starting instead with a cscK reference metric, we prove
Theorem~\ref{thm:sobolev-small} by fixing a Sobolev level above its
Sobolev constant. It follows from compactness and the sharp Euclidean Sobolev
inequality that no sequence at that level can have Calabi energy tending to zero.
By Theorem~\ref{thm:connected}, any metric of sufficiently small energy
can be joined to the reference metric within the same strict sublevel.
Continuity prevents the Sobolev constant from crossing the chosen level
along this path.

For Corollary~\ref{cor:calabi-global}, energy dissipation keeps the flow
in this sublevel. The Sobolev bound therefore holds throughout the flow,
and the extension criterion gives smooth existence for all time.\\

The paper is organized as follows. In Section~\ref{sec:resolvent}, we prove existence and establish estimates for
classical resolvents. The proof of path connectedness in
Section~\ref{sec:connectedness} uses these estimates.
Section~\ref{sec:controlled} identifies the smooth Calabi-energy infimum
in terms of the weak Mabuchi slope and gives applications to K-stability.
Finally, in Section~\ref{sec:sobolev-application}, we establish Sobolev
bounds and global Calabi flow for small initial energy in K\"ahler surface
classes containing a cscK metric.\\

After completing this paper, we learned that Hisamoto~\cite[Theorem 8.1]{HisamotoCalabi}
had also identified the smooth infimum of $\sqrt{\Ca}$ with the infimum of
the weak Mabuchi slope on every compact K\"ahler manifold, as in
Theorem~\ref{thm:smooth-weak-infima}. Hisamoto constructs smooth approximations
to the resolvent by solving a regularized equation. We construct classical
solutions of the resolvent equation and apply their variational inequality
to finite-energy competitors. For a smooth K-unstable polarized manifold,
Hisamoto~\cite[Theorem 1.1]{HisamotoCalabi} further expresses the smooth
infimum as the supremum of $-\Mab^{\mathrm{NA}}/\|\cdot\|_2$ over normal
ample algebraic test configurations of nonzero norm. For cohomologically
K-semistable real ample classes, including irrational classes, we prove in
Corollary~\ref{cor:semistable-small-energy} that the smooth infimum vanishes.\\

\textbf{On the usage of AI.}
The authors used ChatGPT as an assisting tool in the development of
the paper. The authors take full responsibility for the contents of
the paper.

\section{Classical resolvents}\label{sec:resolvent}
\subsection{The resolvent equation}
The classical resolvent is an implicit time step for the negative
Mabuchi gradient flow. Throughout this section, the center is a
smooth normalized potential $u_0\in\HH_0$, and $\tau>0$ is finite.
With the conventions of Section~\ref{sec:intro}, the class average is
\[
 \underline S=\int_XS_u\,d\mu_u
       =\frac nV\int_X\Ric(\omega_*)\wedge\omega_*^{n-1}.
\]
We normalize the Mabuchi functional by $\Mab(0)=0$. Its first variation is
\begin{equation}\label{eq:resolvent-mabuchi-variation}
 d\Mab_u(\xi)=-\int_X\xi(S_u-\underline S)\,d\mu_u.
\end{equation}
For the Mabuchi inner product
$\langle\xi,\eta\rangle_u=\int_X\xi\eta\,d\mu_u$, the gradient of $\Mab$ is $\underline S-S_u$, whose squared norm is $\Ca(u)$.
Thus its negative gradient flow is $\partial_tu=S_u-\underline S$.
The induced distance is denoted by $d_2$.

Darvas~\cite[Theorem 1]{DarvasCompletion} identifies the completion of $(\HH,d_2)$ with the CAT$(0)$ space $\EE^2$.
By a positive classical potential we mean $u\in C^{4,a}$ for some
$0<a<1$, with $\omega_u>0$. For two such endpoints, let
$\Gamma_{u,v}$ be their weak geodesic on $[0,1]$, with constant
speed. Its one-sided endpoint derivatives satisfy
\begin{equation}\label{eq:resolvent-endpoint-speed}
 d_2(u,v)^2=\int_X\dot\Gamma_{u,v}(0)^2\,d\mu_u
          =\int_X\dot\Gamma_{u,v}(1)^2\,d\mu_v.
\end{equation}
For smooth endpoints this is the speed formula of Darvas~\cite[Theorem 1]{Darvas}. For classical endpoints we prove it in Proposition~\ref{prop:classical-comparison}.
It follows from $d\AM_u(\xi)=\int_X\xi\,d\mu_u$ that
$|\AM(u)-\AM(v)|\le d_2(u,v)$. Hence $\AM$ extends continuously to $\EE^2$. Since
$\AM$ is affine along geodesics, the set $\EE^2_0=\AM^{-1}(0)$
is closed and convex. We use the geodesic convexity of $\Mab$ proved by
Berman and Berndtsson~\cite[Theorem 1.1]{BB} and the lower semicontinuous
extension constructed by Berman, Darvas, and Lu~\cite[Theorem 4.7]{BDLflow}.

For each positive classical potential $u$, define
\[
 b_{u_0}(u)=-\dot\Gamma_{u,u_0}(0).
\]
By Lemma~\ref{lem:resolvent-geodesic}, this function is Lipschitz in the spatial
variables. J.~Streets~\cite[Definition 4.4 and Lemma 4.9]{Streets}
introduced the resolvent equation, which in our notation reads
\begin{equation}\label{eq:resolvent-equation}
 b_{u_0}(u)=\tau(S_u-\underline S),\qquad
 u\in C^{4,a}(X,\mathbb R),\quad\omega_u>0,
 \quad 0<a<1.
\end{equation}
When the center is fixed, we write $b(u)$.
The resolvent equation is the Mabuchi analogue of the implicit Euler equation
$(z-z_0)/\tau=-\nabla F(z)$. Here the negative initial geodesic velocity
replaces the linear displacement $z-z_0$.
Put
\[
 D_{u_0}(v)=\tfrac12d_2(v,u_0)^2,\qquad
 \Phi_{\tau,u_0}(v)=\Mab(v)+\tau^{-1}D_{u_0}(v).
\]
For every smooth real direction $\xi$, it follows from the distance variation formula
of Streets~\cite[Lemma 3.13]{Streets} and
\eqref{eq:resolvent-mabuchi-variation} that
\[
 dD_{u_0}|_u(\xi)=\int_Xb_{u_0}(u)\xi\,d\mu_u,
 \qquad
 d\Phi_{\tau,u_0}|_u(\xi)
 =\int_X\left(\tau^{-1}b_{u_0}(u)-(S_u-\underline S)\right)
                                      \xi\,d\mu_u.
\]
In Proposition~\ref{prop:classical-comparison}, we extend the distance formula to
classical endpoints.
Thus \eqref{eq:resolvent-equation} is equivalent to the vanishing of
all smooth affine first variations of $\Phi_{\tau,u_0}$.
For the corresponding gradient flow in the finite-energy space, see
Berman, Darvas, and Lu~\cite[Section 6]{BDLflow}.

The variations include constant directions.
By \eqref{eq:classical-b-mean}, the equation itself imposes
$\AM(u)=0$. By \eqref{eq:resolvent-endpoint-speed}, it also gives
$d_2(u,u_0)^2=\tau^2\Ca(u)$.
In the following theorem we construct classical solutions for every finite positive
parameter and prove that they depend continuously on the data.

\begin{theorem}[Finite-parameter classical resolvents]\label{thm:resolvent}
For every smooth $u_0\in\HH_0$ and every finite $\tau>0$, equation
\eqref{eq:resolvent-equation} has a unique positive classical solution
$J_\tau u_0$. It belongs to $C^{4,\beta}$ for every $0<\beta<1$, and
\begin{equation}\label{eq:intro-energy}
 \Ca(J_\tau u_0)=\frac{d_2(J_\tau u_0,u_0)^2}{\tau^2}
                \le\Ca(u_0).
\end{equation}
For every compact family of smooth centers $\mathcal K\subset\HH_0$ and every interval
with $0<\tau_0\le T<\infty$, the solutions are uniformly positive and bounded in
$C^{4,\beta}$ on $\mathcal K\times[\tau_0,T]$.
The map $(\tau,u_0)\mapsto J_\tau u_0$ is continuous there into
$C^{4,\alpha}$ for every $0<\alpha<1$. For each smooth center,
$J_\tau u_0\to u_0$ in $C^4$ as $\tau\downarrow0$.
\end{theorem}

Here a compact smooth family is compact in $C^\infty$, with all
metrics positive. The metrics in such a family share a positive lower bound, and the
potentials are uniformly bounded in every finite-order norm. After establishing
endpoint regularity and a priori estimates, we construct a solution for small
parameters and continue it by degree. Compactness and uniqueness give continuity.

\subsection{Endpoint regularity and variational comparison}
\subsubsection{The endpoint term}

For every positive classical potential $u$, let $\Gamma_u$ be the
weak geodesic from $u$ to the fixed normalized potential $u_0$ and set
\begin{equation}\label{eq:resolvent-b}
 b(u)=-\partial_t\Gamma_u(0).
\end{equation}
The potential $u$ may have an arbitrary additive constant throughout the construction.
Since the Aubin--Mabuchi energy is affine along this geodesic, we have
\begin{equation}\label{eq:classical-b-mean}
 \int_Xb(u)\,d\mu_u=\AM(u)-\AM(u_0)=\AM(u).
\end{equation}
Indeed, $\frac{d}{dt}\AM(\Gamma_u(t))
=\int_X\dot\Gamma_u(t)\,d\mu_{\Gamma_u(t)}$. Its constant value is
$\AM(u_0)-\AM(u)$. The endpoint identity follows by smooth
approximation, or by taking one-sided difference quotients in the
mixed Monge--Amp\`ere formula for $\AM$.
Thus every solution of \eqref{eq:resolvent-equation} satisfies $\AM(u)=0$.
More explicitly, adding a constant $c$ to the initial endpoint adds
$(1-t)c$ to the whole geodesic, so
$b_{u_0}(u+c)=b_{u_0}(u)+c$, whereas $S_{u+c}=S_u$.
This transformation law fixes the additive constant. To use $b$ in the
degree argument we also need control of its spatial regularity and
its dependence on both endpoints.

\begin{lemma}\label{lem:resolvent-geodesic}
Fix $0<\beta<1$, $B<\infty$, and $\lambda>0$. If the two endpoints
$u,v$ satisfy $\|u\|_{C^{3,\beta}}+\|v\|_{C^{3,\beta}}\le B$ and
$\omega_u,\omega_v\ge\lambda\omega_*$, their weak geodesic has a
$C^{1,1}(X\times[0,1])$ bound depending only on these data and the
background. In particular, both endpoint velocities have a uniform
spatial Lipschitz bound. On uniformly positive, bounded $C^{4,\alpha}$ families, $b$ is
continuous with relatively compact image in $C^{0,\alpha}$ for $0<\alpha<1$. These
continuity statements
also hold when both endpoints vary with the same bounds.
\end{lemma}
\begin{proof}
The proof begins with uniform bounds for the approximate geodesics and then passes to
the endpoint velocities.

\paragraph{Step 1.}
Following Chen~\cite[Section 3]{Chen}, we work first with smooth endpoints on the cylinder
$\mathcal C=X\times\{w=t+\sqrt{-1}\theta:0\le t\le1\}$,
with $\theta$ periodic, and set
\[
 \Omega=\pi_X^*\omega_*+\sqrt{-1}\,dw\wedge d\bar w,
 \qquad \chi(t)=2t^2.
\]
The approximate geodesic is $U_\epsilon+\chi$, where
\begin{equation}\label{eq:resolvent-geodesic-approx}
 (\Omega+\ddc U_\epsilon)^{n+1}=\epsilon\Omega^{n+1},
 \qquad U_\epsilon|_{t=0}=u,\quad
 U_\epsilon|_{t=1}=v-2.
\end{equation}
To specify the uniformity in Chen's estimates, choose a strict
subsolution
\[
 \underline U=(1-t)u+tv-At(1-t)-\chi(t).
\]
We can choose $A$ depending only on $B,\lambda$ and the background so that
$\Omega+\ddc\underline U\ge c\Omega$ for a uniform constant $c>0$.
Indeed, the horizontal block is bounded below by
$\lambda\omega_*$, the mixed block is bounded, and the normal
block is $A/2$. Moreover, $\|\underline U\|_{C^{3,\beta}}\le C$.
For sufficiently small $\epsilon$, solvability of the smooth Dirichlet problems
\eqref{eq:resolvent-geodesic-approx} follows from the theorem of
B{\l}ocki~\cite[Theorem 1.3]{BlockiGeodesics}.
By uniqueness, $U_\epsilon$ is independent of $\theta$. Positivity in
the $w$ direction then implies that $U_\epsilon+\chi$ is convex in $t$.
Using the subsolution for the lower bound and convexity for the upper
bound, we obtain
\[
 \underline U\le U_\epsilon\le (1-t)u+tv-\chi(t).
\]
Both barriers agree with the boundary values and have uniformly bounded
first derivatives. They therefore bound $|U_\epsilon|$ on the cylinder
and $|\nabla U_\epsilon|$ on its boundary. Since the density $\epsilon$
is constant, the gradient estimate of
B{\l}ocki~\cite[Theorem 3.1]{BlockiGeodesics} gives a uniform bound
for $|\nabla U_\epsilon|$ throughout the cylinder.

For the mixed boundary derivatives, we use the barrier argument of
Chen~\cite[Section 3.1]{Chen}. Differentiating the equation once in a
tangential direction uses at most three derivatives of $\underline U$
in the barrier estimate. The gradient bound and the strict lower bound
for $\Omega+\ddc\underline U$ thus control the mixed derivatives.
At either endpoint the horizontal block is prescribed and bounded below
by $\lambda\omega_*$. Solving the determinant equation for the normal
entry now bounds $(U_\epsilon)_{w\bar w}$ there. Since $U_\epsilon$ is
independent of $\theta$, this also bounds $\partial_t^2U_\epsilon$.
The remaining tangential derivatives are fixed by the boundary data.
Applying the Laplacian estimate of
B{\l}ocki~\cite[Theorem 3.3]{BlockiGeodesics} with
$\Delta_\Omega\log\epsilon=0$, we conclude that
\begin{equation}\label{eq:resolvent-geodesic-estimates}
 \|U_\epsilon\|_{C^1(\overline{\mathcal C})}
 +\sup_{\overline{\mathcal C}}\Delta_\Omega U_\epsilon
 +\sup_{\partial\mathcal C}|\nabla^2U_\epsilon|\le C.
\end{equation}
Here the norms use the fixed product background, and $C$ depends only on
$B,\lambda,\beta$ and the background.

\paragraph{Step 2.}
Applying the real Hessian estimate of Chu, Tosatti, and Weinkove~\cite[Theorem
1.2]{CTW} with the bounds in \eqref{eq:resolvent-geodesic-estimates}, we obtain
\[
 \|U_\epsilon\|_{C^{1,1}(\overline{\mathcal C})}\le C.
\]
Its constant is independent of $\inf\log\epsilon$, and the
derivatives of the logarithmic density vanish. 
As $\epsilon\downarrow0$, the approximate geodesics converge to the unique weak
geodesic, which satisfies the same bound. For $C^{3,\beta}$ endpoints, choose smooth
approximations in $C^3$ with uniformly bounded $C^{3,\beta}$ norms and the lower metric
bound $(\lambda/2)\omega_*$. The comparison theorem of
Chen~\cite[Theorem 4]{Chen} identifies the uniform limit with
the weak geodesic. The uniform $C^{1,1}$ bound passes to the limit,
and convergence in $C^{1,\eta}$, $\eta<1$, identifies the endpoint
velocities. Restriction to the boundary therefore gives
\[
 \|\partial_t\Gamma_{u,v}(0)\|_{C^{0,1}(X)}
 +\|\partial_t\Gamma_{u,v}(1)\|_{C^{0,1}(X)}\le C.
\]

\paragraph{Step 3.}
For continuity, take two pairs of endpoints $(u,u_0)$ and
$(v,v_0)$ with these uniform bounds, and put
$\delta=\max\{\|u-v\|_\infty,\|u_0-v_0\|_\infty\}$.
Chen's comparison theorem gives
\[
 \|\Gamma_{u,u_0}-\Gamma_{v,v_0}\|_\infty\le\delta.
\]
The common $C^{1,1}$ bound controls the error in the initial
difference quotient by $Ch$. Thus, writing
$b_{u_0}(u)=-\partial_t\Gamma_{u,u_0}(0)$, we have
\[
 \|b_{u_0}(u)-b_{v_0}(v)\|_\infty
 \le\frac{2\delta}{h}+Ch,\qquad 0<h\le1.
\]
For $0<\delta\le1$, take $h=\sqrt\delta$. For $\delta>1$ the uniform velocity bound
suffices, and for $\delta=0$ the velocities agree by uniqueness. Interpolating with
the spatial Lipschitz bound, we obtain
\begin{equation}\label{eq:resolvent-b-continuity}
 \|b_{u_0}(u)-b_{v_0}(v)\|_{C^{0,\alpha}}
 \le C\delta^{(1-\alpha)/2}.
\end{equation}
The same conclusion holds at the terminal endpoint after reversing
$t$. It follows from the Arzel\`a--Ascoli theorem and the same interpolation that the uniformly Lipschitz velocities are relatively
compact in $C^{0,\alpha}$ for every $\alpha<1$.
Finally, a bounded family in $C^{4,\alpha}$ satisfies the required
$C^{3,\beta}$ bounds.
\end{proof}

\subsubsection{Classical variational comparison}

The endpoint estimates allow us to compare any classical solution with a
competitor. This comparison gives uniqueness and energy control before
we address existence.

\begin{proposition}\label{prop:classical-comparison}
A positive classical solution $u$ has $\AM(u)=0$. For every positive classical competitor $v$, we have
\begin{equation}\label{eq:classical-vi}
 \Mab(v)+\frac{d_2(v,u_0)^2}{2\tau}
 \ge\Mab(u)+\frac{d_2(u,u_0)^2+d_2(u,v)^2}{2\tau}.
\end{equation}
There is at most one classical solution at each parameter. Set $r=d_2(u,u_0)$ and
$E_0=\Ca(u_0)$. Then the energy and distance satisfy
\begin{align}
 r^2/\tau\le\Mab(u_0)-\Mab(u)&\le\sqrt{E_0}\,r,
       \label{eq:classical-distance-energy}\\
 \Ca(u)=r^2/\tau^2&\le E_0.
       \label{eq:classical-energy-identity}
\end{align}
\end{proposition}

\begin{proof}
Put $D(v)=\tfrac12d_2(v,u_0)^2$, let $\gamma=\Gamma_{u,v}$,
and write $w=\dot\gamma(0)$. The normalization follows from
\eqref{eq:classical-b-mean}.

\paragraph{Step 1.}
By geodesic convexity and the one-sided endpoint inequality of Berman and
Berndtsson~\cite[Lemma 3.5 and proof of Corollary 1.2]{BB}, we have
\begin{equation}\label{eq:classical-mabuchi-support}
 \Mab(v)-\Mab(u)\ge\int_X(\underline S-S_u)w\,d\mu_u.
\end{equation}
For classical endpoints, approximate both potentials in $C^4$ by smooth positive
potentials. The scalar curvatures, volume forms, and Mabuchi energies converge, and
Lemma~\ref{lem:resolvent-geodesic}
gives uniform convergence of the endpoint velocities.
The same approximation proves
\eqref{eq:resolvent-endpoint-speed} at this regularity.

For smooth potentials, the distance variation formula of Streets~\cite[Lemma
3.13]{Streets} gives
$dD_u(\xi)=\int_Xb(u)\xi\,d\mu_u$.
To extend the formula to a classical potential $u$, we approximate $u$ by smooth
potentials and pass to the limit in the integrated identity
\[
 D(u+t\xi)-D(u)
   =\int_0^t\!\int_Xb(u+s\xi)\xi\,d\mu_{u+s\xi}\,ds.
\]
The endpoint lemma gives uniform convergence on a fixed short
positive affine segment. For convergence of the distances, we use the estimate of
Darvas~\cite[Theorem 5.5]{Darvas}, namely
$d_2(p,q)\le C_n\|p-q\|_\infty$ for bounded
$\omega_*$-plurisubharmonic potentials.
At $u=u_0$ the formula follows directly from $D(u_0+t\xi)=O(t^2)$.

\paragraph{Step 2.}
To apply the first variation in the Lipschitz direction $w$,
choose smooth $\xi_\delta$ with
$\|\xi_\delta-w\|_\infty\le\delta$.
For fixed $\delta>0$, $u+s\xi_\delta$ is positive for small $s$.
Combining the $C^{1,1}$ geodesic bound with the endpoint distance comparison, we obtain
\[
 d_2(\gamma(s),u+s\xi_\delta)\le C(s\delta+s^2).
\]
Since $D$ is locally $d_2$-Lipschitz, its difference quotients
along these two paths differ by at most $C(\delta+s)$.
The limit as $s\downarrow0$ and then $\delta\downarrow0$ gives
\[
 \left.\frac d{ds}\right|_{0+}D(\gamma(s))
       =\int_Xb(u)w\,d\mu_u.
\]
We now use the CAT$(0)$ inequality
\[
 D(\gamma(s))\le(1-s)D(u)+sD(v)
                      -\tfrac12s(1-s)d_2(u,v)^2
\]
to bound the right derivative at zero. It follows that
\begin{equation}\label{eq:classical-distance-support}
 D(v)-D(u)\ge\int_Xb(u)w\,d\mu_u+\tfrac12d_2(u,v)^2.
\end{equation}
When we add $\tau^{-1}$ times this inequality to \eqref{eq:classical-mabuchi-support},
the linear terms cancel by \eqref{eq:resolvent-equation}. This proves
\eqref{eq:classical-vi}.

\paragraph{Step 3.}
For two solutions, apply \eqref{eq:classical-vi} in both directions and add. This gives $0\ge\tau^{-1}d_2(u,v)^2$, hence $u=v$.
With $u_0$ as the competitor, it follows that
$r^2/\tau\le\Mab(u_0)-\Mab(u)$.
Applying \eqref{eq:classical-mabuchi-support} at $u_0$ and using Cauchy--Schwarz
together with \eqref{eq:resolvent-endpoint-speed}, we obtain
\[
 \Mab(u_0)-\Mab(u)\le\sqrt{\Ca(u_0)}\,r.
\]
Thus $r\le\tau\sqrt{E_0}$, including when $r=0$.
Finally, the endpoint-speed identity and \eqref{eq:resolvent-equation} give
\[
 r^2=\int_Xb(u)^2\,d\mu_u=\tau^2\Ca(u),
\]
which proves the remaining assertions.
\end{proof}

\subsection{A priori estimates on finite parameter intervals}
\subsubsection{Entropy and a one-sided scalar bound}

The variational inequality controls entropy, whereas the first comparison estimate for
the geodesic velocity yields only an upper bound for scalar curvature.
The next lemma turns these bounds into a uniform potential estimate.

\begin{lemma}\label{lem:entropy-upper-scalar}
Fix the background K\"ahler manifold $(X^n,\omega_*)$, $n\ge1$. Suppose
$\phi\in C^{4,\alpha}(X)$ for some $0<\alpha<1$, $\omega_\phi>0$,
$\sup_X\phi=0$, and
\[
 F=\log\frac{\omega_\phi^n}{\omega_*^n},\qquad
 H(\phi)=\int_XF\,d\mu_\phi\le H_0,\qquad S_\phi\le K.
\]
Then there are $p>1$ and $C<\infty$, depending only on
$H_0,K,\omega_*$, such that
\[
 \int_Xe^{pF}\,d\mu_*\le C,\qquad \|\phi\|_\infty\le C.
\]
\end{lemma}

\begin{proof}
We first treat smooth $\phi$. The auxiliary equation of
Chen and Cheng~\cite[Section 5]{CC} supplies the weight in the argument.

\paragraph{Step 1.}
Set $\Phi(s)=\sqrt{1+s^2}$ and
\[
 A=\int_X\Phi(F)\,d\mu_\phi\le1+H_0+2/e.
\]
This follows from $(-F)e^F\le1/e$ on $\{F<0\}$ and
$\int_Xe^F\,d\mu_*=1$. Let $\psi$ solve
\[
 \omega_\psi^n=A^{-1}\Phi(F)\,\omega_\phi^n,
 \qquad \sup_X\psi=0.
\]
The right side has mass $V$, so a smooth solution exists by the theorem of Yau~\cite{Yau}.
Choose $C_R\ge0$ with $\Ric(\omega_*)\ge-C_R\omega_*$ and put
$L=2(1+C_R)$, $T=\tr_{\omega_\phi}\omega_*$.
For $0<a<1$, we combine the scalar equation with the arithmetic--geometric mean inequality to obtain
\begin{equation}\label{eq:entropy-test}
 \Delta_\phi(F+a\psi-L\phi)
 \ge (L-C_R-a)T+naA^{-1/n}\Phi(F)^{1/n}-K-nL.
\end{equation}

\paragraph{Step 2.}
Let $\alpha_*>0$ be the uniform exponential-integrability exponent for
supremum-normalized potentials in Chen and Cheng~\cite[Proposition 5.1]{CC}.
Choose $0<\delta\le1$ with $\delta L\le\alpha_*$ and then
$0<a<1$ with $a(2+\delta)\le\alpha_*$. Set $G=F+a\psi-L\phi$.
The trace coefficient in \eqref{eq:entropy-test} is positive. Since $A$ is bounded, the
right side of that inequality is at least $1$ on $\{|F|>B\}$ for a
uniform $B$ and at least $-C$ elsewhere. Integrating
$\Delta_\phi e^{\delta G}\ge\delta e^{\delta G}\Delta_\phi G$
against $d\mu_\phi$, we obtain
\[
 \int_X e^{\delta G}d\mu_\phi
 \le (1+C)\int_{\{|F|\le B\}} e^{\delta G}d\mu_\phi
 \le C'\int_X e^{-\delta L\phi}d\mu_*\le C''.
\]
Here $\psi\le0$ and $d\mu_\phi=e^Fd\mu_*$. Since $\phi\le0$, this
bounds $\int_X e^{(1+\delta)F+\delta a\psi}d\mu_*$.
Put $p=1+\delta/2$ and $\vartheta=p/(1+\delta)$. H\"older's inequality gives
\[
 \int_X e^{pF}d\mu_*
 \le \left(\int_X e^{(1+\delta)F+\delta a\psi}d\mu_*\right)^\vartheta
     \left(\int_X e^{-a(2+\delta)\psi}d\mu_*\right)^{1-\vartheta}
 \le C.
\]
The same exponential estimate for $\psi$ bounds the second factor.
The estimate of Ko{\l}odziej~\cite{Kolodziej} now bounds $\|\phi\|_\infty$.

For $\phi\in C^{4,\alpha}$, take smooth positive approximants
converging in $C^4$ and normalize their suprema to zero.
For all sufficiently large indices, their entropies are at most $H_0+1$ and their
scalar curvatures are bounded above by $K+1$. The smooth estimates are therefore uniform and
pass to the limit.
\end{proof}

\subsubsection{A priori estimates and compactness}
\label{sec:classical-apriori}

We apply the one-sided estimate to classical resolvents whose centers
and parameters vary in compact sets.

\begin{proposition}[Uniform estimates and compactness]
\label{prop:resolvent-apriori}
Fix $0<\tau_0\le T<\infty$ and a compact smooth family
$\mathcal K\subset\HH_0$. For every positive classical solution
$u$ of \eqref{eq:resolvent-equation} with $u_0\in\mathcal K$ and
$\tau_0\le\tau\le T$, put $F_u=\log(\omega_u^n/\omega_*^n)$.
The following assertions hold with uniform constants $C$ and
$C_\beta$.
\begin{enumerate}[label=\textup{(\roman*)},leftmargin=*]
\item The solution satisfies $\AM(u)=0$ and
\[
 C^{-1}\omega_*\le\omega_u\le C\omega_*.
\]

Moreover, for every $0<\beta<1$,
\[
 \|u\|_{C^{4,\beta}}+\|F_u\|_{C^{2,\beta}}\le C_\beta,
 \qquad
 \|b_{u_0}(u)\|_{C^{0,1}}+\|S_u\|_{C^{0,1}}\le C.
\]

\item Let $\mathscr Z_{\mathcal K,[\tau_0,T]}$ denote the set of
triples $(\tau,u_0,u)$ with $\tau\in[\tau_0,T]$,
$u_0\in\mathcal K$, and $u$ a positive classical solution of
\eqref{eq:resolvent-equation}. This set is compact in
$[\tau_0,T]\times\mathcal K\times C^{4,\alpha}(X)$
for every $0<\alpha<1$, with $\mathcal K$ carrying its $C^\infty$
topology. In particular, limits of such triples are again positive
classical solutions.
\end{enumerate}
The constants depend on the fixed background, $\mathcal K$, $\tau_0$,
$T$, and, for $C_\beta$, $\beta$. All norms use the fixed background.
The bounds may deteriorate as $T\to\infty$. No smallness assumption
on $\Ca(u_0)$, lower bound for $\Mab$, or cscK hypothesis is required.
\end{proposition}

The estimates bound all possible zeros in the degree construction
of Theorem~\ref{thm:resolvent}. Once existence is established, continuous dependence on
the parameters follows from compactness and uniqueness of classical solutions.
The proof first establishes the entropy and scalar bounds needed for the
coupled elliptic estimates.

\begin{proof}
We may assume $\mathcal K\ne\varnothing$.
Put $\Lambda=\tau_0^{-1}$ and choose uniform bounds
\[
 D_0=\sup_{v\in\mathcal K}d_2(v,0),\qquad
 E_0=\sup_{v\in\mathcal K}\Ca(v),\qquad
 M_0=\sup_{v\in\mathcal K}\Mab(v),\qquad
 R_T=D_0+T\sqrt{E_0}.
\]
The centers have uniform bounds in every smooth norm and a common positive lower metric
bound on $\mathcal K$.

\paragraph{Step 1.}
Proposition~\ref{prop:classical-comparison} gives $\AM(u)=0$ and
\begin{equation}\label{eq:classical-bounded-distance}
 d_2(u,u_0)\le T\sqrt{E_0},\qquad
 d_2(u,0)\le R_T,\qquad \Mab(u)\le M_0.
\end{equation}
The endpoint comparison of Darvas~\cite[Theorem 5.5]{Darvas} bounds
$\|u\|_{L^2(\mu_*)}$ by $C_nR_T$.
Since $\Delta_*u\ge-n$ and $\AM(u)=0$, it follows from the background Green-function
estimate that
\begin{equation}\label{eq:classical-sup}
 0\le m:=\sup_Xu\le C(1+R_T)=:B.
\end{equation}
For $q=u-m\le0$, the identity $\AM(q)=-m$ gives
\[
 \frac1V\sum_{j=0}^n\int_X(-q)\,
       \omega_u^j\wedge\omega_*^{n-j}=(n+1)m.
\]
Each normalized mixed measure has mass one. Since $|u|\le -q+m$,
its integral is controlled by the sum above. Hence
\[
 \frac1V\int_X|u|\,\omega_u^j\wedge\omega_*^{n-j}
 \le(n+2)m,\qquad 0\le j\le n.
\]
In our normalizations, the Chen--Tian formula is
\[
 \Mab(u)=H(u)+\underline S\AM(u)
 -\frac1V\sum_{j=0}^{n-1}\int_Xu\,\Ric(\omega_*)
                 \wedge\omega_u^j\wedge\omega_*^{n-1-j}.
\]
Bounding $\Ric(\omega_*)$ above and below by fixed multiples of $\omega_*$, we therefore obtain
\begin{equation}\label{eq:classical-entropy}
 H(u)\le M_0+CB\le\max\{M_0,0\}+CB=:H_0.
\end{equation}

\paragraph{Step 2.}
Set $K_u=\sup_X(u-u_0)$. The function $u-tK_u$ is an
admissible subgeodesic with terminal value at most $u_0$.
Comparison gives $\Gamma_{u,u_0}(t)\ge u-tK_u$, so
$b(u)\le K_u\le B+\|u_0\|_\infty$. Thus
\begin{equation}\label{eq:classical-scalar-upper}
 S_u\le\underline S+\Lambda(B+\|u_0\|_\infty).
\end{equation}
Lemma~\ref{lem:entropy-upper-scalar}, applied to $u-m$, now
bounds $\|u\|_\infty$. Set $A_u=\|u-u_0\|_\infty$. By subgeodesic comparison and
pointwise convexity in $t$, we have
\[
 u-tA_u\le\Gamma_{u,u_0}(t)\le(1-t)u+tu_0\le u+tA_u.
\]
It follows by differentiation at the initial endpoint that $\|b(u)\|_\infty\le A_u$.
The resolvent equation \eqref{eq:resolvent-equation} then gives
\begin{equation}\label{eq:classical-scalar-full}
 \|u\|_\infty+\|S_u\|_\infty\le C,\qquad
 \|S_u\|_\infty\le|\underline S|+\Lambda\|u-u_0\|_\infty.
\end{equation}

\paragraph{Step 3.}
Put $F=\log(\omega_u^n/\omega_*^n)$. The coupled equations are
\[
 \omega_u^n=e^F\omega_*^n,\qquad
 \Delta_uF=-S_u+\tr_{\omega_u}\Ric(\omega_*).
\]
Apply the estimate of Chen and Cheng~\cite[Theorem 5.3]{CC} with potential $u-m$,
$f=S_u$, and $\eta=\Ric(\omega_*)$.
In view of the entropy bound \eqref{eq:classical-entropy} and the scalar bound
\eqref{eq:classical-scalar-full}, for every $1<p<\infty$ we obtain
\begin{equation}\label{eq:classical-cc-sobolev}
 \|u-m\|_{W^{4,p}(\omega_*)}+\|F\|_{W^{2,p}(\omega_*)}\le C_p.
\end{equation}
For all sufficiently large $p$, the estimates also follow from
Lu and Seyyedali~\cite[Theorem 1.3]{LuSeyyedali}. The weighted entropy is bounded
by $1+H_0+2/e$, and all scalar $L^p(\mu_u)$ norms are bounded.
These theorems allow a variable scalar term and require no
bounds for its derivatives.

For a classical $C^{4,a}$ solution, apply the smooth estimates
to positive approximants converging in $C^4$, with their suprema
subtracted. The entropies and scalar suprema converge, so the
estimates have uniform constants and pass to the limit.
The approximants need only satisfy the coupled equations with
their own scalar curvatures.

Taking $p>2n/(1-\beta)$ in \eqref{eq:classical-cc-sobolev} and
using Sobolev embedding and $m\le B$, we obtain
\begin{equation}\label{eq:resolvent-c3}
 \|u\|_{C^{3,\beta}}+\|F\|_{C^{1,\beta}}\le C_\beta,
 \qquad 0<\beta<1.
\end{equation}
The $C^2$ bound gives an upper bound for the metric eigenvalues,
and $\det(g_u)/\det(g_*)=e^F\ge e^{-C}$ bounds them below. Thus
\begin{equation}\label{eq:classical-c2}
 C^{-1}\omega_*\le\omega_u\le C\omega_*.
\end{equation}
Lemma~\ref{lem:resolvent-geodesic}, with the endpoint exponent
fixed at $1/2$, now gives
\[
 \|b(u)\|_{C^{0,1}}+\|S_u\|_{C^{0,1}}\le C.
\]
The scalar equation is uniformly elliptic, with uniformly
controlled coefficients and a $C^{0,\beta}$ right-hand side.
Applying Schauder estimates first to the equation for $F$ and then to the differentiated
Monge--Amp\`ere equation for $u$, we obtain
\begin{equation}\label{eq:resolvent-c4}
 \|F\|_{C^{2,\beta}}+\|u\|_{C^{4,\beta}}\le C_\beta,
 \qquad 0<\beta<1.
\end{equation}
See also the elliptic argument of Chen and Cheng~\cite[Proposition 4.2]{CC}.
This proves \textup{(i)} in every dimension.

\paragraph{Step 4.}
Fix $0<\alpha<\beta<1$. Consider any sequence of solution triples. By compactness of
the parameter interval and the family of centers, together with the uniform
$C^{4,\beta}$ bound, we may pass to a subsequence for which
\[
 \tau_j\to\tau\in[\tau_0,T],\qquad
 u_{0,j}\to u_0\in\mathcal K\quad\text{in }C^\infty,
 \qquad u_j\to u\quad\text{in }C^{4,\alpha}.
\]
The metric lower bound preserves positivity, and
$S_{u_j}\to S_u$ in $C^{0,\alpha}$.
Lemma~\ref{lem:resolvent-geodesic}, with both endpoints varying,
gives $b_{u_{0,j}}(u_j)\to b_{u_0}(u)$ in the same space.
Passing to the limit in \eqref{eq:resolvent-equation}, we see that the limiting triple
is again a solution. Since the spaces in the statement are metrizable, this proves
\textup{(ii)}.
\end{proof}

\subsection{Classical existence by degree}\label{sec:degree}
\subsubsection{A classical solution and local degree near the initial metric}
\label{sec:classical-seed}

We construct a small-parameter solution near the prescribed center,
then use the finite-parameter estimates to continue it by degree. 
Throughout the local construction and its continuation, fix
\[
 0<\alpha_0<1/3,\qquad
 0<\kappa<\frac{1-3\alpha_0}{1+\alpha_0},\qquad
 \mathscr X=C^{4,\alpha_0}(X,\mathbb R),\quad
 \mathscr Y=C^{0,\alpha_0}(X,\mathbb R).
\]
All norms refer to the fixed background. Both spaces contain constants.
In particular, we do not impose the nonlinear normalization $\AM=0$
on a Banach space. Instead, the equation itself will enforce it.

\begin{proposition}[Small-parameter existence and local degree]
\label{prop:resolvent-seed}
For every smooth $u_0\in\HH_0$, with the exponents $\alpha_0,\kappa$
fixed above, there are $M>0$ and $\tau_*>0$ such that the following assertions hold for
every $0<\tau\le\tau_*$.
\begin{enumerate}[label=\textup{(\roman*)},leftmargin=*]
\item Equation \eqref{eq:resolvent-equation} has a unique positive
classical solution $u_\tau$. It satisfies $\AM(u_\tau)=0$ and
\begin{equation}\label{eq:resolvent-zero-endpoint}
 \|u_\tau-u_0\|_0\le M\tau,\qquad
 \|u_\tau-u_0\|_{4,\alpha_0}\le\tau^\kappa.
\end{equation}
In particular, $u_\tau\to u_0$ in $C^{4,\alpha_0}$ as
$\tau\downarrow0$.

\item Define the bounded open convex set
\[
 U_\tau=\{w\in\mathscr X:\|w\|_0<M\tau,
                         \ \|w\|_{4,\alpha_0}<\tau^\kappa\}.
\]
The potentials in $u_0+\overline U_\tau$ have a common positive
metric lower bound. The map $\mathcal T_\tau$ constructed in
\eqref{eq:initial-compact-map} is continuous and compact from
$\overline U_\tau$ into $U_\tau$, and
\[
 w=\mathcal T_\tau(w)
 \quad\Longleftrightarrow\quad
 b_{u_0}(u_0+w)=\tau(S_{u_0+w}-\underline S).
\]
Its Leray--Schauder degree is
\begin{equation}\label{eq:initial-local-degree}
 \deg_{\rm LS}(I-\mathcal T_\tau,U_\tau,0)=1.
\end{equation}
\end{enumerate}
The constants depend only on $u_0$, the fixed background, and the
chosen exponents.
\end{proposition}

\begin{proof}
We construct a compact self-map near the prescribed center.

\paragraph{Step 1.}
Write $u=u_0+w$, $\widetilde b(w)=b_{u_0}(u_0+w)$.
On a fixed small $C^{4,\alpha_0}$ ball, the affine metrics
$\omega_{u_0}+(1-t)\ddc w$ have a common positivity margin.
By the Schur complement criterion, we obtain a subgeodesic of the form
\[
 u_0+(1-t)w-Ct(1-t)\|\nabla w\|_\infty^2
\]
for a sufficiently large uniform constant $C$. Comparing the weak geodesic with this lower barrier and the affine
upper barrier, we obtain
\begin{equation}\label{eq:initial-geodesic-remainder}
 \widetilde b(w)=w+e(w),\qquad
 0\le e(w)\le C\|\nabla w\|_\infty^2.
\end{equation}
The endpoint lemma bounds $e$ in $C^{0,1}$ and makes it continuous
and compact into $\mathscr Y$. Interpolation gives
\begin{equation}\label{eq:initial-holder-remainder}
 \|e(w)\|_{0,\alpha_0}
       \le C\|\nabla w\|_\infty^{2(1-\alpha_0)}.
\end{equation}

\paragraph{Step 2.}
Put $f(w)=S_{u_0+w}-\underline S$, $f_0=f(0)$.
The scalar-curvature decomposition of
Chen and He~\cite[equation (3.1)]{ChenHeStability} gives
\begin{equation}\label{eq:initial-quasilinear}
 f(w)=f_0-L_ww+N(w),\qquad
 \|N(w)\|_{0,\alpha_0}\le C\|w\|_{4,\alpha_0}^2.
\end{equation}
To specify the operator, let $\mathscr P_w$ be the positive
fourth-order principal part of $-f$, using a fixed connection, and set
\[
 L_w=\mathscr P_w+(L_0-\mathscr P_0),\qquad
 L_0v=\Delta_{u_0}^2v+
          \langle\Ric(\omega_{u_0}),\ddc v\rangle_{u_0}.
\]
The principal coefficients depend only on $D^2w$.
Thus $N$ involves at most three derivatives, and its value and
first differential vanish at zero, proving the quadratic estimate.
We retain the varying principal part in $L_w$ to obtain compactness below.

\paragraph{Step 3.}
For a sufficiently small fixed ball $\|w\|_{4,\alpha_0}\le r_0$
and $0<\tau\le\tau_{\rm loc}$, the inverse
$R_{\tau,w}=(I+\tau L_w)^{-1}:\mathscr Y\to\mathscr X$ exists, with
\begin{equation}\label{eq:initial-parameter-estimates}
 \|R_{\tau,w}q\|_0\le C\|q\|_0,\qquad
 \tau\|R_{\tau,w}q\|_{4,\alpha_0}\le C\|q\|_{0,\alpha_0}.
\end{equation}
For the second bound, we apply the sectorial estimate of Lunardi~\cite[Section
3.2]{Lunardi} as in Chen and He~\cite[Sections 3.1--3.2]{ChenHeStability}. This gives
the large-parameter Schauder estimate for the smooth operator $L_0$.
To pass from little H\"older spaces to the stated H\"older spaces, approximate the
right side smoothly in $C^{0,\gamma}$, $\gamma<\alpha_0$, with bounded
$C^{0,\alpha_0}$ norms. The estimate passes to the limit by compactness in
$C^{4,\gamma}$ and lower semicontinuity of the highest-order seminorm.
A homogeneous solution is smooth, so uniqueness also transfers.
Since
\[
 \|\tau R_{\tau,0}(L_w-L_0)\|_{\mathscr X\to\mathscr X}
       \le C\|w\|_{4,\alpha_0},
\]
a Neumann series gives the inverse and its Schauder bound
uniformly on the closed ball, after choosing $Cr_0<1/2$.

For the $C^0$ estimate, suppose instead that
$(I+\tau_jL_{w_j})z_j=q_j$ with
$\|z_j\|_\infty=1$ and $\|q_j\|_\infty\to0$.
If $\tau_j\to0$, choose a point where $|z_j|=1$ and rescale lengths by
$\tau_j^{-1/4}$ in background normal coordinates.
For $p>2n$, apply interior $W^{4,p}$ estimates on each fixed rescaled ball. After
passing to a diagonal subsequence, we obtain a bounded entire limit with $|z(0)|=1$
solving $(I+P_\infty)z=0$ on $\mathbb R^{2n}$.
Here $P_\infty$ has a constant positive quartic symbol because the coefficients
converge locally to constants and all lower-order terms vanish under rescaling.
Taking the Fourier transform in the sense of tempered distributions, we obtain
$(1+p_\infty(\xi))\widehat z=0$, so $z=0$, a contradiction.
If $\tau_j\to\tau_\infty>0$, pass to a subsequence such that $w_j\to w_\infty$ in $C^{4,\gamma}$.
Global $W^{4,p}$ estimates and Schauder regularity give a nonzero
limit in the kernel of $I+\tau_\infty L_{w_\infty}$, contradicting invertibility on the
closed ball. This proves the first estimate. For the higher-order elliptic estimates,
see
Agmon, Douglis, and Nirenberg~\cite{AgmonDouglisNirenberg}.
Finally, since $f_0$ is smooth, it follows from $R_{\tau,w}f_0=f_0-\tau R_{\tau,w}L_wf_0$ that
\begin{equation}\label{eq:initial-smooth-source}
 \|R_{\tau,w}(\tau f_0)\|_{4,\alpha_0}\le C\tau.
\end{equation}

\paragraph{Step 4.}
The resolvent equation \eqref{eq:resolvent-equation} is equivalent to
\begin{equation}\label{eq:initial-compact-map}
 w=\mathcal T_\tau(w)
   :=R_{\tau,w}\{\tau f_0+\tau N(w)-e(w)\}.
\end{equation}
On $\overline U_\tau$, interpolation gives
$\|\nabla w\|_\infty\le C_M\tau^{(3+\kappa)/4}$.
The preceding estimates imply
\begin{align*}
 \|\mathcal T_\tau(w)\|_0
 &\le C\tau\|f_0\|_0+C\tau^{1+2\kappa}
                              +C_M\tau^{(3+\kappa)/2},\\
 \|\mathcal T_\tau(w)\|_{4,\alpha_0}
 &\le C\tau+C\tau^{2\kappa}+C_M\tau^\eta,
 \qquad \eta=\frac{(3+\kappa)(1-\alpha_0)}2-1.
\end{align*}
Choose $M>2C\|f_0\|_0+1$ and then take $\tau$ small.
The error exponents in the first line exceed one, while
$1,2\kappa,\eta>\kappa$, since
\[
 \eta-\kappa
   =\tfrac12\{1-3\alpha_0-\kappa(1+\alpha_0)\}>0.
\]
Thus $\mathcal T_\tau(\overline U_\tau)\subset U_\tau$.

At fixed $\tau>0$, compactness follows from the compact embedding
$\mathscr X\hookrightarrow C^{3,\alpha_0}$, the three-jet
dependence of $N$, and compactness of $e$ in $\mathscr Y$.
The principal coefficients are precompact in $C^{0,\alpha_0}$. By the inverse identity
\[
 R_{\tau,w}-R_{\tau,z}
       =\tau R_{\tau,w}(L_z-L_w)R_{\tau,z}
\]
and the uniform bound on the inverses, the solutions converge in $\mathscr X$ whenever
the sources and coefficients converge.
The same identity proves continuity.

It remains to compute the degree of this self-map. Since $U_\tau$ is convex, contains
zero, and contains the image of its closure, $t\mathcal T_\tau$, $0\le t\le1$, is an
admissible compact homotopy to zero. Normalization of Leray--Schauder degree gives
\eqref{eq:initial-local-degree}, hence a fixed point by
Deimling~\cite[Sections 8.3--8.4]{Deimling}.
Proposition~\ref{prop:classical-comparison} supplies normalization
and uniqueness among all positive classical solutions.
Membership in $U_\tau$ gives \eqref{eq:resolvent-zero-endpoint}
and convergence at zero. Choose $\tau_*$ to satisfy all the
preceding smallness conditions.
\end{proof}

\subsubsection{Continuing solutions of the resolvent equation by degree}
\label{sec:resolvent-degree}

We first show that the global compact map has the same degree as the local map for a
small parameter. The degree then remains constant as the parameter varies, by
Proposition~\ref{prop:resolvent-apriori}.

\begin{proposition}[Classical existence for all finite parameters]
\label{prop:resolvent-existence}
For every smooth center $u_0\in\HH_0$ and every finite $\tau>0$,
equation \eqref{eq:resolvent-equation} has a unique positive
classical solution $u_\tau$. It satisfies $\AM(u_\tau)=0$ and
belongs to $C^{4,\beta}$ for every $0<\beta<1$.
For $0<\tau\le\tau_*$ it is the solution of
Proposition~\ref{prop:resolvent-seed}, so in particular it satisfies
\eqref{eq:resolvent-zero-endpoint}.
\end{proposition}

\begin{proof}
It suffices to continue the degree from a small positive parameter to an arbitrary finite upper bound.

\paragraph{Step 1.}
Choose an initial parameter $\tau_s>0$ small enough for
Proposition~\ref{prop:resolvent-seed} and the coercivity estimate
in Step 3 below. Fix $T>\tau_s$.
Proposition~\ref{prop:resolvent-apriori} bounds all possible zeros
on $[\tau_s,T]$, so we can choose $R,\delta>0$ such that they
lie strictly inside the bounded open set
\begin{equation}\label{eq:resolvent-degree-domain}
 \mathcal O=\{u\in\mathscr X:\|u\|_{4,\alpha_0}<R,
                                  \ \omega_u>\delta\omega_*\},
\end{equation}
with $u_0+\overline U_{\tau_s}\subset\mathcal O$.
All metrics corresponding to potentials in $\overline{\mathcal O}$ have common upper
and lower bounds. Define
\begin{equation}\label{eq:resolvent-full-residual}
 G_\tau(u)=b(u)-\tau(S_u-\underline S),\qquad
 \mathscr F_\tau(u)=\rho_uG_\tau(u),\qquad
 \rho_u=\frac{d\mu_u}{d\mu_*}.
\end{equation}
Every zero is automatically normalized by
\eqref{eq:classical-b-mean}.

\paragraph{Step 2.}
Let $a_u^{abcd}\nabla^*_{abcd}$ be the principal part of
$\rho_u(\underline S-S_u)$, using the fixed background connection, and set
\begin{equation}\label{eq:resolvent-preconditioner}
 A_uv=a_u^{abcd}\nabla^*_{abcd}v+cv.
\end{equation}
Since the quartic symbol is positive and the coefficients are uniformly bounded in
$C^{2,\alpha_0}$, we have the G\aa rding inequality
\[
 \int_XvA_uv\,d\mu_*
 \ge a_0\|v\|_{H^2}^2+(c-C_0)\|v\|_{L^2}^2.
\]
Choose $c>C_0$ uniformly on $\overline{\mathcal O}$.
The operators $(1-s)(\Delta_*^2+I)+sA_u$ remain coercive and
strongly elliptic. Combining the Schauder estimate with interpolation and coercivity, we obtain
\[
 \|v\|_{4,\alpha_0}
 \le C\bigl(\|B_{s,u}v\|_{0,\alpha_0}+\|v\|_{L^2}\bigr)
 \le C'\|B_{s,u}v\|_{0,\alpha_0}.
\]
Here $B_{s,u}=(1-s)(\Delta_*^2+I)+sA_u$.
It follows by the method of continuity, starting at $s=0$, that
$A_u:\mathscr X\to\mathscr Y$ is an isomorphism with a uniform bound on its inverse. We use the higher-order estimates of
Agmon, Douglis, and Nirenberg~\cite{AgmonDouglisNirenberg}.

Define
\begin{equation}\label{eq:resolvent-compact-map}
 T_\tau(u)=u-(\tau A_u)^{-1}\mathscr F_\tau(u)
          =(\tau A_u)^{-1}\{\tau A_uu-\mathscr F_\tau(u)\}.
\end{equation}
The fourth derivatives cancel in braces, leaving at most three
derivatives of $u$ and the compact endpoint term. The coefficients
of $A_u$ are precompact in $C^{0,\alpha_0}$.
As in the local construction, the inverse identity and uniform
Schauder bound show that $T_\tau$ is jointly continuous and
compact on $\overline{\mathcal O}\times[\tau_s,T]$.
Its fixed points are exactly the classical solutions.

\paragraph{Step 3.}
For $u=u_0+w$, $w\in\overline U_{\tau_s}$, put
\[
 P_w=\rho_u(I+\tau_sL_w),\qquad
 Q_{t,w}=(1-t)P_w+t\tau_sA_u,\quad 0\le t\le1.
\]
On this local domain, the G\aa rding estimate gives
\[
 \int_XvP_wv\,d\mu_*
 \ge a\tau_s\|v\|_{H^2}^2+(m-C\tau_s)\|v\|_{L^2}^2
\]
with $a,m>0$ independent of small $\tau_s$.
Choose $\tau_s$ so that $m-C\tau_s>0$.
Both $P_w$ and $\tau_sA_u$ are then coercive.
Starting with $P_w^{-1}=R_{\tau_s,w}\rho_u^{-1}$ and arguing as in Step 2, we find
that every $Q_{t,w}$ is invertible. For the fixed $\tau_s$, the inverses are uniformly
bounded.

Each $Q_{t,w}$ has the same fourth-order coefficient as $\mathscr F_{\tau_s}$. Define
\[
 K_t(w)=w-Q_{t,w}^{-1}\mathscr F_{\tau_s}(u_0+w).
\]
The same cancellation argument shows that this is a continuous compact homotopy.
At $t=0$ it is $\mathcal T_{\tau_s}$, and at $t=1$ it is
$T_{\tau_s}(u_0+w)-u_0$.
All its fixed points solve the equation at the initial parameter and therefore lie
strictly inside $U_{\tau_s}$. By homotopy invariance and translation of the domain, we
obtain
\begin{equation}\label{eq:resolvent-local-degree}
 \deg_{\rm LS}(I-T_{\tau_s},u_0+U_{\tau_s},0)=1.
\end{equation}

\paragraph{Step 4.}
By uniqueness, every solution at the initial parameter lies in $u_0+U_{\tau_s}$.
Excision therefore identifies the degree on $\mathcal O$ with
\eqref{eq:resolvent-local-degree}.
The a priori bounds exclude boundary fixed points throughout
$[\tau_s,T]$. Homotopy invariance therefore gives
\begin{equation}\label{eq:resolvent-global-degree}
 \deg_{\rm LS}(I-T_\tau,\mathcal O,0)=1,
                  \qquad\tau_s\le\tau\le T.
\end{equation}
Since the degree is nonzero, a solution exists throughout this interval. The local
construction covers $0<\tau\le\tau_s$.
Since $T$ is arbitrary, existence holds for every finite
parameter. Normalization and uniqueness follow from
Proposition~\ref{prop:classical-comparison}, and regularity follows from Proposition~\ref{prop:resolvent-apriori}.
\end{proof}

\subsection{Strong continuity and completion of the resolvent theorem}\label{sec:resolvent-continuity}
\begin{proof}[Proof of Theorem~\ref{thm:resolvent}]
Existence and uniqueness follow from Proposition~\ref{prop:resolvent-existence}. The
energy estimate \eqref{eq:intro-energy} follows from
Proposition~\ref{prop:classical-comparison}, and the uniform bounds on finite
parameter intervals follow from Proposition~\ref{prop:resolvent-apriori}\textup{(i)}.
The small-parameter estimate \eqref{eq:resolvent-zero-endpoint} gives
convergence in $C^4$ as $\tau\downarrow0$.

For joint continuity, let $(\tau_j,u_{0,j})\to(\tau,u_0)$ with
$\tau>0$ and centers in a compact smooth family.
By Proposition~\ref{prop:resolvent-apriori}\textup{(ii)}, the solution triples are
precompact in the product topology, with the $C^{4,\alpha}$ topology on the solution
component for every $\alpha<1$. Every limiting triple solves the equation with data
$(\tau,u_0)$. Uniqueness identifies its last component as
$J_\tau u_0$. Thus the full sequence converges, which proves continuity.
\end{proof}

\section{Connectedness of Calabi-energy sublevels}\label{sec:connectedness}
To join two low-energy endpoints, we start with a smooth path of centers.
The energy along this path may exceed the prescribed threshold. The
following estimate lets us bring the whole path below it with a single
resolvent parameter.

\begin{lemma}\label{lem:contraction}
For smooth normalized centers $x,y$ and every $\tau>0$, we have
\begin{align}
 d_2(J_\tau x,J_\tau y)&\le d_2(x,y),\label{eq:contraction}\\
 \sqrt{\Ca(J_\tau x)}
 &\le\sqrt{\Ca(y)}+\frac{2d_2(x,y)}{\tau}.
 \label{eq:two-center}
\end{align}
\end{lemma}

\begin{proof}
By Theorem~\ref{thm:resolvent}, the resolvent equation
\eqref{eq:resolvent-equation} has positive classical solutions
$p=J_\tau x$ and $q=J_\tau y$ at the respective centers.
Write $d=d_2$. Taking each solution as a competitor for the other in
\eqref{eq:classical-vi}, we obtain
\begin{align*}
 2\tau(\Mab(q)-\Mab(p))
 &\ge d(p,x)^2+d(p,q)^2-d(q,x)^2,\\
 2\tau(\Mab(p)-\Mab(q))
 &\ge d(q,y)^2+d(p,q)^2-d(p,y)^2.
\end{align*}
When the two inequalities are added, the Mabuchi terms cancel. This gives the first inequality below.
\[
 2d(p,q)^2
 \le d(q,x)^2+d(p,y)^2-d(p,x)^2-d(q,y)^2
 \le 2d(p,q)d(x,y).
\]
The second inequality is the quadrilateral Cauchy--Schwarz inequality
of Berg and Nikolaev~\cite{BN}, applied in the CAT$(0)$ space
$(\EE^2,d)$ of Darvas~\cite[Theorem~1]{DarvasCompletion}.
This proves \eqref{eq:contraction}, including the case $p=q$.

It follows from the triangle inequality and \eqref{eq:contraction} that
$d(p,x)\le d(p,q)+d(q,y)+d(y,x)\le d(q,y)+2d(x,y)$.
Dividing by $\tau$ and using \eqref{eq:intro-energy} at both centers, we
obtain \eqref{eq:two-center}, since
$d(q,y)/\tau=\sqrt{\Ca(q)}\le\sqrt{\Ca(y)}$.
\end{proof}

\begin{proof}[Proof of Theorem~\ref{thm:connected}]
Fix $0<\alpha<1$ and $a,b\in\mathcal U_c^\alpha$, so $c>0$.
Write $N(u)=u-\AM(u)$. The mixed Monge--Amp\`ere formula gives
$\AM(u+k)=\AM(u)+k$ for each constant $k$, since each mixed
volume has mass $V$. Thus $N(u)$ is normalized and has the same
metric and Calabi energy as $u$. The same formula gives continuity
of $\AM$ in $C^2$, so $N$ preserves
continuity of paths in $C^4$, $C^{4,\alpha}$, and $C^\infty$.

\paragraph{Step 1.}
Suppose first that $a,b$ are smooth. The path
$\eta(s)=N((1-s)a+sb)$, $0\le s\le1$, joins them through smooth
positive potentials, since its metrics are convex combinations of
$\omega_a$ and $\omega_b$. Let $D$ be its Mabuchi length, which is
finite by smoothness. Then $d_2(\eta(s),a)\le D$ for all $s$.
Choose a finite $T>0$ with
$2D/T<\sqrt c-\sqrt{\Ca(a)}$. Applying \eqref{eq:two-center} with centers $\eta(s)$ and $a$, we obtain
\[
 \sup_{0\le s\le1}\sqrt{\Ca(J_T\eta(s))}
 \le \sqrt{\Ca(a)}+2D/T<\sqrt c.
\]
The centers $\eta([0,1])$ form a compact family in $C^\infty$.
It follows from Theorem~\ref{thm:resolvent} that $s\mapsto J_T\eta(s)$ is a
$C^{4,\alpha}$-continuous path from $J_Ta$ to $J_Tb$.
For $u=a,b$, the same theorem gives continuity of $\tau\mapsto J_\tau u$
in $C^{4,\alpha}$ on $(0,T]$ and in $C^4$ at $\tau=0$ upon setting $J_0u=u$.
Moreover, for $\tau>0$,
$\Ca(J_\tau u)\le\Ca(u)<c$ by \eqref{eq:intro-energy}.
Following the path from $a$ to $J_Ta$, then $J_T\eta$, and finally the return path
from $J_Tb$ to $b$, we obtain a $C^4$-continuous path
$\gamma:[0,1]\to\mathcal U_c^\alpha$.

\paragraph{Step 2.}
It remains in the smooth case to replace $\gamma$ by a polygonal path with smooth
vertices and the same endpoints. In $C^4(X,\mathbb R)$, positivity is open because
$u\mapsto\omega_u$ is continuous in $C^2$ and $X$ is compact.
On this open set, convergence of the potentials in $C^4$ implies convergence of the
scalar curvatures and volume densities in $C^0$. It follows from \eqref{eq:energy-definition}
that $O_c=\{u:\omega_u>0,\ \Ca(u)<c\}$ is open.

By compactness of the image of $\gamma$, choose $r>0$ so that
$B_{C^4}(\gamma(t),2r)\subset O_c$ for all $t$. By uniform continuity, we can choose a
partition $0=t_0<\cdots<t_m=1$
so that $\gamma$ has $C^4$ oscillation less than $r$ on each
subinterval. At each partition point $t_i$, choose a smooth $v_i$
within $r$ of $\gamma(t_i)$ in $C^4$ by spatial mollification,
keeping $v_0=a$ and $v_m=b$. Both $v_i$ and $v_{i+1}$ lie in
$B_{C^4}(\gamma(t_i),2r)$ by the triangle inequality. Since the ball is convex, the
whole segment $w_i(s)=(1-s)v_i+sv_{i+1}$,
$0\le s\le1$, lies in $O_c$. Each segment is continuous in every $C^k$
norm because its vertices are smooth. The normalized segments
$N(w_i(s))$ agree at shared vertices and retain the endpoints
$N(a)=a$ and $N(b)=b$. They give the required $C^\infty$-continuous path.

\paragraph{Step 3.}
For arbitrary $a,b\in\mathcal U_c^\alpha$, choose a $C^4$
ball about each endpoint contained in $O_c$ and a smooth approximant
in each ball. The affine segment from each endpoint to its approximant
lies in that ball and is continuous in $C^{4,\alpha}$, since its
endpoints belong to $C^{4,\alpha}$. Normalize the two segments and
join their smooth ends by Steps 1 and 2. This middle path is also
$C^{4,\alpha}$-continuous by the continuous inclusion
$C^5(X)\hookrightarrow C^{4,\alpha}(X)$. The assertion follows by concatenating the
three paths.
\end{proof}

\section{Calabi-energy infima and stability}
\label{sec:controlled}
\label{sec:semistability-zero-infimum}
We first identify the smooth Calabi-energy infimum with the weak
Mabuchi-slope infimum in an arbitrary compact K\"ahler class.
Entropy recovery in real ample classes then gives a zero infimum
under K-semistability. For classes in the controlled cone of a compact
K\"ahler surface,
we deduce cscK existence under K-polystability from the class-preserving degeneration
theorem of Li, Wang, and Zheng~\cite{ControlledCompanion}.

Let $\Mab$ also denote its lower semicontinuous extension to
$\EE^2$, and write $\operatorname{Dom}(\Mab)$ for its domain of
finite Mabuchi energy. For $w\in\operatorname{Dom}(\Mab)$, its
descending metric slope is
\[
 |\partial\Mab|(w)=\sup_{v\in\operatorname{Dom}(\Mab),\ v\ne w}
 \frac{(\Mab(w)-\Mab(v))_+}{d_2(w,v)}.
\]
This agrees with the local slope by convexity, as in
Xia~\cite[equation (3.2)]{Xia}. Indeed, along the geodesic from
$w$ to $v$, convexity gives
$\Mab(w)-\Mab(\gamma(t))\ge t(\Mab(w)-\Mab(v))$ and
$d_2(w,\gamma(t))=t\,d_2(w,v)$. Taking positive parts and letting $t\downarrow0$, we
see that the local slope bounds each global quotient.
The reverse inequality follows from the definitions.
By Xia~\cite[Proposition 3.6]{Xia}, at a smooth
potential $w$ this slope equals $\sqrt{\Ca(w)}$. Put
\[
 s_*:=\inf_{w\in\EE^2}|\partial\Mab|(w).
\]
The slope is assigned the value $+\infty$ outside the domain of $\Mab$.
In particular $0\le s_*<\infty$, since any smooth potential has
finite Mabuchi energy and finite slope.

At smooth potentials the slope identity gives one inequality between the two infima.
The reverse inequality follows by comparing the energy of a classical resolvent with
the slope at a finite-energy competitor.

\begin{proposition}
\label{prop:weak-classical-slope}
For every smooth normalized center $x$, every $\tau>0$, and every $w\in\EE^2$ with
finite Mabuchi energy and finite slope, we have
\begin{equation}\label{eq:weak-center-estimate}
 \sqrt{\Ca(J_\tau x)}
 \le |\partial\Mab|(w)+\frac{2d_2(x,w)}{\tau}.
\end{equation}
Consequently
\begin{equation}\label{eq:smooth-weak-infimum}
 \inf_{u\in\HH}\sqrt{\Ca(u)}=s_*.
\end{equation}
In particular, if $s_*=0$, then the fixed class contains smooth positive metrics with arbitrarily small Calabi energy.
\end{proposition}
\begin{proof}
Fix $x$ and $\tau$ as in the statement, and let $p=J_\tau x$ be the
positive classical solution given by Theorem~\ref{thm:resolvent}.

\paragraph{Step 1.}
We first extend the classical variational inequality to weak competitors.
Let $w\in\EE^2$ with $\Mab(w)<\infty$. To justify the
approximation, write the extended Chen--Tian formula as
\[
 \Mab(z)=H_*(z)+\mathcal E_*(z),\qquad
 H_*(z)=\int_X\log\left(\frac{d\mu_z}{d\mu_*}\right)d\mu_z,
\]
where $H_*(z)=+\infty$ if $\mu_z$ is not absolutely continuous
with respect to $\mu_*$ or the entropy integral is infinite.
The nonentropy term $\mathcal E_*$ is a fixed
linear combination of the Aubin--Mabuchi energy and the energy
contracted with the smooth form $\Ric(\omega_*)$, plus a
normalizing constant. These terms are finite and $d_1$-continuous
on $\EE^1$ by Berman, Darvas, and Lu~\cite[Propositions 4.1 and 4.4]{BDLflow}.
Since $\EE^2\subset\EE^1$, it follows that $H_*(w)$ is finite whenever $\Mab(w)$ is
finite.

By the approximation theorem of Berman, Darvas, and Lu~\cite[Theorem 3.2]{BDLflow},
applied with exponent $2$ and the untwisted background measure, there exist
$w_j\in\HH$ such that
\[
 d_2(w_j,w)\longrightarrow0,\qquad H_*(w_j)\longrightarrow H_*(w).
\]
Our probability normalization gives $d_1\le d_2$, first on
smooth paths by H\"older's inequality and then on the completions.
Hence $\mathcal E_*(w_j)\to\mathcal E_*(w)$ and
$\Mab(w_j)\to\Mab(w)$, as in
Berman, Darvas, and Lu~\cite[Theorem 4.7]{BDLflow}.

Applying \eqref{eq:classical-vi} at center $x$ with competitor $w_j$
and passing to the limit using $d_2(w_j,w)\to0$, we obtain
\begin{equation}\label{eq:weak-competitor-vi}
 d_2(p,x)^2+d_2(p,w)^2-d_2(w,x)^2
 \le2\tau(\Mab(w)-\Mab(p)).
\end{equation}

\paragraph{Step 2.}
Write $s=|\partial\Mab|(w)$, $a=d_2(p,x)$,
$r=d_2(p,w)$, and $D=d_2(x,w)$. The global slope inequality gives
$\Mab(w)-\Mab(p)\le sr$, including $r=0$. Combining this with
\eqref{eq:weak-competitor-vi}, we obtain
\[
 a^2+r^2-D^2\le2\tau sr.
\]
The reverse triangle inequality gives $a\ge|r-D|$, so
\[
 2r(r-D)=(r-D)^2+r^2-D^2
 \le a^2+r^2-D^2\le2\tau sr.
\]
If $r>0$, division by $2r$ gives $r\le D+\tau s$.
For $r=0$ the same bound holds directly.
Thus $a\le r+D\le\tau s+2D$. By
\eqref{eq:intro-energy},
$\sqrt{\Ca(p)}=a/\tau$, proving
\eqref{eq:weak-center-estimate}.

\paragraph{Step 3.}
Put $I=\inf_{u\in\HH}\sqrt{\Ca(u)}$. Choose smooth real potentials
$p_k\to p$ in $C^4$. By compactness of $X$ and strict positivity, we have $\omega_p\ge c_p\omega_*$ for some $c_p>0$, and hence
$\omega_{p_k}\ge(c_p/2)\omega_*$ for all sufficiently large $k$.
Thus $p_k\in\HH$ for all sufficiently large $k$. The scalar curvature is a continuous
expression in the derivatives of the potential up to order four
and the inverse metric on this positive set. Consequently
\[
 S_{p_k}\longrightarrow S_p,\qquad
 \frac{d\mu_{p_k}}{d\mu_*}\longrightarrow
 \frac{d\mu_p}{d\mu_*}\quad\hbox{uniformly},\qquad
 \Ca(p_k)\longrightarrow\Ca(p).
\]
The average scalar curvature is fixed by the class. Thus, for every finite
$\tau>0$,
\[
 I\le\sqrt{\Ca(J_\tau x)}
 \le |\partial\Mab|(w)+\frac{2d_2(x,w)}{\tau}.
\]
Fix $w$ before letting $\tau\to\infty$. This numerical inequality
implies
\[
 \inf_{u\in\HH}\sqrt{\Ca(u)}\le |\partial\Mab|(w).
\]
Taking the infimum over the finite-slope points, we obtain $I\le s_*$.
Points with infinite slope do not affect this infimum, since
$s_*<\infty$. Conversely, for every $u\in\HH$,
\[
 s_*\le|\partial\Mab|(u)=\sqrt{\Ca(u)}
\]
by Xia~\cite[Proposition 3.6]{Xia}, in the normalization fixed above.
Taking the infimum over smooth $u$, we obtain $s_*\le I$, proving
\eqref{eq:smooth-weak-infimum}.
If $s_*=0$, choose $u\in\HH$ with $\sqrt{\Ca(u)}<\sqrt\varepsilon$.
This proves the final assertion and completes the proof of
Theorem~\ref{thm:smooth-weak-infima}.
\end{proof}

For real ample classes on a projective manifold, we use the
non-Archimedean pluripotential theory of Boucksom and Jonsson~\cite{BJDiv}, Mesquita-Piccione~\cite{MP}, and Mesquita-Piccione and Witt Nystr\"om~\cite{MPN}.

With the non-Archimedean Monge--Amp\`ere measures normalized to have mass one, our
convention is
\[
 H^{\mathrm{NA}}_\Omega(\phi)
 =\int_{X^{\mathrm{NA}}}A_X(v)\,
              \operatorname{MA}_\Omega(\phi)(v).
\]
For a divisorial valuation $v=c\operatorname{ord}_F$ on a smooth
birational model $Y\to X$, the log discrepancy is
$A_X(v)=c(1+\operatorname{ord}_F K_{Y/X})$. The trivial valuation
has discrepancy zero. The space $\mathcal E^{1,\mathrm{NA}}(\Omega)$
consists of finite-energy non-Archimedean psh potentials.
Its strong topology is the refinement of the weak topology by
continuity of the Monge--Amp\`ere energy, with additive constants
fixed when applying the Monge--Amp\`ere homeomorphism. See
Boucksom and Jonsson~\cite{BJDiv} and Mesquita-Piccione and Witt Nystr\"om~\cite[Sections 5--6]{MPN}.

To apply semistability, we need model potentials in the original real class
whose entropies converge to that of the given potential. The proof uses
the local discrepancy estimate in Boucksom's appendix to Trusiani~\cite{Trusiani}.

\begin{lemma}
\label{lem:real-entropy-recovery}
Let $X$ be smooth and projective and let
$\Omega\in\operatorname{NS}(X)_{\mathbb R}$ be ample, with no
rationality assumption. Suppose that $\phi\in\mathcal E^{1,\mathrm{NA}}(\Omega)$ has
finite non-Archimedean entropy. Then there are model potentials $\phi_j$ in the same
class, represented by relatively K\"ahler cohomological test configurations, such that
$\phi_j\to\phi$
strongly and
\[
 H^{\mathrm{NA}}_\Omega(\phi_j)\longrightarrow
 H^{\mathrm{NA}}_\Omega(\phi).
\]
Consequently the same approximations recover
$\Mab^{\mathrm{NA}}_\Omega$.
\end{lemma}
\begin{proof}
We extend the base-ideal argument in Boucksom's appendix to
Trusiani~\cite{Trusiani}. Neither the discrepancy estimate nor the graded-ideal lemma
requires a polarization. It remains to check that, after rounding a real ample class, the resulting family is graded and admits a fixed auxiliary divisor that ensures global generation. Write $n=\dim X$ and
$V=\Omega^n$, and normalize all Monge--Amp\`ere measures to have
mass one. All models and resolutions below are projective and
equivariant for the $\mathbb C^*$-action. Relative base ideals are
invariant, so their blowups and equivariant resolutions remain
models of test configurations.

\paragraph{Step 1.}
Write $P_\Omega(f)$ for the largest $\Omega$-psh potential below
a model function $f$. It suffices to consider $\phi=P_\Omega(f)$ with $f$ rational, as
the following approximation shows.
Applying the approximation theorem of Boucksom and Jonsson~\cite[Theorem 3.10]{BJDiv}
to $\mu=\operatorname{MA}_\Omega(\phi)$, we obtain divisorial probability measures
that satisfy $\mu_j\to\mu$ strongly and
$\int A_X\,d\mu_j\to\int A_X\,d\mu$. By the
Monge--Amp\`ere homeomorphism of
Mesquita-Piccione and Witt Nystr\"om~\cite[Theorem B]{MPN}, their
sup-normalized potentials converge strongly to the sup-normalization
of $\phi$. Constants can be restored at the end. Mesquita-Piccione and Witt
Nystr\"om~\cite[Proposition 8.3.2]{MPN} show that each such potential is the envelope
of a real model function on an snc model.

On this fixed model, approximate the coefficients of its vertical
real divisor by rational coefficients to obtain model functions $f_k$ that converge
uniformly to $f$. The envelope operator is nonexpansive in the uniform norm. Indeed,
if $f-\varepsilon\le f_k\le f+\varepsilon$, then $P_\Omega(f)-\varepsilon\le
P_\Omega(f_k) \le P_\Omega(f)+\varepsilon$. Thus their envelopes converge
uniformly, hence strongly. By the restricted-volume formula of
Mesquita-Piccione and Witt Nystr\"om~\cite[Proposition 4.2.2]{MPN}, after a common constant shift if
needed, all their Monge--Amp\`ere measures are supported on the
same finite set of central-fibre vertices. By weak convergence, the masses converge at
each vertex. Since $A_X$ is finite on this finite set, the entropies converge as well.
The reduction follows by taking a diagonal sequence with a consistent normalization.

\paragraph{Step 2.}
Choose very ample line bundles $A_1,\ldots,A_r$ and positive real
numbers $a_i$ such that
\[
 \Omega=\sum_{i=1}^r a_i c_1(A_i).
\]
For this, choose a rational simplex in the ample cone whose interior contains $\Omega$. Let
$b:\mathcal X\to X\times\mathbb P^1$ be a smooth dominating snc
model on which $f$ is represented by a vertical $\mathbb Q$-divisor
$D$. Write $P_i=b^*p_X^*A_i$ and
$\mathcal L=\sum_i a_iP_i+D$ as a real divisor class. By rescaling
the class and the model function by the same positive integer, we
may assume that $D$ is integral. This scaling will be undone at the end.
We work over a neighborhood of $0\in\mathbb P^1$ and denote the structure map by $\pi$.
Adding an integral multiple of the central fibre shifts the model
function by a constant. It leaves its Monge--Amp\`ere measure,
entropy, and Mabuchi invariant unchanged. Choose the shift so that
$D\ge\mathcal X_0$. On the projective compactification,
$\sum_i a_iP_i+\mathcal X_0$ is the pullback of an ample class on
$X\times\mathbb P^1$. Adding the effective divisor $D-\mathcal X_0$
shows that $\mathcal L$ is big.

After passing to a higher fixed snc model, we can choose an effective vertical $\mathbb Q$-divisor $E$, supported on the relative augmented
base locus, such that $\mathcal L-E$ is relatively ample. To see this for a real
class, proceed as follows. Choose a rational class $\mathcal L'$
arbitrarily close to $\mathcal L$, with $\mathcal L-\mathcal L'$
relatively ample and with the same relative augmented base locus.
The latter equality follows from
Ein, Lazarsfeld, Musta\c t\u a, Nakamaye, and Popa~\cite[Proposition 1.5 and Corollary 1.6]{ELMNPBase}. Their proof by
finitely many Kodaira decompositions applies over a neighborhood
of $0$. Such rational classes exist because rational points are
dense in the finite-dimensional N\"eron--Severi space and the
relative ample cone is open. Approximate a small ample
subtraction from $\mathcal L$ closely enough that the difference
remains ample. Apply the rational construction in the proof of
Boucksom~\cite[Appendix A, Theorem A.6]{Trusiani} to $\mathcal L'$. On its higher
model, $\mathcal L'-E$ is relatively ample and $E$ is supported
on its augmented base locus. The pullback of
$\mathcal L-\mathcal L'$ is relatively nef, so
$\mathcal L-E$ is relatively ample as well. The birational
transformation rule for augmented base loci preserves their equality.
We now use this higher model as $\mathcal X$. Write
$\mathcal X_0=\sum_G b_GG$ and
$E=\sum_G e_GG$, with $e_G\ge0$, and let $v_G$ be the
valuation on $X$ induced by $b_G^{-1}\operatorname{ord}_G$.
By the restricted-volume formula of Mesquita-Piccione and Witt
Nystr\"om~\cite[Proposition 4.2.2]{MPN}, we have
\[
 \operatorname{MA}_\Omega(\phi)
 =\frac1V\sum_G b_G
   \operatorname{vol}_{\mathcal X\mid G}(\mathcal L)\,
   \delta_{v_G}.
\]
If $e_G>0$, then $G$ lies in the relative augmented base locus.
Every decomposition of $\mathcal L$ into an ample class and an effective
real divisor on the compactification restricts to a relative Kodaira
decomposition near $0$. Thus $G$ lies in the absolute augmented base
locus as well. This locus is proper, so the prime divisor $G$ is one
of its irreducible components. By Ein, Lazarsfeld, Musta\c t\u a,
Nakamaye, and Popa~\cite[Theorem 5.7 and Remark 5.8]{ELMNPVol},
the restricted volume of $\mathcal L$ along $G$ vanishes.
Since $v_G(E)=e_G/b_G$, we obtain
\begin{equation}\label{eq:real-envelope-contact}
 \int e\,\operatorname{MA}_\Omega(\phi)=0,
 \qquad e(v):=v(E)\ge0.
\end{equation}
The function $e$ is a continuous model function. 

\paragraph{Step 3.}
For each positive integer $m$, put
\begin{equation}\label{eq:real-rounded-bundles}
 \mathcal L_m=\sum_i\lfloor ma_i\rfloor P_i+mD,
 \qquad R_m=m\mathcal L-\mathcal L_m
       =\sum_i\{ma_i\}P_i.
\end{equation}
Let $\mathfrak a_m$ be the relative base ideal of $\mathcal L_m$.
These are vertical ideals, since their restrictions to every fibre
away from $0$ are base-point-free. Locally near $0$, each
$\mathfrak a_m$ therefore contains a power of the base parameter
$t$, so it defines a finite continuous model function. Consider the difference bundle
\[
 \mathcal L_{m+k}-\mathcal L_m-\mathcal L_k
 =\sum_i\bigl(\lfloor(m+k)a_i\rfloor-
       \lfloor ma_i\rfloor-\lfloor ka_i\rfloor\bigr)P_i.
\]
Every coefficient on the right is either zero or one, so this bundle is globally generated.
By multiplying sections and then using local nonvanishing sections of this difference bundle, we obtain
\begin{equation}\label{eq:real-rounded-grading}
 \mathfrak a_m\mathfrak a_k\subset\mathfrak a_{m+k}.
\end{equation}

On the normalized blowup $\rho_m$ of $\mathfrak a_m$, write
$\mathfrak a_m\mathcal O=\mathcal O(-F_m)$ and let
\begin{equation}\label{eq:real-rounded-potentials}
 \psi_m=f+m^{-1}\log|\mathfrak a_m|.
\end{equation}
The corresponding model class over $\Omega$ is
\[
 \rho_m^*\mathcal L-m^{-1}F_m
 =m^{-1}(\rho_m^*\mathcal L_m-F_m)+m^{-1}\rho_m^*R_m,
\]
which is relatively nef. The first summand is relatively generated.
The second is a nonnegative combination of pullbacks of the $A_i$.
Thus $\psi_m$ is $\Omega$-psh and
$\psi_m\le P_\Omega(f)\le f$.

To prove uniform convergence, first fix rational numbers $0<b_i<a_i$.
Set $B=\sum_i b_i c_1(A_i)$ and $\mathcal L_B=\sum_i b_iP_i+D$.
Choose $d$ so that all $db_i$ are integral. For $d\mid m$, denote by
$\mathfrak c_m$ the relative base ideal of $m\mathcal L_B$.
Over the trivially valued field $\mathbb C$, Boucksom and
Jonsson~\cite[Theorem 5.20 and Lemma 5.22]{BJGlobal} prove continuity
of the envelope and identify it as the uniform limit of the base-ideal
approximants. Applying their result to $d\mathcal L_B$ and dividing by
$d$, we obtain
\[
 f+m^{-1}\log|\mathfrak c_m|\longrightarrow P_B(f)
 \quad\text{uniformly as }m\to\infty,\quad d\mid m.
\]
For these degrees, $\mathcal L_m-m\mathcal L_B$ is a nonnegative
integral combination of the globally generated bundles $P_i$.
Multiplying sections gives $\mathfrak c_m\subset\mathfrak a_m$,
so each approximant on the left is bounded above by $\psi_m$.

It remains to compare the two envelopes. Put
$\eta=\max_i(1-b_i/a_i)$, so that
$B-(1-\eta)\Omega$ is a nonnegative combination of ample
classes. If $C=\|f\|_\infty$, then $-C\le P_\Omega(f)\le C$
and
\[
 (1-\eta)P_\Omega(f)-\eta C
       \le P_B(f)\le P_\Omega(f).
\]
The left side is $B$-psh and at most $f$, which proves the first
inequality. The second follows because $\Omega-B$ is nef. Thus
$\|P_B(f)-P_\Omega(f)\|_\infty\le2\eta C$.
Given $\varepsilon>0$, choose the $b_i$ so that $2\eta C<\varepsilon$.
Since $d\mid k!$ for all large $k$, the rational approximation and the
upper bound $\psi_{k!}\le P_\Omega(f)$ give
$0\le P_\Omega(f)-\psi_{k!}\le 2\varepsilon$ for all sufficiently
large $k$. Thus $\psi_{k!}\to P_\Omega(f)$ uniformly and hence strongly.

\paragraph{Step 4.}
Set $\mathfrak b_m=\mathcal J(\mathfrak a_\bullet^m)$.
Asymptotic multiplier ideals satisfy
$\mathfrak a_m\subset\mathcal J(\mathfrak a_m)
\subset\mathfrak b_m$ and are antigraded, that is,
$\mathfrak b_{m+k}\subset\mathfrak b_m\mathfrak b_k$.
See Trusiani~\cite[Section 2.3]{Trusiani}.
There is a single sufficiently divisible integer $m_0$, independent
of $m$, such that
\begin{equation}\label{eq:real-fixed-buffer}
 \mathfrak b_m\otimes
 \mathcal O_{\mathcal X}(\mathcal L_{m+m_0}-m_0E)
 \quad\text{is relatively globally generated.}
\end{equation}
We verify this using ordinary multiplier ideals of integral bundles.
Fix a relatively very ample line bundle $H$ and put $d=\dim\mathcal X$.
Every $R_j$ lies in the fixed compact box
$\{\sum_i t_iP_i:0\le t_i\le1\}$. Since $\mathcal L-E$ is
relatively ample, choose one integer $m_0$, divisible enough that
$m_0E$ is integral, for which
\begin{equation}\label{eq:real-buffer-positivity}
 m_0(\mathcal L-E)-R_j-K_{\mathcal X/\mathbb P^1}-qH
 \quad\text{is relatively ample}
 \quad(j\ge1,\ 1\le q\le d).
\end{equation}
For each fixed $m$, choose a sufficiently divisible $k$ realizing
$\mathfrak b_m=\mathcal J(\mathfrak a_k^{m/k})$. No uniform
choice of $k$ in $m$ is required. The evaluation map for the
integral bundle $\mathcal L_k$ generates
$\mathcal O_{\mathcal X}(\mathcal L_k)\otimes\mathfrak a_k$
relatively. If
$N_m=\mathcal L_{m+m_0}-m_0E$, then
\[
 N_m-qH-K_{\mathcal X/\mathbb P^1}-\frac m k\mathcal L_k
 =m_0(\mathcal L-E)-R_{m+m_0}
       -K_{\mathcal X/\mathbb P^1}-qH+\frac m k R_k.
\]
The right side is relatively ample by
\eqref{eq:real-buffer-positivity}, because $R_k$ is relatively nef.
Applying relative Nadel vanishing to $\mathcal J(\mathfrak a_k^{m/k})$, we obtain
$R^i\pi_*(\mathfrak b_m\otimes\mathcal O(N_m-qH))=0$
for every $i>0$ and $1\le q\le d$.
Relative Castelnuovo--Mumford regularity proves
\eqref{eq:real-fixed-buffer}. See
Lazarsfeld~\cite[Vol.~I, Section 1.8, and Vol.~II, Sections 9.4 and 11.2]{Lazarsfeld}.
Here relative Nadel vanishing is applied on a log resolution of
$\mathfrak a_k$ and descends by local vanishing for multiplier ideals.
Multiplying by the section of $m_0E$, we obtain
\[
 \mathfrak b_m(-m_0E)\subset\mathfrak a_{m+m_0}.
\]
To apply the lemma of Boucksom~\cite[Appendix A, Lemma A.7]{Trusiani}, use
\eqref{eq:real-rounded-grading}, antigrading,
$\mathfrak a_m\subset\mathfrak b_m$, and the last inclusion,
with fixed divisor $m_0E$.
For $m\ge2m_0$, the lemma yields
\begin{equation}\label{eq:real-ideal-defect}
 0\le v(\mathfrak a_m)-v(\mathfrak b_m)\le2m_0e(v).
\end{equation}
The valuation inequalities can also be seen directly. Fix a Gauss
valuation $v$ and write $A_m=v(\mathfrak a_m)$,
$B_m=v(\mathfrak b_m)$, and $F=m_0e(v)$. Grading gives the
nonnegative limit $a=\lim_m A_m/m=\inf_m A_m/m$.
By applying antigrading and $\mathfrak a_m\subset\mathfrak b_m$ at multiples of $m$
and passing to the limit, we obtain $B_m\le ma$. The preceding ideal inclusion gives
$A_{m+m_0}\le B_m+F$. Consequently, for $m\ge2m_0$,
\[
 A_m\le B_{m-m_0}+F\le(m-m_0)a+F,\qquad
 B_m\ge A_{m+m_0}-F\ge(m+m_0)a-F.
\]
Subtracting these inequalities and using $a\ge0$, we obtain
$0\le A_m-B_m\le2F$, uniformly in $m$ and $v$.

\paragraph{Step 5.}
Let $A_{\mathcal X}$ be the log discrepancy of the log smooth pair
$(\mathcal X,\mathcal X_{0,\mathrm{red}})$, evaluated at Gauss
extensions of valuations. For a central component $G$ on the
normalized blowup, put $b_G=\operatorname{ord}_G(t)$ and denote
by $v_G$ the restriction of $b_G^{-1}\operatorname{ord}_G$ to
$\mathbb C(X)$. Values on ideals and divisors on $\mathcal X$
are taken using its Gauss extension. To relate ordinary discrepancy to log discrepancy,
we use the following identity. If
$a_G=\operatorname{ord}_G(K_{\mathcal X_m/\mathcal X})$, then
\[
 A_{\mathcal X}(v_G)
 =\frac{a_G+1-
       \operatorname{ord}_G(\rho_m^*\mathcal X_{0,\mathrm{red}})}
       {b_G}
 \le\frac{a_G}{b_G},
\]
because a central component has
$\operatorname{ord}_G(\rho_m^*\mathcal X_{0,\mathrm{red}})\ge1$.
Nonnegativity follows from the snc condition. At an exceptional component, divide the local Rees-valuation estimate of Boucksom~\cite[Proposition A.1]{Trusiani} by $b_G$. Using $\mathcal J(\mathfrak a_m)\subset\mathfrak b_m$ and writing $v=v_G$, we obtain
\[
 0\le A_{\mathcal X}(v)
 \le(\dim\mathcal X)
       \bigl(v(\mathfrak a_m)-v(\mathfrak b_m)\bigr)
 \le2m_0(\dim\mathcal X)e(v).
\]
At strict transforms of the original central components,
$A_{\mathcal X}=0$. Write
$\mathcal B_m=\rho_m^*\mathcal L-m^{-1}F_m$. Its
Monge--Amp\`ere measure is
\[
 \operatorname{MA}_\Omega(\psi_m)
 =\frac1V\sum_{G\subset\mathcal X_{m,0}}
    b_G(\mathcal B_m^n\cdot G)\,\delta_{v_G}.
\]
In particular, addition of the nef remainder $R_m/m$ changes
the masses but introduces no valuations outside this finite set.
Integrating the preceding estimate, we obtain
\begin{equation}\label{eq:real-entropy-defect}
 0\le\int A_{\mathcal X}\operatorname{MA}_\Omega(\psi_m)
 \le2m_0(\dim\mathcal X)
       \int e\operatorname{MA}_\Omega(\psi_m)\longrightarrow0.
\end{equation}
This limit is taken along the factorial degrees of Step 3 and follows from uniform convergence of $\psi_m$, weak continuity of Monge--Amp\`ere measures, and \eqref{eq:real-envelope-contact}.

To recover the entropy, take $v=v_G$ and denote its Gauss extension by
$\widetilde v=b_G^{-1}\operatorname{ord}_G$, so that $\widetilde v(t)=1$.
As in Mesquita-Piccione and Witt Nystr\"om~\cite[proof of Proposition 9.1.3]{MPN},
we have $A_X(v)=A_{X\times\mathbb P^1}(\widetilde v)-1$.
To express the right side in terms of the pair
$(\mathcal X,\mathcal X_{0,\mathrm{red}})$, we use the relative
canonical divisor formula and $\widetilde v(\mathcal X_0)=1$. This gives
\[
 A_X(v)=A_{\mathcal X}(v)+\chi_{\mathcal X}(v),\qquad
 \chi_{\mathcal X}(v)
 =\widetilde v\bigl(K_{\mathcal X/(X\times\mathbb P^1)}
             +\mathcal X_{0,\mathrm{red}}-\mathcal X_0\bigr).
\]
The function $\chi_{\mathcal X}$ is continuous, so its integrals converge
under weak convergence of the measures. Since $\operatorname{MA}_\Omega(\phi)$ is supported
on the original central-fibre vertices, where $A_{\mathcal X}=0$,
\eqref{eq:real-entropy-defect} proves entropy convergence.

\paragraph{Step 6.}
For each $m$, choose a fixed equivariant snc resolution
$\sigma_m:\mathcal Y_m\to\mathcal X_m$, obtained by projective
modifications over $0$, and put
$\widetilde{\mathcal B}_m=\sigma_m^*\mathcal B_m$.
There is a relatively K\"ahler class $\mathcal C_m$ on
$\mathcal Y_m$ restricting to $\Omega$. Start with the pullback
of $\Omega$ from $X\times\mathbb P^1$ and add sufficiently
small vertical exceptional classes that are relatively ample
for the successive modifications, as in
Sj\"ostr\"om Dyrefelt~\cite[proof of Proposition 3.12]{SDSemi}. These perturbations
are supported over $0$, so their restrictions to every noncentral
fibre vanish. Thus the general-fibre class remains exactly
$\Omega$. For $0<\epsilon<1$, set
\[
 \mathcal B_{m,\epsilon}
 =(1-\epsilon)\widetilde{\mathcal B}_m
       +\epsilon\mathcal C_m.
\]
This is relatively K\"ahler and still restricts to $\Omega$.
On this fixed model, its potential converges uniformly to
$\psi_m$ as $\epsilon\downarrow0$, since the difference is
$\epsilon$ times a fixed continuous model function.

To compute its entropy, we sum the intersection masses of $\mathcal B_{m,\epsilon}$
with the fixed weights $A_X(v_G)$. There are finitely many masses, each polynomial in
$\epsilon$. A component introduced
by $\sigma_m$ has zero mass at $\epsilon=0$. Its pushforward
as an $n$-cycle is zero, so
$(\widetilde{\mathcal B}_m^n\cdot G)=0$ by the projection
formula. Thus entropy converges to that of $\psi_m$.
Choose $\epsilon_m$ so that both the uniform potential error
and the entropy error are less than $1/m$ along the factorial
subsequence. No discrepancy bound uniform in $m$ is needed for
this choice. This proves the lemma for rational model envelopes.

Undo the simultaneous scaling of the class and potential. Indeed,
for $c>0$ the normalized measure satisfies
$\operatorname{MA}_{c\Omega}(c\psi)=\operatorname{MA}_\Omega(\psi)$,
as both the numerator intersection masses and the volume scale
by $c^n$. Entropy is therefore unchanged, and the inverse scaling
preserves strong convergence. Then take the diagonal sequence
from Step 1, so that the potentials converge in the strong topology and their
entropies converge as well. Restore the
original additive normalization. All final model classes restrict
to the original $\Omega$. Finally, write
$\Mab^{\mathrm{NA}}_\Omega=H^{\mathrm{NA}}_\Omega+
\mathcal R^{\mathrm{NA}}_\Omega$, where the nonentropy energy
term $\mathcal R^{\mathrm{NA}}_\Omega$ is finite and strongly
continuous by the real-class energy-pairing theory of
Boucksom and Jonsson~\cite{BJDiv} and Mesquita-Piccione~\cite{MP}. It follows from entropy convergence that $\Mab^{\mathrm{NA}}_\Omega$ also converges.
\end{proof}

\begin{proof}[Proof of Corollary~\ref{cor:semistable-small-energy}]
\begin{samepage}
The first part of the argument takes place in the class $\Omega$.

\paragraph{Step 1.}
Let $\Mab^{\mathrm{NA}}_\Omega$ denote the non-Archimedean Mabuchi
functional, normalized compatibly with $\Mab$. We first prove its
nonnegativity on relatively K\"ahler model metrics in the exact
class $\Omega$. If $(\mathcal X,\mathcal A)$ represents such a
metric $\phi$, then, with $n=\dim X$ and $V=\Omega^n$,
\[
 \Mab^{\mathrm{NA}}_\Omega(\phi)
 =\DF(\mathcal X,\mathcal A)
  +V^{-1}\bigl((\mathcal X_{0,\mathrm{red}}-\mathcal X_0)
                    \cdot\mathcal A^n\bigr).
\]
\end{samepage}
To eliminate the nonpositive correction term, make a sufficiently
divisible degree-$d$ base change and semistable reduction of an
equivariant resolution dominating $X\times\mathbb P^1$.
The resulting smooth model $\mathcal Y$ dominates the product
after base change and has reduced central fibre.
Let $\mathcal B$ be the pulled-back class. By homogeneity under base change and
birational invariance, we have
\[
 d\,\Mab^{\mathrm{NA}}_\Omega(\phi)
 =\Mab^{\mathrm{NA}}(\mathcal Y,\mathcal B)
 =\DF(\mathcal Y,\mathcal B).
\]
See Sj\"ostr\"om Dyrefelt~\cite[Lemma 3.16 and Section 3.6]{SDSemi} and Boucksom, Hisamoto, and Jonsson~\cite{BHJ}.
Although $\mathcal B$ is only relatively nef, choose a relatively
K\"ahler class $\mathcal C$ on $\mathcal Y$ with the same
general-fibre restriction $\Omega$, as in Step 6 of
Lemma~\ref{lem:real-entropy-recovery}. Then
$\mathcal B_\epsilon=(1-\epsilon)\mathcal B+\epsilon\mathcal C$
is relatively K\"ahler in that exact class. By semistability and continuity of intersection numbers, we obtain
\[
 \DF(\mathcal Y,\mathcal B)
 =\lim_{\epsilon\downarrow0}
       \DF(\mathcal Y,\mathcal B_\epsilon)\ge0.
\]
Birational invariance follows from the projection formula for
exceptional divisors. The perturbation takes place on the fixed
reduced model, where the correction term vanishes and the
Donaldson--Futaki invariant is a continuous intersection polynomial.

\paragraph{Step 2.}
For $\phi\in\mathcal E^{1,\mathrm{NA}}(\Omega)$ with finite
entropy, choose the relatively K\"ahler model approximations
$\phi_j$ of Lemma~\ref{lem:real-entropy-recovery}. By the lemma, the Mabuchi invariants converge. It therefore follows from Step 1 that
\[
 \Mab^{\mathrm{NA}}_\Omega(\phi)
 =\lim_j\Mab^{\mathrm{NA}}_\Omega(\phi_j)\ge0.
\]
If the entropy is infinite, the nonentropy energy terms are still
finite on $\mathcal E^{1,\mathrm{NA}}(\Omega)$, and hence
$\Mab^{\mathrm{NA}}_\Omega(\phi)=+\infty$. Consequently
\[
 \Mab^{\mathrm{NA}}_\Omega(\phi)\ge0
 \qquad(\phi\in\mathcal E^{1,\mathrm{NA}}(\Omega)).
\]

\paragraph{Step 3.}
Fix a smooth potential $x$ and let $\ell$ be a nonconstant
$\EE^2$ geodesic ray starting at $x$, parametrized with constant
speed $\|\ell\|_{d_2}=d_2(\ell_0,\ell_1)>0$.
By convexity, $(\Mab(\ell_t)-\Mab(x))/t$ is nondecreasing in $t$. The global
slope inequality at the smooth point $x$ gives
\[
 \Mab(\ell_t)\ge\Mab(x)
       -\sqrt{\Ca(x)}\,t\|\ell\|_{d_2}.
\]
Thus its asymptotic slope exists in $\mathbb R\cup\{+\infty\}$.
In particular, a slope of $-\infty$ cannot occur. If this slope
is finite, convexity also ensures that $\Mab(\ell_t)<\infty$
for every finite $t$. Otherwise all later secant slopes would
be infinite. Regard $\ell$ as an $\EE^1$ ray and apply the maximality criterion of
Mesquita-Piccione~\cite[Proposition 6.3.1]{MP}. The ray is maximal and corresponds
to a potential $\ell^{\mathrm{NA}}\in
\mathcal E^{1,\mathrm{NA}}(\Omega)$. Combining the real-class slope
comparison of Mesquita-Piccione~\cite[Corollary 6.2.3]{MP} with Step 2, we obtain
\[
 \Mab^\infty(\ell)
 :=\lim_{t\to\infty}\frac{\Mab(\ell_t)}{t}
 \ge\Mab^{\mathrm{NA}}_\Omega(\ell^{\mathrm{NA}})\ge0.
\]
The comparison extends the polarized results of Li~\cite{LiRays}.
A ray of slope $+\infty$ also has nonnegative slope. Thus every ray in the metric
duality formula has nonnegative Mabuchi slope.

\paragraph{Step 4.}
The duality theorem of Xia~\cite[Theorem 1.1]{Xia}, valid in every compact
K\"ahler class, now gives
\[
 s_*=
 \max\left\{0,\sup_{\ell\ne0}
       \frac{-\Mab^\infty(\ell)}{\|\ell\|_{d_2}}\right\}=0.
\]
Here $\ell\ne0$ means a nonconstant ray based at $x$. The
constant-shift ray $\ell_t=x+t$ has zero Mabuchi slope, so the
supremum is at least zero. In Xia's
notation the extended Calabi functional is the metric slope.
Our smooth Calabi energy is its square. Therefore
Proposition~\ref{prop:weak-classical-slope} gives
\[
 \inf_{u\in\HH}\sqrt{\Ca(u)}=s_*=0,
 \qquad \inf_{u\in\HH}\Ca(u)=0.
\]
Scaling a K\"ahler metric by $\lambda>0$ multiplies its scalar
curvature by $\lambda^{-1}$ and its probability-normalized Calabi
energy by $\lambda^{-2}$. Thus the infimum is also zero in
$2\pi\Omega$. The $\partial\bar\partial$-lemma represents every
smooth K\"ahler metric in that class by a smooth potential relative
to a fixed reference form. This proves \eqref{eq:intro-semistable-small-energy} and
yields a smooth metric below each positive energy threshold. For $\Omega=c_1(L)$, algebraic and cohomological
K-semistability agree by Sj\"ostr\"om Dyrefelt~\cite[Proposition 3.14]{SDSemi}.
\end{proof}

We now use these small-energy metrics as initial data for Calabi flow. Assume
that $X$ is a compact connected K\"ahler surface and
$\Omega\in\mathcal U$ is a real K\"ahler class. These conditions
force $X$ to be rational. Indeed,
$K_X\cdot\Omega=-c_1(X)\cdot\Omega<0$, so no positive multiple
of $K_X$ has a nonzero section. The zero divisor of such a section would be effective
and would have nonnegative intersection with $\Omega$.
Thus $X$ has Kodaira dimension $-\infty$. By the classification of compact K\"ahler surfaces, $X$ is rational or a blowup of a minimal ruled surface over a curve of genus $g\ge1$. See Barth, Hulek, Peters, and Van de Ven~\cite{BHPV}. In the latter case
$c_1(X)^2=8(1-g)-k\le0$, where $k\ge0$ is the number of
blowups, contrary to $\delta(\Omega)>0$. Hence $X$ is rational,
and in particular projective with $H^{2,0}(X)=0$ and
$\chi(\mathcal O_X)=1$. The Lefschetz $(1,1)$ theorem then gives
$H^{1,1}(X,\mathbb R)=\operatorname{NS}(X)_{\mathbb R}$,
and its K\"ahler cone is the real ample cone.

For
$[\omega]=2\pi\Omega$ and
$dV_\omega=\omega^2/2$, set
\[
 \mathcal V_\Omega=\int_XdV_\omega=2\pi^2\Omega^2,
 \qquad \underline S=\frac{2c_1(X)\cdot\Omega}{\Omega^2}.
\]
Since $d\mu_\omega=dV_\omega/\mathcal V_\Omega$,
\[
 \int_X(S_\omega-\underline S)^2dV_\omega
       =\mathcal V_\Omega\Ca(\omega).
\]
In Li, Wang, and Zheng~\cite[Corollary 1.3]{LWZSobolev}, the class is $[\omega]$
and the volume convention is $V=\int_X\omega^2=4\pi^2\Omega^2$.
Its centered normalized Calabi energy is our $\Ca$ and its
$\delta_{[\omega]}$ equals $\delta(\Omega)$. Consequently its
threshold $16\pi^2\delta_{[\omega]}/V$ becomes
\begin{equation}\label{eq:normalized-cone-threshold}
 \kappa_\Omega=\frac{4\delta(\Omega)}{\Omega^2}>0.
\end{equation}

\begin{corollary}
\label{cor:semistable-global}
If $(X,\Omega)$ is cohomologically K-semistable, then for any
$\varepsilon>0$ there is a smooth K\"ahler metric
$\omega_0\in2\pi\Omega$ with
\[
 \Ca(\omega_0)<\min\{\varepsilon,\kappa_\Omega\}.
\]
Its Calabi flow exists smoothly for all time and has uniformly
bounded Sobolev constants.
In particular this applies when $\Omega=c_1(L)$ and $(X,L)$ is
K-semistable.
\end{corollary}
\begin{proof}
By the preceding classification argument, $X$ is projective and
$\Omega$ is a real ample N\"eron--Severi class.
Corollary~\ref{cor:semistable-small-energy} supplies the smooth
initial metric. Its unnormalized centered energy satisfies
\[
 \int_X(S_{\omega_0}-\underline S)^2dV_{\omega_0}
       <8\pi^2\delta(\Omega).
\]
Along the maximal smooth Calabi flow the class is fixed and
$\Ca(\omega(t))\le\Ca(\omega_0)$, by the Calabi-energy
dissipation identity of Chen and He~\cite{ChenHeStability}. The strict initial
gap therefore persists. The argument of Li, Wang, and Zheng~\cite[Corollary
1.3]{LWZSobolev} supplies a time-independent
Sobolev bound determined by this gap and the class. If the maximal smooth time were
finite, the solution would extend past it by the extension theorem of Li, Wang, and
Zheng~\cite[Theorem 1.2]{LWZSobolev}. This proves global smooth existence
and the asserted uniform Sobolev bound.
For the last assertion, the polarized case of
Corollary~\ref{cor:semistable-small-energy} supplies the initial
metric, and the same flow criterion applies.
\end{proof}

To identify subsequential limits of this flow, we need the limiting
energy. Corollary~\ref{cor:semistable-small-energy} and
Theorem~\ref{thm:smooth-weak-infima} give $s_*=0$. By
Berman, Darvas, and Lu~\cite[Proposition 6.1]{BDLflow}, the global smooth flow is the
weak Mabuchi gradient flow with the same initial potential.
Combining the gradient-flow slope limit in
Xia~\cite[Theorem 4.1 and Remark 4.2]{Xia} with the smooth-point
identity of Xia~\cite[Proposition 3.6]{Xia}, we obtain
\begin{equation}\label{eq:controlled-flow-energy-limit}
 \lim_{t\to\infty}\sqrt{\Ca(\omega(t))}=s_*=0.
\end{equation}
\begin{proof}[Proof of Theorem~\ref{thm:stability}]
The first two steps address the existence assertion in a real class. The last two
establish the equivalences for discrete automorphism groups and for polarizations.

\paragraph{Step 1.}
Suppose that $(X,\Omega)$ is cohomologically K-polystable. It is
K-semistable, so Corollary~\ref{cor:semistable-global} supplies a global
smooth flow in $2\pi\Omega$ with uniformly bounded Sobolev constants.
Its energy tends to zero by \eqref{eq:controlled-flow-energy-limit}.
For any sequence $t_j\to\infty$, apply the compactness theorem of Li, Wang, and
Zheng~\cite[Theorem 1.1(I)--(IV)]{WeakCompanion} to obtain a compact extremal K\"ahler
orbifold limit $(Z,\omega_Z)$ along a subsequence.
By the same theorem, the limit has the same volume and average scalar curvature, and
the unnormalized Calabi functional converges to its value at the limit.
Since $\mathcal V_\Omega$ is fixed,
\[
 \int_Z(S_{\omega_Z}-\underline S)^2dV_{\omega_Z}
 =\lim_j\mathcal V_\Omega\Ca(\omega(t_j))=0.
\]
The scalar curvature is continuous in the orbifold charts, so
$S_{\omega_Z}\equiv\underline S$. The limiting eternal flow is stationary.

\paragraph{Step 2.}
We apply the real-class degeneration theorem. The initial
metric chosen above satisfies
\[
 \int_X(S_{\omega_0}-\underline S)^2dV_{\omega_0}
 =\mathcal V_\Omega\Ca(\omega_0)
 <\mathcal V_\Omega\kappa_\Omega=8\pi^2\delta(\Omega).
\]
Together with $\Omega\in\mathcal U$, this verifies precisely
the hypotheses of Li, Wang, and Zheng~\cite[Theorem 1.1]{ControlledCompanion}. Their
conventions for scalar curvature, volume, and the K\"ahler class agree with ours.
Put $\Omega_Z=[\omega_Z]/(2\pi)$. By
Li, Wang, and Zheng~\cite[Proposition 4.4]{ControlledCompanion}, there is a normal
relatively K\"ahler cohomological test configuration
$(\mathcal X,\mathcal A)$ for $(X,\Omega)$ with central fibre
$(Z,\Omega_Z)$. Since $\omega_Z$ is cscK, the test configuration has $\DF=0$.
Its invariant differs from our volume-normalized convention by
a positive class-dependent factor, which preserves its sign and zeros.

By cohomological K-polystability, this configuration is a product
over $\mathbb C$. To identify the class, follow
Li, Wang, and Zheng~\cite[Section 4.5]{ControlledCompanion}. Compose a product
identification with a constant automorphism of $X$ so that its
identification on the fibre at $1$ agrees with the prescribed
general-fibre marking. The inclusions of the fibres at $0$ and
$1$ in $X\times\mathbb C$ are homotopic. The pullback of
$\mathcal A$ therefore restricts to the same de Rham class
$\Omega$ on both fibres. The $\partial\bar\partial$-lemma on
$X$ identifies this with equality of real $(1,1)$-classes.
Thus the induced biholomorphism $F:X\to Z$ satisfies
$F^*\Omega_Z=\Omega$.

Under the product identification, the coarse space $Z$ is smooth.
Choose an effective local chart $B/\Gamma$ with $\Gamma$ acting
linearly and freely on $B\setminus\{0\}$.
In smooth coordinates on the quotient, the quotient map
$\Psi:B\to\mathbb C^2$ has nonzero Jacobian determinant away
from $0$. By Hartogs' theorem, $1/\det D\Psi$ extends to a holomorphic function
$h$ with $h\det D\Psi=1$ on the punctured ball and hence on $B$.
Thus $D\Psi(0)$ is invertible. Differentiating $\Psi\circ\gamma=\Psi$ at $0$, we obtain $\gamma=I$ for every
$\gamma\in\Gamma$. Consequently the orbifold-smooth metric
$\omega_Z$ is an ordinary smooth metric. Its pullback is cscK
and satisfies
\[
 [F^*\omega_Z]=2\pi F^*\Omega_Z=2\pi\Omega.
\]

\paragraph{Step 3.}
When $\Aut^0(X,J)=\{1\}$, every homomorphism
$\mathbb C^*\to\Aut(X,J)$ has image in the trivial identity
component. Thus the product actions are trivial and
cohomological K-polystability is K-stability. Conversely, cscK
existence implies cohomological K-stability in this case by
Dervan and Ross~\cite[Theorem 1.2]{DRKahler}. Together with Steps 1 and 2,
this gives the stated equivalence under the
automorphism hypothesis.

\paragraph{Step 4.}
For the polarized assertion, algebraic K-polystability implies
algebraic K-semistability. By Corollaries~\ref{cor:semistable-small-energy} and~\ref{cor:semistable-global}, there is a global smooth flow in $2\pi c_1(L)$ with uniformly bounded Sobolev constants.
The Sobolev convention in Li, Wang, and Zheng~\cite{PolarizedCompanion} uses
$dV=\omega^2/2$ and the fixed coefficient
$\mathcal V_\Omega^{-1/2}$. It is equivalent to our convention
by fixed factors in this class. Since $X$ is rational,
$\chi(\mathcal O_X)=1$ and $c_1(X)\cdot c_1(L)>0$.
These are exactly the analytic and numerical hypotheses of
Li, Wang, and Zheng~\cite[Theorem 1.1]{PolarizedCompanion}.
Part (4) gives a smooth cscK limit $(Z,\omega_Z)$ and a polarized
isomorphism $F:(X,L^r)\to(Z,A_Z)$ for some common integer
$r>0$, with $[\omega_Z]=(2\pi/r)c_1(A_Z)$. In particular,
\[
 [F^*\omega_Z]=\frac{2\pi}{r}c_1(F^*A_Z)
             =\frac{2\pi}{r}c_1(L^r)=2\pi c_1(L).
\]
Thus $F^*\omega_Z$ is the required smooth cscK metric on $(X,J)$.
For the reverse implication, rescale a cscK
metric in $2\pi c_1(L)$ by $1/(2\pi)$ and apply the theorem of Berman, Darvas, and
Lu~\cite[Theorem 1.6]{BDL} in $c_1(L)$.
\end{proof}

\begin{remark}
\label{rem:controlled-cone-role}
In the proof of Theorem~\ref{thm:stability}, we use the assumption $\Omega\in\mathcal U$ three times.
First, it makes $X$ rational, so $\Omega$ is a real ample
N\"eron--Severi class and
Corollary~\ref{cor:semistable-small-energy} applies.
Second, since $\delta(\Omega)>0$, the threshold
$\kappa_\Omega=4\delta(\Omega)/\Omega^2$ is positive. Under the strict energy bound of
Corollary~\ref{cor:semistable-global}, the flow exists for all time with uniformly
bounded Sobolev constants.
Third, it verifies the hypotheses of the class-preserving
degeneration theorem of
Li, Wang, and Zheng~\cite[Proposition 4.4]{ControlledCompanion}.

For $\delta(\Omega)\le0$, the threshold $\kappa_\Omega$ is
nonpositive, so this flow criterion cannot be met.

In the polarized case, the following assumptions suffice without a restriction to the controlled cone.
Suppose that $(X,L)$ is K-polystable and
\[
 \chi(\mathcal O_X)=1,\qquad c_1(X)\cdot c_1(L)>0.
\]
If $2\pi c_1(L)$ admits a global smooth Calabi flow with uniformly
bounded Sobolev constants, then it contains a smooth cscK metric,
by Li, Wang, and Zheng~\cite[Theorem 1.1(4)]{PolarizedCompanion}.

The cone is not a necessary condition for cscK existence itself.
For example, on $\mathbb P^1\times\mathbb P^1$ let $H_1,H_2$
be the two hyperplane classes and take $\Omega=4H_1+H_2$.
Then
\[
 \delta(\Omega)=8-\frac{2\cdot10^2}{3\cdot8}=-\frac13,
\]
whereas the product of appropriately scaled Fubini--Study metrics
is cscK in $2\pi\Omega$.

\end{remark}

\section{Small Calabi energy and global Calabi flow}
\label{sec:sobolev-application}\label{sec:surface-setup}
Throughout this section, $X$ is a compact connected K\"ahler surface
with cscK reference form $\omega_*$. We retain
\eqref{eq:energy-definition} and \eqref{eq:sobolev-definition},
with $V=\int_X\omega_*^2$. By geodesic convexity, a cscK potential minimizes the
Mabuchi functional. See Berman and Berndtsson~\cite{BB}. Thus we may take
$\Mab(0)=\inf\Mab=0$.
The static estimates below precede the local compactness argument, which also requires
a critical Sobolev bound. In the proof of Theorem~\ref{thm:barrier}, the assumption
$\CS=A$ gives exactly this bound.

For the associated Riemannian metric $g_u$, fix $c,\gamma>0$ by
$d\mu_u=c\,dV_{g_u}$ and
$D_u(f)=\gamma\int_X|df|_{g_u}^2d\mu_u$.
With $dV_g=\omega^2/2$ and the usual real gradient, these constants
are $c=2/V$ and $\gamma=1/2$. Let $S_4$ be the optimal constant in
$\|f\|_{L^4(dx)}^2\le S_4\int_{\mathbb R^4}|df|^2dx$ for
$f\in C_c^\infty(\mathbb R^4)$. See Talenti~\cite{Talenti}.
With the measure $c\,dx$ and Dirichlet energy $\gamma\int|df|^2c\,dx$, the homogeneous
inequality has optimal constant
$\CE=S_4/(\gamma\sqrt c)$. By concentrating test functions in a coordinate ball, we obtain
\begin{equation}\label{eq:euclidean-threshold}
 \CS(g)\ge\CE.
\end{equation}
Indeed, take as a test function in normal coordinates
$f_\epsilon(x)=\epsilon^{-1}f(x/\epsilon)$, extended by zero.
As $\epsilon\downarrow0$, the integral of its fourth power and its Dirichlet energy
tend to their Euclidean values, while its squared $L^2$
norm is $O(\epsilon^2)$. The Sobolev quotient therefore tends
to the homogeneous Euclidean quotient of $f$. The assertion follows by taking the supremum over $f$.
The rescaled measures retain the factor $c$ under blow-up. The constants below may depend on these fixed normalizations.

\label{sec:static}

Write $V=\int_X\omega_*^2$, $d\mu_*=V^{-1}\omega_*^2$, and
$d\mu_\varphi=e^F d\mu_*$. Both measures have mass one. Set
\begin{equation}\label{eq:static-notation}
 f=S_\varphi-\underline S,\quad E=\Ca(\varphi)=\int f^2d\mu_\varphi,
 \quad H=\int Fd\mu_\varphi,\quad
 Q=\tr_{\omega_*}\omega_\varphi,\quad T=\tr_{\omega_\varphi}\omega_*.
\end{equation}
Here $F=\log(\omega_\varphi^2/\omega_*^2)$ and
$D_\varphi(u)=\int|\partial u|_\varphi^2d\mu_\varphi$.
The trace identities and the scalar curvature equation read
\begin{gather}
 T=Qe^{-F},\qquad QT\ge4,\qquad \int Qd\mu_*=2,
 \label{eq:static-surface}\\
 \Delta_\varphi F=\tr_{\omega_\varphi}\Ric(\omega_*)-\underline S-f.
 \label{eq:scalar-static}
\end{gather}
We first choose representatives modulo holomorphic automorphisms.
For these representatives, small Calabi energy gives the entropy bound
needed in the estimates below without a prior bound on the Sobolev constant.

\begin{lemma}\label{lem:entropy-gauge}
Normalize $\AM(u)=0$ and $\inf\Mab=\Mab(0)=0$, and let
$G=\Aut_0(X,J)$ and $d_{2,G}(0,u)=\inf_{\sigma\in G}d_2(0,\sigma\cdot u)$.
There are background constants $C,H_0,\varepsilon_g>0$ such that
\[
 d_{2,G}(0,u)\le C(1+\Mab(u)),\qquad
 \Ca(u)\le\varepsilon_g\ \Longrightarrow\
 0\le\Mab(u)\le C\sqrt{\Ca(u)}.
\]
Each potential in this sublevel has a representative with $H\le H_0$.
If $\Ca(u_j)\to0$, then, after passing to a subsequence and taking suitable pullbacks, we have
\begin{equation}\label{eq:entropy-gauge-limit}
 H(u_j)\longrightarrow0,\qquad d_2(0,u_j)\longrightarrow0.
\end{equation}
\end{lemma}
\begin{proof}
The normalized action of $G$ is defined by $\omega_{\sigma\cdot u} =\sigma^*\omega_u$ and $\AM(\sigma\cdot u)=0$. We first bound the distance to the reference in terms of entropy.

\paragraph{Step 1.}
Put $\psi=u-\sup u-1$ and $K=\int|\psi|^2d\mu_u\ge1$.
Combining the Moser--Trudinger inequality of Di Nezza, Guedj, and Lu~\cite[Theorem B]{DGL}, with $n=p=2$, and the entropy variational inequality, we obtain
\[
 \int e^{c|\psi|^2/K^{1/2}}d\mu_*\le C,
 \qquad cK^{1/2}\le H(u)+\log C.
\]
By exponential integrability of normalized plurisubharmonic functions, $\psi$ has
bounded background $L^2$ norm. Moreover,
$0\le-\AM(\psi)\le\int(-\psi)d\mu_u\le K^{1/2}$.
We use the endpoint comparison of Darvas~\cite[Theorem 3, equation (7)]{Darvas} in the form
\[
 d_2(v,w)^2\le C\int_X|v-w|^2(d\mu_v+d\mu_w).
\]
Together with the identity $u=\psi-\AM(\psi)$, this gives
\begin{equation}\label{eq:entropy-dtwo}
 d_2(0,u)\le C(1+H(u)),\qquad \AM(u)=0.
\end{equation}
Here the entropy variational inequality is
$\int a\,d\nu\le\operatorname{Ent}(\nu\mid\mu_*)+
\log\int e^a d\mu_*$ for probability measures. It follows
from Jensen's inequality applied to $e^a/(d\nu/d\mu_*)$,
and is applied above with $a=c|\psi|^2/K^{1/2}$.

\paragraph{Step 2.}
The properness theorem of Berman, Darvas, and Lu~\cite[Theorem 1.5]{BDL} bounds the infimum of $J_*$ on $G\cdot u$. Choose $v\in G\cdot u$ with $J_*(v)<1+\inf_{G\cdot u}J_*$. Then
$J_*(v)=\int v\,d\mu_*-\AM(v)\le C(1+\Mab(u))$.
Since the Futaki character vanishes, $\Mab$ is $G$-invariant. We use the Chen--Tian formula
\begin{equation}\label{eq:chen-tian-static}
 \Mab(v)=H(v)+\underline S\AM(v)
 -V^{-1}\int v\,\Ric(\omega_*)\wedge(\omega_*+\omega_v)
\end{equation}
to bound the absolute value of the nonentropy part by $C(1+J_*(v))$.
Replace $v$ by $w=v-\sup v-1$, use $|\Ric(\omega_*)|_*\le C$, and
bound the mixed integrals of $-w$ by $3[-\AM(w)]$.
We have $-\AM(w)=\sup v+1\le J_*(v)+C$. Thus $H(v)\le C(1+\Mab(u))$.
Combining this with \eqref{eq:entropy-dtwo}, we obtain
$d_{2,G}(0,u)\le C(1+\Mab(u))$.
The endpoint inequality \eqref{eq:classical-mabuchi-support}, applied at
$v$ with competitor $0$, gives $\Mab(v)\le\sqrt{\Ca(v)}\,d_2(v,0)$.
Taking the infimum over the orbit and using invariance of both energies, we obtain
\[
 0\le\Mab(u)\le\sqrt{\Ca(u)}\,d_{2,G}(0,u)
 \le C\sqrt{\Ca(u)}(1+\Mab(u)).
\]
Absorb the last term when $C\sqrt{\varepsilon_g}\le1/2$.
The chosen representatives then have uniformly bounded entropy.

\paragraph{Step 3.}
For $\Ca(u_j)\to0$, these representatives have bounded entropy and
$J_*$ and satisfy $\Mab(u_j)\to0$. By the $d_1$ entropy compactness theorem of Di Nezza, Guedj, and Lu~\cite[Theorem 3.5]{DGL} and lower semicontinuity, a subsequence converges to a minimizer, which is smooth and cscK by Berman, Darvas, and Lu~\cite[Theorem 1.4]{BDL}.
The normalizing constants are bounded, so the compactness theorem,
stated for supremum-normalized potentials, applies also in our
$\AM=0$ normalization. Denote the limit by $u_\infty$.
By the same theorem, a fixed $\sigma\in G$ satisfies
$\sigma^*\omega_{u_\infty}=\omega_*$. Pull back the entire
subsequence by this $\sigma$. The normalized action is a
$d_1$-isometry, hence the new potentials tend to $0$ in $d_1$.
The nonentropy terms in \eqref{eq:chen-tian-static} are
$d_1$-continuous by Berman, Darvas, and Lu~\cite[Propositions 4.1 and 4.4]{BDLflow} and
vanish at $0$. Since $\Mab(u_j)\to0$, the same formula now
gives $H(u_j)\to0$ relative to the original reference $\omega_*$.
By the preceding exponential estimate and the boundedness of the normalizing constants, we have $\int e^{c|u_j|^2}d\mu_*\le C$. Convergence in $d_1$ implies convergence in background measure.
The functions $e^{c|u_j|^2/2}$ have uniformly bounded $L^2(\mu_*)$
norms, so Vitali's theorem gives
$\int e^{c|u_j|^2/2}d\mu_*\to1$.
It follows that
\[
 (c/2)\int|u_j|^2d\mu_{u_j}
 \le H(u_j)+\log\int e^{c|u_j|^2/2}d\mu_*\longrightarrow0.
\]
Since $x\le e^x-1$ for $x\ge0$, the same exponential integral also gives
$\int|u_j|^2d\mu_*\to0$. Applying the comparison of Darvas~\cite[Theorem 5.5]{Darvas}
to these two endpoint integrals, we obtain $d_2(0,u_j)\to0$.
\end{proof}

Once the entropy is bounded, we can use the small-energy assumption to control moments
of the volume density and its reciprocal. These bounds will enter the estimates for the
two metric traces.

\begin{lemma}\label{lem:density-moments}
Given $H_0<\infty$ and a finite $P\ge1$, there are
$\varepsilon_P,C_P>0$, depending only on these data and the background,
such that, with $\sup\varphi=0$,
\begin{equation}\label{eq:density-moments}
 H\le H_0,\quad E\le\varepsilon_P
 \quad\Longrightarrow\quad
 \|\varphi\|_\infty+\int(e^{PF}+e^{-PF})d\mu_*\le C_P.
\end{equation}
\end{lemma}
\begin{proof}
The auxiliary potentials will control both tails of the volume density.

\paragraph{Step 1.}
Fix a background exponent $\alpha_0>0$ for which
\begin{equation}\label{eq:alpha-static}
 \sup_{\sup v=0,\ \omega_v>0}\int e^{-\alpha_0v}d\mu_*<\infty.
\end{equation}
Such an exponent exists by uniform exponential integrability of
normalized K\"ahler potentials. See Chen and Cheng~\cite[Proposition 5.1]{CC}.
For $E>0$, solve the smooth normalized Monge--Amp\`ere equations
\begin{equation}\label{eq:auxiliary-static}
 \omega_\psi^2=A_F^{-1}\sqrt{1+F^2}\,\omega_\varphi^2,
 \quad \omega_\rho^2=\frac{f^2+E}{2E}\omega_\varphi^2,
 \quad \sup\psi=\sup\rho=0,
\end{equation}
where $1\le A_F=\int\sqrt{1+F^2}d\mu_\varphi\le1+H_0+2/e=:A_0$.
Both right sides are smooth, strictly positive, and have the
same total mass as $\omega_*^2$, so solvability follows from
Yau~\cite{Yau}. For the bound on $A_F$, use
$\sqrt{1+F^2}\le1+|F|$ and
$\int_{\{F<0\}}(-F)e^F d\mu_*\le1/e$.
By the arithmetic--geometric mean inequality, we have
\begin{equation}\label{eq:aux-traces}
 \Delta_\varphi\psi\ge2A_F^{-1/2}(1+F^2)^{1/4}-T,
 \qquad\Delta_\varphi\rho\ge\sqrt{2/E}\sqrt{f^2+E}-T.
\end{equation}
Choose $\lambda$ so that $\pm\Ric(\omega_*)+(\lambda-2)\omega_*\ge0$.
For $0<a\le1$ and $U_\pm=\pm F+a(\psi+\rho)-\lambda\varphi$,
\begin{equation}\label{eq:aux-absorption}
 \Delta_\varphi U_\pm\ge k(1+F^2)^{1/4}-C_0,
 \quad k=2aA_0^{-1/2},\quad E\le2a^2,
\end{equation}
since $\mp f+a\sqrt{2/E}\sqrt{f^2+E}\ge0$.

\paragraph{Step 2.}
Choose $0<\delta\le\min(1,\alpha_0/(4\lambda))$ and
$0<a\le\min(1/4,\alpha_0/[8(2+\delta)])$.
Integrating $\Delta_\varphi e^{\delta U_+}$, we bound the integral of $e^{\delta U_+}$
over $\{k(1+F^2)^{1/4}\ge C_0+1\}$ by $C_0$ times the integral
over the complement, with both integrals taken against $d\mu_\varphi$.
On the complement $|F|$ is bounded. Since $\psi,\rho\le0$,
the latter integral is bounded by $C\int e^{-\delta\lambda\varphi}d\mu_*$.
Since $\varphi\le0$, we may drop the factor $e^{-\delta\lambda\varphi}\ge1$ to obtain
\begin{equation}\label{eq:first-weighted-density}
 \int e^{(1+\delta)F+a\delta(\psi+\rho)}d\mu_*\le C.
\end{equation}
To remove the auxiliary factor, apply H\"older's inequality with
$\vartheta=(1+\delta/2)/(1+\delta)$ and use $2a(2+\delta)\le\alpha_0$. This gives
$\int e^{(1+\delta/2)F}d\mu_*\le C$.
The estimate of Ko\l odziej~\cite{Kolodziej} now bounds $\|\varphi\|_\infty$.

\paragraph{Step 3.}
The same argument now applies to either sign. For any fixed $\delta>0$, it gives
\begin{equation}\label{eq:both-weighted-density}
 \int e^{(1\pm\delta)F+a\delta(\psi+\rho)}d\mu_*\le C_\delta.
\end{equation}
Take $\delta=2P+2$ and $a\le\min(1/4,\alpha_0/[8(P+2)])$.
For the plus sign, the preceding H\"older exponent gives the
$(P+2)$th positive density moment. For the minus sign take
$\vartheta=P/(2P+1)$. The auxiliary cost is admissible because
$2a\delta\vartheta/(1-\vartheta)=4Pa\le\alpha_0$.
This proves \eqref{eq:density-moments}, with the minimum of the
finitely many thresholds $E\le2a^2$. If $E=0$, omit $\rho$.
\end{proof}

The density controls the product of the two eigenvalues of the metric
relative to $\omega_*$. To estimate the metric and its inverse separately,
we now bound the traces $Q$ and $T$ in finite $L^r$ norms.

\begin{proposition}\label{prop:finite}
Fix $H_0<\infty$ and a finite exponent $r\ge1$. Suppose that $H\le H_0$ and that $E$
is sufficiently small in terms of $r,H_0$ and the background. Then
$\|Q\|_{L^r(\mu_*)}+\|T\|_{L^r(\mu_*)}\le C_r$.
\end{proposition}
\begin{proof}
We estimate a weighted trace and then remove the weight.

\paragraph{Step 1.}
Normalize $\sup\varphi=0$.
The logarithmic trace calculation of Chen and Cheng~\cite[Section 3]{CC} gives
\begin{equation}\label{eq:log-trace-static}
 \Delta_\varphi\log Q\ge Q^{-1}\Delta_*F-C_0T.
\end{equation}
The third-derivative term cancels the negative gradient term in
$\Delta_\varphi\log Q$. The background zeroth-order terms are estimated using $Q^{-1}\le T/4$.
Fix $p\ge2$, $\alpha=p+1$, $0<b\le1$, and put
$A=-\alpha(F+\lambda\varphi)+b\rho$, $v=e^AQ$, $B=A+F$.
For a sufficiently large background $\lambda$, it follows from the scalar equation and \eqref{eq:aux-traces} that
\begin{equation}\label{eq:weighted-trace-static}
 \Delta_\varphi v\ge\kappa Tv-C_Av+e^A\Delta_*F
                   +v^{-1}|\partial v|_\varphi^2,
 \quad\kappa\ge2,\quad E\le2b^2/\alpha^2.
\end{equation}
Indeed, the only remaining scalar contribution is nonnegative because
$\alpha f+b\sqrt{2/E}\sqrt{f^2+E}\ge0$. Here $C_A$ depends on $p$.
For $E=0$ set $b=0$.

\paragraph{Step 2.}
In this step all unmarked integrals are taken against $d\mu_\varphi$. Set
$G_p=\int v^{p-2}|\partial v|_\varphi^2$, $I=\int Tv^p$, $J=\int v^p$.
Integration gives
\begin{align}
 pG_p+\kappa I&\le C_AJ+R,\label{eq:trace-ibp-static}\\
 R&=\frac1{\alpha-1}\int e^Bv^{p-1}
 (\Delta_*B+\alpha\lambda\Delta_*\varphi-b\Delta_*\rho)d\mu_*.
 \notag
\end{align}
Completing the square and using $|\partial v|_*^2\le Q|\partial v|_\varphi^2$, we obtain $\int e^Bv^{p-1}\Delta_*B\,d\mu_*\le (p-1)^2G_p/4$.
In detail, integration by parts makes this integral equal to
\[
 -\int e^Bv^{p-1}|\partial B|_*^2d\mu_*
 -(p-1)\operatorname{Re}\int e^Bv^{p-2}
                \langle\partial v,\partial B\rangle_*d\mu_*.
\]
It is at most
$\frac{(p-1)^2}{4}\int e^Bv^{p-3}|\partial v|_*^2d\mu_*$.
Since $e^{B-F}Q=e^AQ=v$, conversion to $d\mu_\varphi$
gives exactly the asserted multiple of $G_p$.
After factoring out their coefficients, we can bound the other two terms by
$J$ and $I/2$, respectively, because $\Delta_*\varphi=Q-2$,
$-\Delta_*\rho\le2$, and $Q^{-1}\le T/4$. It follows that
\[
 \left(p-\frac{(p-1)^2}{4(\alpha-1)}\right)G_p
 +\left(\kappa-\frac b{2(\alpha-1)}\right)I
 \le\left(C_A+\frac{\alpha\lambda}{\alpha-1}\right)J.
\]
Thus $G_p+I\le C_pJ$. To bound $J$ itself, use $T=Qe^{-F}=ve^{-B}$ and H\"older's inequality to obtain
\begin{equation}\label{eq:trace-GIJ}
 I=\int e^{-B}v^{p+1},\quad
 J\le I^{p/(p+1)}K^{1/(p+1)},\quad
 K=\int e^{pB}d\mu_\varphi\le C_p.
\end{equation}
The bound for $K$ uses $\rho\le0$, bounded $\varphi$, and the density
moment with exponent $|1-p(\alpha-1)|$. Hence $G_p+I+J\le C_p$.

\paragraph{Step 3.}
To remove the weight, set $p=2r$ and $b\le\min(1,\alpha_0/(8r))$.
Applying H\"older's inequality twice, we obtain
\[
 \int Q^rd\mu_*\le C_r
 \left(\int e^{(4\alpha r-2)F}d\mu_*\right)^{1/4}
 \left(\int e^{-4br\rho}d\mu_*\right)^{1/4}\le C_r.
\]
Indeed, the first Cauchy--Schwarz factor is bounded using $J=\int e^{pA+F}Q^p d\mu_*\le
C_p$. In the other factor $e^{p\alpha\lambda\varphi}\le1$,
and \eqref{eq:alpha-static} controls the exponential of $-\rho$. Finally, using $T=Qe^{-F}$ and the $2r$th moment of $Q$, we obtain a bound for $T$ by H\"older's inequality.
Only finitely many moments and smallness conditions were used.
\end{proof}

For a sequence with Calabi energy tending to zero, we can choose
representatives whose entropies tend to zero. The trace estimates then
allow us to pass from convergence in measure to strong convergence of
the metrics and their inverses.

\begin{proposition}\label{prop:strong}
Suppose that $\Ca(\omega_j)\to0$. After passing to a subsequence and pulling back by
holomorphic automorphisms, we have, for every fixed finite $p$,
\begin{gather}
 g_j\to g_*,\quad g_j^{-1}\to g_*^{-1},\quad e^{\pm F_j}\to1
 \quad\hbox{strongly in }L^p(d\mu_*),\label{eq:static-strong}\\
 \int|\Ric(\omega_j)-\Ric(\omega_*)|_j^2d\mu_j\longrightarrow0.
 \label{eq:static-ricci-error}
\end{gather}
\end{proposition}
\begin{proof}
Use Lemma~\ref{lem:entropy-gauge} to choose a subsequence and pullbacks
for which $H_j\to0$.

\paragraph{Step 1.}
Put $A_j=g_*^{-1}g_j$, $q_j=\tr A_j$, and $a_j=\sqrt{\det A_j}$. Then we have $|A_j-a_jI|_*^2=q_j(q_j-2a_j)$, $\int q_j=2$, and $\int a_j^2=1$.
It now follows from $t\log t-t+1\ge(\sqrt t-1)^2$ that
\begin{equation}\label{eq:entropy-matrix}
 \int|A_j-I|_*d\mu_*\le2\sqrt2\sqrt{H_j}\longrightarrow0.
\end{equation}
Indeed, $\int(a_j-1)^2=2(1-\int a_j)\le H_j$, while
$\int(q_j-2a_j)=2(1-\int a_j)$. Cauchy--Schwarz gives
$\int|A_j-a_jI|_*\le\sqrt{2H_j}$, and
$\int|(a_j-1)I|_*\le\sqrt{2H_j}$.
Thus $A_j\to I$ in measure. Inversion and the determinant are continuous
at $I$, so $A_j^{-1}$ and $e^{\pm F_j}$ converge in measure as well.
The norms of $A_j$ and $A_j^{-1}$ are bounded by $Q_j$ and $T_j$.
For each $p\ge1$, apply Lemma~\ref{lem:density-moments} and
Proposition~\ref{prop:finite} with exponent $2p$. It follows that the $p$th powers are
uniformly integrable. Vitali's theorem gives \eqref{eq:static-strong}.
For each fixed exponent, the estimates hold for all sufficiently large indices.

\paragraph{Step 2.}
Since $\tr_{\omega_*}\Ric(\omega_*)=\underline S$, the scalar equation gives
$\Delta_jF_j\to0$ in $L^2(d\mu_j)$. Here we use strong $L^4$ convergence of the inverse
metric with respect to the background measure, the uniform $L^2$ bound for the density,
and
$\|f_j\|_{L^2(\mu_j)}\to0$.
More explicitly, if
$b_j=(g_j^{i\bar k}-g_*^{i\bar k})(\Ric_*)_{i\bar k}$, then
\[
 \|b_j\|_{L^2(\mu_j)}^2
 \le\|b_j\|_{L^4(\mu_*)}^2\|e^{F_j}\|_{L^2(\mu_*)}
 \longrightarrow0,\qquad \Delta_jF_j=b_j-f_j.
\]

To pass from the Laplacian to the full Ricci error, consider the exact
real $(1,1)$-form $\eta_j=\ddc F_j$. Stokes' theorem and the
trace--norm identity give
$\int|\eta_j|_j^2d\mu_j=\int(\Delta_jF_j)^2d\mu_j$.
To see the normalization directly, diagonalize $\eta_j$ relative
to $\omega_j$ at a point. If its eigenvalues are $\lambda_1,
\lambda_2$, then
\[
 \eta_j\wedge\eta_j
 =\bigl((\tr_{\omega_j}\eta_j)^2-|\eta_j|_j^2\bigr)
       \frac{\omega_j^2}{2}.
\]
The integral on the left vanishes because $\eta_j$ is exact and $X$ is closed.
As $\Ric_j-\Ric_*=-\eta_j$, this proves the second assertion.
\end{proof}

The compactness argument below requires control of the full curvature on
sets of small volume. We obtain this from the Ricci estimate by comparing
the Chern curvatures on the fixed holomorphic tangent bundle.

\begin{proposition}\label{prop:nc}
For the preceding representatives, let $\mathcal R_j$ denote the
Chern curvature of $T^{1,0}X$. We regard it as a $(1,1)$-form
with values in $\operatorname{End}(T^{1,0}X)$ on this fixed
holomorphic bundle. Its squared norm is the sum of squared
components in a unitary frame. Curvature differences are taken
in the fixed bundle before lowering indices with a metric.
Then
\begin{equation}\label{eq:strong-chern}
 \int|\mathcal R_j-\mathcal R_*|_j^2d\mu_j\longrightarrow0.
\end{equation}
There exist a sequence $\eta_j\to0$ and a background constant $C$ such that every measurable
$\Omega\subset X$ satisfies
\begin{equation}\label{eq:static-nc}
 \int_\Omega|\Rm(g_j)|_j^2dV_j
 \le\eta_j+C\Vol_{g_j}(\Omega)^{1/2}.
\end{equation}
\end{proposition}
\begin{proof}
We compare the Chern connections before passing to their curvature energies.

\paragraph{Step 1.}
The holomorphic identity map $(X,g_j)\to(X,g_*)$ has energy density
$T_j$ and Bochner identity
\begin{equation}\label{eq:bochner-identity-map}
 \begin{aligned}
 \Delta_jT_j&=P_j+g_j^{i\bar l}g_j^{k\bar q}
 \bigl((\Ric_j)_{i\bar q}(g_*)_{k\bar l}-(R_*)_{i\bar l k\bar q}\bigr),\\
 P_j&=|\nabla d\operatorname{id}|_{j,j,*}^2\ge0.
 \end{aligned}
\end{equation}
This is the K\"ahler Bochner formula for a holomorphic map.
See Li, Wang, and Zheng~\cite[Section 3.2]{LWZSobolev}. Since the manifold is closed,
the integral of $\Delta_jT_j$ vanishes.
When we replace $\Ric_j$ by $\Ric_*$, the error in the integral is bounded by
$\|\Ric_j-\Ric_*\|_{L^2(\mu_j)}\|T_j\|_{L^2(\mu_j)}=o(1)$.
The remaining contractions tend in $L^1(\mu_j)$ to their cancelling
background values. Integration gives $\int P_jd\mu_j\to0$.
For $\mathcal A_j=\nabla_j^{1,0}-\nabla_*^{1,0}$, we have
$|\mathcal A_j|_*^2\le Q_j^2P_j$. For $1<q<2$,
\begin{equation}\label{eq:connection-strong}
 \int|\mathcal A_j|_*^qd\mu_*
 \le\left(\int P_jd\mu_j\right)^{q/2}
 \left(\int Q_j^{2q/(2-q)}e^{-2F_j/(2-q)}d\mu_j\right)^{(2-q)/2}
 \longrightarrow0.
\end{equation}
The second factor uses only finite moments. The fixed holomorphic
structure makes the $(0,1)$ connections identical, so
$\mathcal R_j-\mathcal R_*=\bar\partial\mathcal A_j$ distributionally.
There is no quadratic connection term of type $(1,1)$.

\paragraph{Step 2.}
By the Chern--Weil formulas on a K\"ahler surface, we have
\begin{equation}\label{eq:chern-energy-identity}
 \int|\mathcal R_j|_j^2d\mu_j
 -\int|\mathcal R_*|_*^2d\mu_*=\Ca(\omega_j)\longrightarrow0.
\end{equation}
Indeed, for the complex tensor norms and $dV=\omega^2/2$,
\begin{align*}
 \int(S^2-|\Ric|^2)dV&=4\pi^2c_1^2[X],\\
 \int(|\mathcal R|^2-2|\Ric|^2+S^2)dV&=8\pi^2c_2[X].
\end{align*}
Subtracting twice the first identity from the second, we obtain
$\int|\mathcal R|^2dV=\int S^2dV+8\pi^2(c_2-c_1^2)[X]$.
See also Chen and He~\cite[Section 2]{ChenHe}. Taking the difference from the cscK reference and using $d\mu=2dV/V$, we obtain
\eqref{eq:chern-energy-identity}. Real curvature conventions change
only a fixed positive factor.

\paragraph{Step 3.}
To pass to the limit in the mixed term, write
$|\Xi|_j^2d\mu_j=\langle W_j\Xi,\Xi\rangle_*d\mu_*$.
The matrices $W_j^{\pm1}$ depend algebraically on $g_j^{\pm1},e^{\pm F_j}$. Hence
$W_j^{\pm1}\to I$ in every finite background $L^p$.
By weighted H\"older's inequality and \eqref{eq:chern-energy-identity}, $\mathcal R_j$ is bounded in background $L^{3/2}$. Hence its distributional convergence is also weak $L^{3/2}$ convergence. Explicitly,
\[
 \int|\mathcal R_j|_*^{3/2}d\mu_*
 \le\left(\int\langle W_j\mathcal R_j,
                         \mathcal R_j\rangle_*d\mu_*\right)^{3/4}
       \left(\int|W_j^{-1}|^3d\mu_*\right)^{1/4}.
\]
Every weakly convergent subsequence has the same distributional
limit $\mathcal R_*$, which identifies the weak limit. In
\[
 \int\langle W_j\mathcal R_j,\mathcal R_*\rangle_*
 =\int\langle\mathcal R_j,\mathcal R_*\rangle_*
 +\int\langle\mathcal R_j,(W_j-I)\mathcal R_*\rangle_*
\]
all integrals use $\mu_*$. The first term converges by the weak convergence of the curvature tensors.
The second vanishes by strong $L^3$ convergence of $W_j-I$.
By expanding the square and taking twice the real part of this Hermitian mixed pairing, we obtain \eqref{eq:strong-chern}.
Finally $|\mathcal R_*|_j^2\le C Q_jT_j^3$, so its intrinsic fourth
moment is bounded. To obtain \eqref{eq:static-nc}, bound the squared norm of the sum by twice the sum of the squared norms and apply H\"older's inequality to the background curvature term. The real and complex curvature norms differ only by a fixed factor.
By restricting the global error integral to $\Omega$, we obtain the estimate for every
measurable set, including balls for $g_j$.
\end{proof}

\label{sec:smoothing}

We next use local Ricci-flow smoothing to prove flat compactness
when the local curvature energy vanishes. The smoothing theorem
is due to Lee and Lee~\cite[Theorem 1.1]{LL}.
Write $h_j$ for the initial metric, $h_j(t)$ for its smoothing,
and $H_j(t)=F_j^*h_j(t)$ on a fixed coordinate domain, with $F_j$
independent of time. Here $dx$ denotes coordinate Lebesgue measure,
while the Riemannian density is $\sqrt{\det H_j}\,dx$.
To pass from positive times to time zero, we first estimate how much the metric and its
inverse can change during a short smoothing interval.

\begin{lemma}\label{lem:zero-time-coefficients}
Let $H_j(t)$, $0\le t\le T$, be smooth local Ricci flows in a fixed
four-dimensional coordinate domain. Suppose that the metrics $Q_j=H_j(T)$ are
uniformly equivalent to the Euclidean metric, that $H_j(t)\to H(t)$ smoothly on
compact subsets for $t>0$, and that
\[
 |\Ric(H_j(t))|_{\rm op}\le a/t,\qquad
 \norm{\Ric(H_j(t))}_{L^2(H_j(t))}\le L.
\]
Assume that these estimates hold on a neighborhood of the domain under consideration.
If $P\ge2$ and $a<1/(P+2)$, then, with
$\theta=2/P-2a-4a/P>0$,
\begin{equation}\label{eq:time-zero-modulus}
 \norm{H_j(0)-H_j(t)}_{L^P(dx)}
 +\norm{H_j(0)^{-1}-H_j(t)^{-1}}_{L^P(dx)}\le Ct^\theta.
\end{equation}
The initial metrics and inverses consequently converge strongly in
local $L^P$ to mutually inverse positive definite fields $H_0,H_0^{-1}$.
\end{lemma}

\begin{proof}
Put $k=(T/t)^{2a}$. Integrating the Ricci bound, we obtain
\begin{equation}\label{eq:reference-metric-comparison}
 k^{-1}Q_j\le H_j(t)\le kQ_j,
 \qquad dV_{Q_j}\le k^2dV_{H_j(t)}.
\end{equation}
Interpolation gives $\norm{\Ric}_{L^P(H_j(t))}\le Ct^{-1+2/P}$.
For a covariant two-tensor, the change to the $Q_j$ norm introduces a factor of at most $k$. The same factor applies to the evolution of the inverse metric, since
$\partial_tH_j^{-1}=2H_j^{-1}\Ric H_j^{-1}$ and its intrinsic norm
is $2|\Ric|$. Thus we have
\begin{equation}\label{eq:reference-time-derivative}
 \norm{\partial_tH_j}_{L^P(Q_j)}
 +\norm{\partial_tH_j^{-1}}_{L^P(Q_j)}
 \le C(T/t)^{2a+4a/P}t^{-1+2/P}.
\end{equation}
The right side is integrable at zero. First integrate from a positive lower endpoint.
By smoothness of each flow at zero and Fatou's lemma, we obtain \eqref{eq:time-zero-modulus}.
The positive-time limits have the same modulus. In the triangle inequality, first let
$j\to\infty$ and then let $t\downarrow0$. This proves strong convergence. The product
identity is preserved under almost-everywhere convergence of both matrix fields. Since
the limiting fields are positive semidefinite and invertible, they are positive
definite.
\end{proof}

For $G>0$, put $\rho=c\sqrt{\det G}$ and $B=\gamma\rho G^{-1}$,
where $c,\gamma>0$ account for the fixed volume and Dirichlet conventions.
By a direct calculation with the eigenvalues in real dimension four, we have
\begin{equation}\label{eq:determinant-cancellation}
 \rho^{\pm1}\le C|G^{\pm1}|^2,\quad
 |B|\le C|G|^{3/2}|G^{-1}|^{1/2},\quad
 |B^{-1}|\le C|G|^{1/2}|G^{-1}|^{3/2}.
\end{equation}
Consequently, if $G_j,G_j^{-1}$ converge strongly in $L^P$, $P>2$, then
\begin{equation}\label{eq:strong-density-coefficients}
 \rho_j^{\pm1}\to\rho_0^{\pm1},\qquad
 B_j^{\pm1}\to B_0^{\pm1}\quad\hbox{strongly in local }L^{P/2}.
\end{equation}
Indeed, uniform integrability at exponent $P/2$ follows from the displayed bounds and strong convergence. The conclusion then follows from continuity of the matrix expressions and Vitali's theorem.

We apply these estimates to sequences whose curvature energy tends to
zero on bounded balls. Local smoothing identifies the flat limit, and
the estimate at time zero recovers the original metrics in its charts.

\begin{lemma}\label{lem:flat-package}
Let $(M_j,h_j,p_j)$ be complete connected smooth four-manifolds with
fixed constants $S,\kappa,B>0$ and $0\le\lambda_j\le1$ such that
\begin{align}
 \norm{u}_{L^4(dV_{h_j})}^2
 &\le S\int(|du|_{h_j}^2+\lambda_j u^2)dV_{h_j},
 &&u\in C_c^\infty(M_j),\label{eq:flat-sobolev-input}\\
 \kappa r^4\le\Vol_{h_j}B(x,r)&\le Br^4,
 &&x\in M_j,\ 0<r\le1,\label{eq:flat-volume-input}\\
 \int_{B(p_j,R)}|\Rm(h_j)|^2dV_{h_j}&\longrightarrow0,
 &&R<\infty.\label{eq:vanishing-full-curvature}
\end{align}
A subsequence converges in the pointed measured Gromov--Hausdorff topology to a
complete connected smooth flat manifold $(N,h,p)$ with its Riemannian volume
measure. For every prescribed finite $P>2$ there are compatible local
maps $F_j$ with
\begin{equation}\label{eq:flat-strong-metric-convergence}
 F_j^*h_j\to h,\qquad (F_j^*h_j)^{-1}\to h^{-1}
 \quad\hbox{strongly in local }L^P.
\end{equation}
Compact subcharts cover balls of a fixed small radius in the original metrics. Densities
and Dirichlet matrices, including their inverses, converge in $L^{P/2}$.
The same conclusions hold if all the measures are multiplied by the same positive constant.
The local maps and subsequence may depend on $P$. On compact overlaps
the transition maps converge uniformly to the transitions of the
limiting atlas.
\end{lemma}

\begin{proof}
Fix a finite $P>2$. We work on bounded regions about $p_j$ and
use a diagonal subsequence at the end. All chart maps constructed
below are independent of time.

\paragraph{Step 1.}
Fix $0<\ell<1/10$ and take a maximal $\ell$-separated net in
each $M_j$, including $p_j$. Join two net points when their
distance is at most $3\ell$. The balls of radius $\ell/3$
about neighbors of one vertex are disjoint and contained in its
$4\ell$-ball. The two volume bounds in
\eqref{eq:flat-volume-input} therefore give a uniform bound for
the degree of this graph. By sampling a minimizing geodesic from $p_j$ to any net point in $B(p_j,R)$ and using maximality of the net, we obtain a graph path of at most $C(1+R/\ell)$ edges.
Hence the number of net points in each bounded ball is uniformly bounded. The
$\ell$-balls centered at net points in
$B(p_j,R+\ell)$ cover $B(p_j,R)$, giving
$\Vol_{h_j}B(p_j,R)\le C_R$.

For any $0<r\le1$, disjoint $r/2$-balls with centers in
$B(p_j,R)$ have volume at least $\kappa(r/2)^4$ each and lie
in $B(p_j,R+1)$. Hence their number is at most
$C_{R+1}/[\kappa(r/2)^4]$. These uniform covering bounds give
pointed Gromov--Hausdorff precompactness. After extraction, write
$(Z,d,p)$ for the limit. It is a proper complete connected length
space. By the local volume bounds, the measures also have a weak limit $\mu_Z$ along a subsequence in a common realization of this convergence.
We identify it in Step 7.

\paragraph{Step 2.}
For any sequence of centers $q_j$ at bounded distance from $p_j$, apply the local
smoothing theorem of Lee and Lee~\cite[Theorem 1.1]{LL} on a larger initial ball.
With $\bar r=0$, the theorem requires precompactness of the initial ball $B(q_j,3)$,
the upper volume bound for centers in
$B(q_j,2)$ at radii at most one, and the Sobolev and curvature
bounds on their unit balls. All the balls in these hypotheses lie in a single bounded
ball about $p_j$.
More explicitly, choose $R_0>0$ so that the fixed-radius ball about $p_j$ contains the
larger initial ball and all the required unit balls, and put
\[
 \epsilon_j=
 \left(\int_{B(p_j,R_0)}|\Rm(h_j)|^2dV_{h_j}\right)^{1/2}
 \longrightarrow0.
\]
Each required unit ball has curvature norm at most $\epsilon_j$.
The volume upper bound is \eqref{eq:flat-volume-input}. By extending test functions by zero and using $\lambda_j\le1$, we obtain the required local Sobolev inequality from \eqref{eq:flat-sobolev-input}.
Since the initial metric is complete, the larger initial ball is precompact. The construction
time $\widehat T$ and all constants in the smoothing theorem are
therefore independent of $j$ and of the chosen bounded sequence
of centers.
In the interior of each fixed initial domain, the resulting flows satisfy
\begin{equation}\label{eq:lee-lee-smoothing-estimates}
 \begin{gathered}
 |\Rm(h_j(t))|\le\Lambda(\epsilon_j+\sqrt t)t^{-1},\qquad
 \norm{\Rm(h_j(t))}_{L^2(B_{j,t}(x,1/2))}
       \le\Lambda(\epsilon_j+\sqrt t),\quad\epsilon_j\to0,\\
 \operatorname{inj}_{h_j(t)}\ge\iota\sqrt t,
 \qquad |\nabla^k\Rm(h_j(t))|\le C_kt^{-1-k/2}.
 \end{gathered}
\end{equation}
For the derivative bounds, apply the interior estimates of Shi~\cite{Shi} on $[t/2,t]$
in a slightly larger spatial region. See also the use of these
estimates in Lee and Lee~\cite[proof of Theorem 1.3]{LL}. Choose
$0<a<1/(P+2)$. After discarding finitely many terms and shortening the time interval to $0<t\le T\le\widehat T$, we have $|\Ric|_{\rm op}\le a/t$. Keep the original construction time
$\widehat T$ fixed during this shortening and set
$r_e=2L\sqrt{\widehat T}$, $\tau_e=L^2\widehat T$.
Theorem 1.1(4) of Lee and Lee~\cite{LL} gives
$\nu(B_{j,t}(q_j,r_e),h_j(t),\tau_e)\ge-\Lambda$.
Restriction to a smaller interior ball preserves this lower bound.
For $z$ compactly supported in $B_{j,t}(q_j,r)$, $r<\min(r_e,1/2)$,
Martens~\cite[Lemma 3.3(2)]{Martens} therefore gives
\begin{equation}\label{eq:positive-time-local-sobolev}
 \norm{z}_4^2\le S_0\int(|dz|^2+L_0z^2)dV_{h_j(t)}.
\end{equation}
Here $L_0=\tau_e^{-1}$ and $S_0$ are independent of $t$ and $j$.
The additional scalar term in the entropy-to-Sobolev estimate is
bounded by
$C\norm{\Rm}_{L^2(B_{j,t}(q_j,r))}\norm{z}_4^2$.
It is absorbed by making $\Lambda(\epsilon_j+\sqrt T)$ smaller
than the fixed absorption threshold. See also
Lee and Lee~\cite[proof of Lemma 2.1, equation (2.2)]{LL}.
Because $r_e$ and $\tau_e$ remain fixed, this shortening does not
change $S_0$ or $L_0$. We choose the smallness threshold to absorb
the cubic curvature term in Step 4 as well.

\paragraph{Step 3.}
Write $d_{j,t}$ for intrinsic distances in the flow domain. On a fixed initial
comparison neighborhood, the local smoothing theorem gives
\begin{equation}\label{eq:local-holder-distance-comparisons}
 d_{j,t}(x,y)\le Cd_{j,0}(x,y)^\alpha,
 \qquad d_{j,0}(x,y)\le Cd_{j,t}(x,y)^\beta
\end{equation}
at small distances, with fixed $C$ and $\alpha,\beta\in(0,1)$
for all sufficiently large $j$. These follow by weakening the exponents
in Lee and Lee~\cite[Theorem 1.1]{LL} uniformly as $\epsilon_j\to0$.
For central points, original minimizing geodesics
stay inside the larger initial ball, so $d_{j,0}$ is the original distance.
Choose $\overline B_{j,0}(q_j,2r_b)$ inside the comparison neighborhood.
The shrinking-balls lemma of Lee and Lee~\cite[Lemma 3.2]{LL} gives
$B_{j,t}(q_j,2r_b-c_b\sqrt t)\subset B_{j,0}(q_j,2r_b)$.
Reduce $T$ so that $c_b\sqrt T<r_b/2$. Thus positive-time balls of
radius $r_b$ lie in the region where both comparisons are valid.

Choose fixed radii $r_{\rm in}<r_{\rm out}/4<r_b/8$, with
$r_{\rm out}<\min(r_e,1/2)$. Take
$r_{\rm work}$ with $Cr_{\rm work}^\alpha<r_{\rm in}/4$.
At time $T$ choose normal-coordinate maps $F_j:B(0,2\rho)\to M_j$,
where $2\rho$ is below the uniform injectivity radius and $r_b/4$,
and $C(2\rho)^\beta<r_{\rm work}$. By the shrinking-balls lemma, their images lie in
the comparison neighborhood. Then \eqref{eq:local-holder-distance-comparisons}
gives, for every $0<t\le T$,
\begin{equation}\label{eq:fixed-chart-moving-ball-containment}
 F_j(B(0,2\rho))\subset B_{j,0}(q_j,r_{\rm work})
 \subset B_{j,t}(q_j,r_{\rm in}/4).
\end{equation}
Conversely, choose $s_0$ with $Cs_0^\alpha<\rho/2$. Then
\begin{equation}\label{eq:original-ball-chart-coverage}
 B_{j,0}(q_j,s_0)\subset F_j(B(0,\rho/2)).
\end{equation}
By the positive-time curvature and injectivity bounds, we obtain smooth convergence in these normal charts at $T$ after choosing orthonormal frames and passing to a subsequence. In particular $H_j(T)$ and its
inverse are uniformly bounded in every local coordinate norm.
On each interval $[\tau,T]$, we first integrate $|\Ric|\le a/t$ to bound $H_j(t)$ and
its inverse relative to $H_j(T)$. Bounds for the Christoffel symbols in these
coordinates then follow from the evolution equation for the connection and the
estimates for the covariant derivatives of curvature.
By differentiating these equations and using metric compatibility, we obtain all higher coordinate derivative bounds by induction.
Along a diagonal subsequence, the pullbacks satisfy $H_j(t)=F_j^*h_j(t)\to H(t)$
smoothly for $t>0$ on smaller fixed charts.

\paragraph{Step 4.}
Put $f_j=|\Rm(h_j(t))|$. For fixed nested radii $r<R$ between
$r_{\rm in}$ and $r_{\rm out}$, the shifted distance cutoffs of
Lee and Lee~\cite[Lemma 2.1]{LL} have the form
$\chi=e^{-C_{r,R}t}\varphi_{r,R}(d_{j,t}(q_j,\cdot)+b\sqrt t)$,
with $(\partial_t-\Delta)\chi\le0$ and
$|d\chi|^2\le C_{r,R}\chi$ in the barrier sense. Their supports
are precompact in the above entropy region.

For the cutoff inequality, apply the distance comparison of
Perelman~\cite[Lemma 8.3(a)]{Perelman} at scale $\sqrt t$. The bound
$\Ric\le a/t$ gives
\[
 (\partial_t-\Delta)d_{j,t}(q_j,\cdot)\ge-C/\sqrt t
\]
on the transition annulus, after time has been shortened. Choosing
$b\ge2C$ makes $\eta=d_{j,t}+b\sqrt t$ satisfy
$(\partial_t-\Delta)\eta\ge0$. Choose a nonincreasing profile
$\varphi_{r,R}$ with $\varphi''\ge-C_{r,R}\varphi$ and
$|\varphi'|^2\le C_{r,R}\varphi$, equal to one before $r$ and
zero after $R$. The asserted cutoff inequalities follow from the chain rule and the exponential factor. The transition stays away
from the center. Distance barriers and approximation handle the
cut locus. 

Combining the curvature evolution equation, Kato's inequality, and $\partial_tdV=-\operatorname{Scal}\,dV$, we obtain
\[
 \frac d{dt}\int\chi^2f_j^2+c_0\int|d(\chi f_j)|^2
 \le C\int f_j(\chi f_j)^2+C_{r,R}\int_{\operatorname{supp}\chi}f_j^2.
\]
Distributional approximation justifies the calculation at the cut locus.
To check the absorption, put $w=\chi f_j$ and use
\[
 \int f_jw^2
 \le\norm{f_j}_{L^2(\operatorname{supp}\chi)}\norm{w}_4^2
 \le S_0\Lambda(\epsilon_j+\sqrt t)
            \int(|dw|^2+L_0w^2)dV_{h_j(t)}.
\]
The smallness condition in Step 2 allows us to absorb the gradient term. The
remaining $L_0$ term is bounded by the outer curvature energy.
Define
$E_j(r,t)=\int_{d_{j,t}+b\sqrt t<r}f_j^2dV_{h_j(t)}$ and
$E_j(r,0)=\int_{B_{j,0}(q_j,r)}|\Rm(h_j)|^2dV_{h_j}$.
Since $\chi$ is supported in the outer region and equals $e^{-C_{r,R}t}$ on the inner
region, integration gives
\begin{equation}\label{eq:nested-curvature-energy}
 E_j(r,t)\le C_{r,R}\left(E_j(R,0)+\int_0^tE_j(R,s)\,ds\right).
\end{equation}
We first integrate from a positive lower time and then let that time tend to zero for each $j$, using smoothness of the individual flow at zero.

Fix an integer $N\ge1$ and choose
$r_{\rm in}=r_0<r_1<\cdots<r_N<r_{\rm out}$. Set
\[
 \delta_j^2=\int_{B_{j,0}(q_j,r_{\rm out})}
                    |\Rm(h_j)|^2dV_{h_j}\longrightarrow0.
\]
All initial terms in the $N$ applications of
\eqref{eq:nested-curvature-energy} are bounded by $\delta_j^2$,
and the outermost positive-time energy is bounded by a uniform
constant $K$. Repeated time integration produces a sum of terms
$C_N\delta_j^2t^k/k!$, $0\le k<N$, and a remainder at most
$C_NKt^N/N!$. Consequently we have
\begin{equation}\label{eq:iterated-curvature-energy}
 E_j(r_{\rm in},t)\le C_N\delta_j^2+C_Nt^N
 \qquad(0<t\le T_N).
\end{equation}
Here $T_N\le T$ may decrease to retain every cutoff margin.
Decrease $T_N$ so that $b\sqrt{T_N}<r_{\rm in}/4$.
Then \eqref{eq:fixed-chart-moving-ball-containment} places the
fixed chart images strictly inside the energy region.

For a point in a smaller fixed chart, take a parabolic cylinder
with spatial radius $c\sqrt t$ and time interval $[t/2,t]$.
The cylinder lies in the energy region for all sufficiently
small $t$. On this interval the metrics are uniformly comparable,
and $f_j$ satisfies
$(\partial_s-\Delta)f_j\le C t^{-1}f_j$ in the weak sense.
When the metrics are rescaled by $t^{-1}$, the Sobolev inequality
\eqref{eq:positive-time-local-sobolev} has a bounded lower-order
coefficient $tL_0$. Parabolic Moser iteration, as in
Lee and Lee~\cite[proof of Proposition 2.15]{LL}, therefore yields
\begin{equation}\label{eq:curvature-moser-from-energy}
 \sup_{\rm inner}|\Rm(h_j(t))|
 \le C_Nt^{-1}\sup_{t/2\le s\le t}E_j(r_{\rm in},s)^{1/2}
 \le C_N(\delta_j/t+t^{N/2-1}).
\end{equation}
The factor $t^{-1}$ is the four-dimensional scaling factor for
a spatial $L^2$ curvature norm. Equivalently, the parabolic
$L^2$ estimate is $Ct^{-3/2}$ times the spacetime norm, whose
square is at most $t\sup_s E_j(r_{\rm in},s)$.

\paragraph{Step 5.}
First let $j\to\infty$ at a fixed positive time. Applying interior derivative estimates on a slightly larger cylinder, we obtain
\[
 |\nabla^k\Rm(H(t))|
 \le C_{k,N}t^{N/2-1-k/2}.
\]
Given $k,m$, choose $N\ge2m+2+k$. Thus every covariant
curvature derivative is $O(t^m)$ on a smaller fixed chart, for
sufficiently small time depending on $k,m$.

In particular $\int_0^{t_0}\sup|\Ric(H(t))|\,dt<\infty$.
By integrating $\partial_tH=-2\Ric$, we obtain upper and lower bounds for $H(t)$ in terms of the positive metric $H(t_0)$, uniformly down to zero.
The right side of $\partial_t\Gamma=-H^{-1}*\nabla\Ric$ is now integrable in time in
these coordinates. Hence $\Gamma$ remains bounded down to time zero. Metric
compatibility, $\partial H=H*\Gamma$, bounds the
first spatial derivatives of $H$. All spatial derivatives are controlled by differentiating the connection equation and repeating this argument. The corresponding time derivatives are integrable.
Hence $H(t)$ converges in every local $C^k$ norm to a smooth
positive metric $h_0$. Curvature convergence gives $\Rm(h_0)=0$.
For each derivative order, $N$ is fixed before $t\downarrow0$.

\paragraph{Step 6.}
It remains to identify the distance, which does not follow from bi-H\"older comparison
alone. We first show that the curves under consideration stay in the chart. Choose
$\eta,s_1>0$ with
\[
 4\eta<\min(r_b,1/2),\qquad
 C(4\eta)^\beta<s_0,\qquad Cs_1^\alpha<\eta.
\]
Endpoints in $B_{j,0}(q_j,s_1)$ lie in $B_{j,t}(q_j,\eta)$.
Their short minimizing sequences, and hence their minimizing
curves, stay in $B_{j,t}(q_j,3\eta)$. Sampling balls of sufficiently
small radius lie in $B_{j,t}(q_j,4\eta)$. By the shrinking-balls lemma, this last ball
lies in the comparison neighborhood. The reverse
comparison then places it in $B_{j,0}(q_j,s_0)$, which is covered
by \eqref{eq:original-ball-chart-coverage}. Its closure is contained
in a larger compact subchart. This precompact containment also
justifies existence of the minimizers in the incomplete flow domain.

Fix two such endpoints and let $\sigma$ be any minimizing
$h_j(t)$-geodesic between them. Its length is at most $2\eta$. For $r\le c\sqrt t$, it follows from the curvature and injectivity estimates that the sampling balls have volume at least $c_1r^4$, and $|\nabla\Ric|\le Ct^{-3/2}$ on a neighborhood
of them. Partition $\sigma$ into segments of length $r$, with
at most one shorter final segment.
For a segment $I$, let $B$ be its sampling ball of radius comparable to $r$. Averaging
the Lipschitz function $|\Ric|$ over $B$, we obtain
\[
 \int_I|\Ric|\,ds
 \le Cr\,\Vol(B)^{-1/2}\norm{\Ric}_{L^2(B)}
       +Cr^2t^{-3/2}
 \le Cr^{-1}\norm{\Ric}_{L^2(B)}+Cr^2t^{-3/2}.
\]
There are at most $C/r$ segments. Since $\sigma$ is minimizing, the distance between
any two of its points equals the length of the intervening segment. The
sampling balls can therefore be split into a fixed number of
disjoint families. Summing and applying the Cauchy--Schwarz inequality, we obtain
\begin{equation}\label{eq:geodesic-ricci-averaging}
 \int_\sigma|\Ric|\,ds
 \le C\left(r^{-3/2}\norm{\Ric}_{L^2(\rm region)}
                +rt^{-3/2}\right).
\end{equation}
The regional $L^2$ norm is uniformly bounded by
\eqref{eq:lee-lee-smoothing-estimates}. For sufficiently small
$t$, choose $r=ct^{3/5}\le c\sqrt t$. Both terms are then
bounded by $Ct^{-9/10}$. The first variation formula for length gives $|\partial_t d_{j,t}(x,y)|\le Ct^{-9/10}$ for almost every
positive time (or the corresponding upper and lower Dini bounds).
Integrating between positive times and using smoothness at zero for each fixed $j$, we obtain the uniform estimate
\begin{equation}\label{eq:uniform-distance-time-modulus}
 |d_{j,t}(x,y)-d_{j,0}(x,y)|\le Ct^{1/10}.
\end{equation}
Any further shortening of time here leaves the maps $F_j$ unchanged.

Restrict to a smaller central subchart whose image lies in
$B_{j,0}(q_j,s_1)$. By the two distance comparisons, the chart images still cover a fixed smaller original ball. If $\Phi_j$ are
pointed Gromov--Hausdorff approximation maps, it follows from the comparisons and the normal-coordinate control at time $T$ that there are constants
$c_1,C_1>0$ such that, on compact smaller subcharts,
\[
 c_1|x-y|^{1/\alpha}
 \le d_{j,0}(F_j(x),F_j(y))
 \le C_1|x-y|^\beta.
\]
Indeed, $d_{j,T}$ is uniformly comparable with coordinate
distance there, and both inequalities follow by applying
\eqref{eq:local-holder-distance-comparisons} at $T$.
Thus $\Phi_jF_j$ are equicontinuous up to errors tending to
zero and separate distinct coordinate points. After extraction,
these maps converge locally
uniformly to an injective map $F$ into $Z$. The ball coverage
shows that its image contains a neighborhood of the limit of
$q_j$. Approximate a point in a smaller limit ball by points in
the corresponding original balls and take a convergent subsequence
of their coordinates in a compact subchart.

For central endpoints, all minimizing curves lie in the larger
represented chart. By smooth convergence at fixed $t$, their distances therefore
converge to the corresponding distances for $H(t)$. Now use
\eqref{eq:uniform-distance-time-modulus} and let $t\downarrow0$.
Since $H(t)\to h_0$ smoothly, this proves
$d(F(x),F(y))=d_{h_0}(x,y)$ for sufficiently close coordinate
points. Thus $F$ is a local isometry onto a neighborhood.
Repeating this argument at every limit point, we obtain flat charts. Their transition
maps are local distance isometries of smooth flat metrics, hence
are smooth. The prelimit transitions converge locally uniformly
on compact overlaps by the same equicontinuity and inverse
equicontinuity. These are the compatible charts in the statement.
We can now write $(N,h,p)$ for $(Z,d,p)$ with this smooth flat
structure. Completeness and connectedness follow from Step 1.

\paragraph{Step 7.}
Recall from Step 2 that $a<1/(P+2)$. By
\eqref{eq:fixed-chart-moving-ball-containment}, the uniform $L^2$ bound on the Ricci
tensor holds on a neighborhood of each retained chart, as required by
Lemma~\ref{lem:zero-time-coefficients}. That
lemma gives strong local $L^P$ convergence of $F_j^*h_j$ and
its inverse. The limiting metric is $h_0$, since the positive-time
limit $H(t)$ converges smoothly to $h_0$. This proves
\eqref{eq:flat-strong-metric-convergence}. The determinant
estimates \eqref{eq:strong-density-coefficients} then give strong
local $L^{P/2}$ convergence of the densities and Dirichlet
matrices, including their inverses. In particular, because $P/2>1$,
the densities converge in $L^1$ on compact subcharts.

For a chart domain $U$, $\chi\in C_c(U)$, and a bounded continuous
function $v$ on $N$, this yields
\[
 \int\chi v(\Phi_jF_j)\rho_jdx\longrightarrow
 \int\chi v(F)\rho_0dx.
\]
Indeed, the density error tends to zero in $L^1$, and
$v(\Phi_jF_j)\to v(F)$ uniformly on the support of $\chi$.
Here $c=1$ and $\rho_0dx=dV_{h_0}$. To identify the measure locally,
take a test function supported in a smaller limit ball. By
\eqref{eq:original-ball-chart-coverage}, all points of the approximating
spaces mapping into its support lie in a compact subchart for large $j$.
We may therefore take $\chi=1$ there. This identifies $\mu_Z$ with
$dV_h$ on the smaller ball. A finite partition of unity gives the same
identity on every compact region. If all measures are multiplied by the same positive
constant, the limiting measure acquires the same factor.

Finally, every distance sphere in a complete smooth Riemannian
manifold has volume zero. It is contained in the image under the
exponential map of the radius-$R$ sphere in a tangent space.
By weak convergence of the measures, the volumes of fixed balls converge,
even when their centers vary along a convergent sequence. Taking a diagonal subsequence over
an exhaustion by bounded regions and a countable family of charts, we obtain both the
pointed measured convergence and the asserted local convergences.
\end{proof}

To pass to the limit in the equation for a bounded function on these charts, it remains
to establish compactness of the functions and control the Dirichlet form as its
coefficients vary.

\begin{lemma}
\label{lem:bounded-function-compactness}
On fixed charts, let $\rho_j,B_j$ be smooth positive fields, with $B_j$
symmetric. Suppose $\rho_j\to\rho_0>0$ strongly in local $L^1$ and
$\norm{B_j^{-1}}_{L^r}$ is locally bounded for some $r>1$.
If $w_j\in W^{1,1}_{\rm loc}$, $|w_j|\le M$, and
$\int_K\nabla w_j^TB_j\nabla w_j\le C_K$, then, after passing to a subsequence, we have
\begin{equation}\label{eq:bounded-function-strong-convergence}
 w_j\to w\ \hbox{a.e. and in every finite local }L^b(dx),\qquad
 f(w_j)\rho_j\to f(w)\rho_0\ \hbox{in local }L^1
\end{equation}
for continuous $f$ on $[-M,M]$. If additionally $B_j^{\pm1}\to
B_0^{\pm1}$ in local $L^1$, with $B_0,\rho_0$ smooth positive, and
\begin{equation}\label{eq:approximating-semilinear-equation}
 -\operatorname{div}(B_j\nabla w_j)=\rho_j f_j(w_j),
 \qquad f_j\to f\ \hbox{uniformly on }[-M,M],
\end{equation}
then $w\in H^1_{\rm loc}$ solves
\begin{equation}\label{eq:limiting-semilinear-equation}
 -\operatorname{div}(B_0\nabla w)=\rho_0f(w),
\end{equation}
and
\begin{equation}\label{eq:local-energy-lower-semicontinuity}
 \int\chi\nabla w^TB_0\nabla w
 \le\liminf_j\int\chi\nabla w_j^TB_j\nabla w_j
 \quad(0\le\chi\in C_c^\infty).
\end{equation}
\end{lemma}

\begin{proof}
Set $q=2r/(r+1)>1$.

\paragraph{Step 1.}
By the weighted H\"older inequality,
\begin{equation}\label{eq:weighted-gradient-rellich}
 \int_K|\nabla w_j|^q
 \le\left(\int_K\nabla w_j^TB_j\nabla w_j\right)^{q/2}
       \left(\int_K|B_j^{-1}|^r\right)^{(2-q)/2}.
\end{equation}
By Rellich's theorem on a slightly larger domain, a subsequence converges almost everywhere.
Since the functions are uniformly bounded, they converge strongly in every finite
local Lebesgue norm. Split the nonlinear
density difference as $f(w_j)(\rho_j-\rho_0)+(f(w_j)-f(w))\rho_0$.
By strong $L^1$ convergence and the dominated convergence theorem, we obtain
\eqref{eq:bounded-function-strong-convergence}, also for $f_j$.
In particular, for the critical fourth-power mass on a compact chart,
\[
 \int|w_j^4\rho_j-w^4\rho_0|\,dx
 \le M^4\|\rho_j-\rho_0\|_{L^1}
       +\int|w_j^4-w^4|\rho_0\,dx\longrightarrow0.
\]

\paragraph{Step 2.}
Under the additional assumptions on the coefficients, we may test the dual quadratic form against any smooth vector field $\zeta$ to obtain
\[
 \int\chi\nabla w_j^TB_j\nabla w_j
 \ge2\int\chi\nabla w_j\cdot\zeta-\int\chi\zeta^TB_j^{-1}\zeta.
\]
Pass to the limit distributionally and optimize in $\zeta$ to obtain
\eqref{eq:local-energy-lower-semicontinuity} and $w\in H^1_{\rm loc}$.
To pass to the equation, consider the positive semidefinite matrix
\begin{equation}\label{eq:coefficient-duality-identity}
 Q_j=(B_j-B_0)B_j^{-1}(B_j-B_0)
     =B_j-2B_0+B_0B_j^{-1}B_0\longrightarrow0\quad\hbox{in }L^1.
\end{equation}
By Cauchy--Schwarz, the square of the flux error against a smooth test function $\phi$ is bounded by
\[
 \left(\int_{\operatorname{supp}\phi}\nabla w_j^TB_j\nabla w_j\right)
 \int\nabla\phi^TQ_j\nabla\phi\longrightarrow0.
\]
After integration by parts, the remaining flux term involves $w_j$ paired with the smooth function $\operatorname{div}(B_0\nabla\phi)$. Convergence on the right follows from
\eqref{eq:bounded-function-strong-convergence}.
\end{proof}

On the charts of Lemma~\ref{lem:flat-package}, we may take $P=4$,
$r=2$, and $q=4/3$. Write $d\mu_j=c\,dV_{h_j}$ and $d\mu=c\,dV_h$.
By \eqref{eq:bounded-function-strong-convergence} and the chart
convergence, the measures $w_jd\mu_j$ converge locally to $w\,d\mu$.
Their weak limits agree on overlaps, so the local limits define a
single function $w$ on $N$.

Fix a compact set $K\subset N$. Choose finitely many nonnegative smooth
cutoffs $\chi_a$, each compactly supported in a limit chart, with
$\sum_a\chi_a\le1$ on $N$ and equality near $K$. Transport each cutoff
to $M_j$ using the corresponding chart and extend it by zero. Denote
the resulting function by $\chi_{a,j}$. On its support,
$\chi_{a,j}=\chi_a\circ\Phi_j+o(1)$ uniformly, since the chart maps
converge locally uniformly after composition with the pointed
approximations $\Phi_j$. As all cutoffs are nonnegative, it follows
that $\sum_a\chi_{a,j}\le1+o(1)$ on $M_j$.
Applying \eqref{eq:local-energy-lower-semicontinuity} in each chart
and summing, we obtain
\[
 \gamma\int_K |dw|_h^2\,d\mu
 \le \liminf_j\sum_a\gamma\int_{M_j}\chi_{a,j}|dw_j|_{h_j}^2\,d\mu_j
 \le \liminf_j D_{h_j}(w_j).
\]
Exhausting $N$ by compact sets gives
\begin{equation}\label{eq:global-dirichlet-liminf}
 D_h(w)\le\liminf_j D_{h_j}(w_j).
\end{equation}
Applying the same cutoffs to the fourth-power densities in
\eqref{eq:bounded-function-strong-convergence}, we obtain
\[
 \int_N |w|^4\,d\mu\le\liminf_j\int_{M_j}|w_j|^4\,d\mu_j.
\]

\label{sec:gap}

Fix a Sobolev level $A$ above the value at the reference metric.
To control the Sobolev constant throughout a connected low-energy
sublevel, it suffices to exclude metrics at level $A$ from that sublevel.
The compactness estimates give the required energy gap.

\begin{theorem}\label{thm:barrier}
For each $A>\CS(\omega_*)$, there is $\delta_A>0$ such that every smooth
K\"ahler metric $\omega\in[\omega_*]$ with $\CS(\omega)=A$ satisfies
$\Ca(\omega)\ge\delta_A$.
\end{theorem}

\begin{proof}
Suppose, to the contrary, that $\CS(g_j)=A$ and $\Ca(g_j)\to0$.
By \eqref{eq:euclidean-threshold},
$A>\CS(g_*)\ge\CE$. We use the same normalizations
$d\mu_j=c\,dV_{g_j}$ and
$D_j(w)=\gamma\int|dw|_{g_j}^2d\mu_j$. The operators
$-\Delta_j$ and $-\Delta_{h_j}$ below have this Dirichlet
normalization. Here $\Delta_{\mathrm E}$ denotes the ordinary Euclidean Laplacian, so
the limiting operator is
$-\gamma\Delta_{\mathrm E}$.

\paragraph{Step 1.}
By Propositions~\ref{prop:strong} and~\ref{prop:nc}, we may pull back by holomorphic
automorphisms and pass to a subsequence. The metric coefficients, inverses, and
densities then converge strongly to those of the background metric in every fixed
finite background $L^p$ space. Moreover,
\begin{equation}\label{eq:gap-nc}
 \int_\Omega|\Rm(g_j)|^2dV_j
 \le\eta_j+C\Vol_{g_j}(\Omega)^{1/2},\qquad\eta_j\to0.
\end{equation}
This holds for every measurable $\Omega\subset X$, including
balls for $g_j$. Pullback preserves $\CS$ and $\Ca$. In
particular, the critical Sobolev bound used below is supplied by
$\CS(g_j)=A$ itself.

Write $v_j(x,r)=\Vol_{g_j}B(x,r)$. Apply the Sobolev inequality
to a Lipschitz distance cutoff supported in $B(x,r)$, equal to
one on $B(x,r/2)$, with gradient at most $2/r$. Smooth
approximation makes it an admissible test. For $0<r\le1$,
\[
 v_j(x,r/2)^{1/2}\le Cr^{-2}v_j(x,r),
\]
where $C$ depends only on $A,c,\gamma$. Thus
$q_j(x,r)=r^{-4}v_j(x,r)$ satisfies
$q_j(x,r)\ge c_0q_j(x,r/2)^{1/2}$. Iterating $k$ times, we obtain
\[
 q_j(x,r)\ge c_0^{\,2(1-2^{-k})}
                  q_j(x,2^{-k}r)^{2^{-k}}.
\]
For each fixed $j,x$, smoothness gives
$q_j(x,2^{-k}r)\to\beta_4$, the Euclidean unit-ball volume.
By letting $k\to\infty$, we obtain the uniform lower bound
\begin{equation}\label{eq:gap-lower-volume}
 \Vol_{g_j}B(x,r)\ge\kappa r^4\quad(0<r\le r_0)
\end{equation}
with $r_0\le1$, which we may decrease below. For this argument, the small-ball volume
ratios need only converge for each fixed metric.

\paragraph{Step 2.}
We next prove a matching upper bound, uniform over all sufficiently small radii. Fix
$B>\beta_4$. If no such interval exists, choose,
after taking a subsequence, radii tending to zero where
$\max_xq_j(x,r)>B$. On each fixed smooth compact manifold,
$q_j(x,r)\to\beta_4$ uniformly in $x$ as $r\downarrow0$,
and the maximum over $x$ is continuous in $r$. There is therefore
a first radius $\ell_j\downarrow0$ at which this maximum is
$B$. Choose a maximizing center $p_j$ and set
$h_j=\ell_j^{-2}g_j$. By the choice of this first radius and
\eqref{eq:gap-lower-volume}, we have
\[
 \kappa r^4\le\Vol_{h_j}B(x,r)\le Br^4\quad(0<r\le1),
 \qquad
 \Vol_{h_j}B(p_j,1)=B.
\]
It follows from the separated-net argument in Step 1 of Lemma~\ref{lem:flat-package} that
$\Vol_{h_j}B(p_j,R)\le C_R$ for every fixed $R$. The unit-scale
volume bounds control the degree of the net graph, and minimizing
geodesics control the number of edges needed to reach a bounded
ball. In real dimension four, the integral of $|\Rm|^2$ is
unchanged by constant rescaling. Hence \eqref{eq:gap-nc} gives
\[
 \int_{B_{h_j}(p_j,R)}|\Rm(h_j)|^2dV_{h_j}
 \le\eta_j+CC_R^{1/2}\ell_j^2\longrightarrow0.
\]
With $d\nu_j=\ell_j^{-4}d\mu_j=c\,dV_{h_j}$, the rescaled
Sobolev inequality is
\[
 \left(\int|w|^4d\nu_j\right)^{1/2}
 \le A\left(D_{h_j}(w)+\ell_j^2\int w^2d\nu_j\right).
\]
Up to the fixed normalization factors, these estimates verify the hypotheses of
Lemma~\ref{lem:flat-package}. It gives a
complete smooth flat limit with its Riemannian volume measure.
Since the volumes of the unit balls converge, the limit unit ball has volume $B$.
Bishop--Gromov comparison gives volume at most $\beta_4$
for every unit ball in a complete flat four-manifold, a
contradiction. Indeed, the exponential map is a surjective local isometry from
the Euclidean unit tangent ball onto the metric unit ball. Thus, after decreasing $r_0$, we have
\[
 \kappa r^4\le\Vol_{g_j}B(x,r)\le Br^4
 \qquad(x\in X,\ 0<r\le r_0).
\]

\paragraph{Step 3.}
For any $r_j\downarrow0$ and centers $q_j$, set
$\widetilde h_j=r_j^{-2}g_j$. For each fixed $R$, the original
radius $r_jR$ is eventually at most $r_0$. The preceding bounds
give
\[
 \kappa R^4\le\Vol_{\widetilde h_j}B(q_j,R)\le BR^4,
 \qquad
 \int_{B_{\widetilde h_j}(q_j,R)}|\Rm(\widetilde h_j)|^2
             dV_{\widetilde h_j}
 \le\eta_j+CB^{1/2}r_j^2R^2\longrightarrow0.
\]
The rescaled Sobolev lower-order coefficient is $r_j^2$.
Lemma~\ref{lem:flat-package} gives a complete connected smooth
flat limit $(N,h,q)$, with its natural measure and compatible
charts, and
$\Vol_hB(q,R)\ge\kappa R^4$ for every $R>0$.

The quartic lower bound rules out every nontrivial flat quotient.
Indeed, the universal Riemannian cover of $N$ is $\mathbb R^4$.
If its deck group contains a nonidentity isometry
$T(x)=Ox+b$, freeness implies that the projection of $b$ onto
$\ker(I-O)$ is nonzero. Otherwise $(I-O)x=b$ has a solution.
After translating the origin and choosing an axis in that
projection direction, we can write
\[
 T(s,y)=(s+L,O'y),\qquad L>0,\quad y\in\mathbb R^3.
\]
The cyclic quotient by $\langle T\rangle$ has a fundamental
slab $0\le s<L$. A radius-$R$ ball in this quotient is contained
in the portion of this slab with $|y|\le R+C$, where $C$ depends
on its center, and so has volume $O(R^3)$. The covering onto
$N$ maps a ball about a chosen lift onto the corresponding ball
in $N$ and cannot increase its volume. This contradicts the quartic lower bound for
volume growth. The deck group is trivial, so every such
subsequential limit is $(\mathbb R^4,g_{\mathrm E})$. The
measures $r_j^{-4}d\mu_j$ converge to $c\,dx$.

\paragraph{Step 4.}
We justify attainment of the Sobolev quotient on each fixed
$(X,g_j)$. In this step constants may depend on $j$. For every $\varepsilon>0$, the
Euclidean inequality on sufficiently small coordinate balls gives
\begin{equation}\label{eq:fixed-metric-near-sharp}
 \left(\int|w|^4d\mu_j\right)^{1/2}
 \le(\CE+\varepsilon)D_j(w)
       +K_{j,\varepsilon}\int w^2d\mu_j.
\end{equation}
Indeed, choose a finite partition with $\sum_a\chi_a^2=1$
subordinate to balls on which the metric and density are as close
to their Euclidean values as needed. The triangle inequality in
$L^2$ gives $\|w\|_4^2\le\sum_a\|\chi_aw\|_4^2$.
Apply the local Euclidean inequality and use
\[
 \sum_aD_j(\chi_aw)
 =D_j(w)+\gamma\int\sum_a|d\chi_a|_{g_j}^2w^2d\mu_j
\]
to obtain \eqref{eq:fixed-metric-near-sharp}.

Put $Q_j(w)=D_j(w)+\int w^2d\mu_j$, and take a nonnegative
minimizing sequence $w_k$ with
$\int w_k^4d\mu_j=1$ and $Q_j(w_k)\to A^{-1}$.
By weak $H^1$ compactness and Rellich's theorem, we can pass to a subsequence such that
$w_k\rightharpoonup w$ in $H^1$ and the sequence converges to the same limit strongly
in $L^2$ and almost everywhere. Write $z_k=w_k-w$ and $\theta=\int w^4d\mu_j$.
By the lemma of Br\'ezis and Lieb~\cite[Theorem 1]{BrezisLieb} and the quadratic energy splitting, we have
\[
 \int|z_k|^4d\mu_j\longrightarrow1-\theta,\qquad
 Q_j(w_k)=Q_j(w)+Q_j(z_k)+o(1).
\]
Choose $\varepsilon>0$ so that $\CE+\varepsilon<A$.
Apply \eqref{eq:fixed-metric-near-sharp} to $z_k$, using its
strong $L^2$ convergence to zero, and use the definition of $A$
on $w$. We obtain
\[
 \frac1A\ge\frac{\sqrt\theta}{A}
                 +\frac{\sqrt{1-\theta}}{\CE+\varepsilon}.
\]
For $0\le\theta<1$ the right side is strictly greater than
$A^{-1}$. Thus $\theta=1$, $Q_j(w)=A^{-1}$, and the infimum
is attained. The Euler--Lagrange equation has multiplier
$A^{-1}$, as is seen by testing it with $w$.

To prove regularity at the critical exponent, set $w_L=\min(w,L)$ and test the equation
with $w w_L^2$.
Combining the resulting energy estimate with the fixed-metric Sobolev inequality and putting $W_L=w w_L$, we obtain
\[
 \|W_L\|_4^2\le C_j\int(1+w^2)W_L^2d\mu_j.
\]
Split $w^2$ at a level $K$. Since $w^2\in L^2$, choose $K$
so that $\|w^2\|_{L^2(\{w^2>K\},\mu_j)}$ is
small enough to absorb its contribution. The remaining integral
is bounded independently of $L$, because $W_L^2\le w^4$.
Letting $L\to\infty$, we find that $w\in L^8$. By elliptic $W^{2,8/3}$ regularity, $w$ is bounded and H\"older continuous. Further elliptic estimates then give smoothness. See also
Gilbarg and Trudinger~\cite{GilbargTrudinger}. The strong maximum principle gives
positivity. Denote this extremizer by $u_j$. It satisfies
\begin{equation}\label{eq:extremizer}
 \int u_j^4d\mu_j=1,\qquad
 -\Delta_ju_j+u_j=A^{-1}u_j^3,\qquad
 D_j(u_j)+\int u_j^2d\mu_j=A^{-1}.
\end{equation}

\paragraph{Step 5.}
We claim that $M_j:=\max_Xu_j\to\infty$. If the claim fails, there is a subsequence on
which $M_j$ is bounded. Proposition~\ref{prop:strong}
gives strong convergence of the densities and Dirichlet
coefficients, including their inverses, in fixed background charts.
Lemma~\ref{lem:bounded-function-compactness}, applied in a finite
chart cover, then gives a limit $u\in H^1(X,g_*)$ with
\[
 \int u^4d\mu_*=1,\qquad
 D_*(u)+\int u^2d\mu_*\le A^{-1}.
\]
Both mass terms converge because $u_j$ is uniformly bounded and the densities converge strongly. Summing the local lower bounds for the Dirichlet
energies against a partition of unity, we obtain the global energy
inequality. The Sobolev quotient of
$u$ is at least $A$, contradicting $\CS(g_*)<A$.

Choose a maximum point $p_j$ and put
\begin{equation}\label{eq:maximum-rescaling}
 h_j=M_j^2g_j,\qquad v_j=M_j^{-1}u_j,
 \qquad d\nu_j=M_j^4d\mu_j=c\,dV_{h_j}.
\end{equation}
Four-dimensional scaling gives
\[
 0<v_j\le1,\quad v_j(p_j)=1,\quad
 \int v_j^4d\nu_j=1,\quad D_{h_j}(v_j)=D_j(u_j).
\]
Put $a_j=M_j^{-2}\to0$. In these variables the equation and energy identity become
\begin{equation}\label{eq:rescaled-equation}
 -\Delta_{h_j}v_j+a_jv_j=A^{-1}v_j^3,
 \qquad D_{h_j}(v_j)+a_j\int v_j^2d\nu_j=A^{-1}.
\end{equation}
The rescaled Sobolev inequality has constant $A$ and lower-order
coefficient $a_j$. Step 3, with $r_j=M_j^{-1}$, gives the
pointed measured limit $(\mathbb R^4,g_{\mathrm E},0,c\,dx)$
and the coefficient convergence of Lemma~\ref{lem:flat-package}.

\paragraph{Step 6.}
To ensure that the limit is nonzero, we prove a uniform lower bound for the local mass. Since $0<v_j\le1$,
\eqref{eq:rescaled-equation} implies
$-\Delta_{h_j}v_j\le A^{-1}v_j$. Testing this inequality with $\eta^2v_j^{p-1}$, $p\ge2$, and using Young's inequality, we obtain
\[
 D_{h_j}(\eta v_j^{p/2})
 \le Cp^2\int(\eta^2+|d\eta|_{h_j}^2)v_j^p d\nu_j.
\]
For concentric balls with radii $s<t\le R$, take a distance
cutoff equal to one on the inner ball and with gradient at most
$2/(t-s)$. Combining the preceding estimate with the rescaled Sobolev inequality and $a_j\le1$, we obtain
\[
 \|v_j\|_{L^{2p}(B(p_j,s),\nu_j)}
 \le\bigl[Cp^2(1+(t-s)^{-2})\bigr]^{1/p}
           \|v_j\|_{L^p(B(p_j,t),\nu_j)}.
\]
Iterate with $p=4\cdot2^k$ and radii decreasing from $R$ to
$R/2$. The product of constants is finite and independent of $j$,
so
\begin{equation}\label{eq:rescaled-local-mass}
 1=\sup_{B(p_j,R/2)}v_j
 \le C_R\left(\int_{B(p_j,R)}v_j^4d\nu_j\right)^{1/4}.
\end{equation}

Use Lemma~\ref{lem:flat-package} with $P=4$ and apply
Lemma~\ref{lem:bounded-function-compactness} on its compatible
charts. The bound $D_{h_j}(v_j)\le A^{-1}$ gives the required
local energy bounds, and the nonlinearities
$A^{-1}s^3-a_js$ converge uniformly on $[-1,1]$.
Passing to a diagonal subsequence, we obtain a limit $v\in H^1_{\rm loc}(\mathbb R^4)$
with $0\le v\le1$. It satisfies $-\gamma\Delta_{\mathrm E}v=A^{-1}v^3$
weakly and is smooth by local elliptic regularity. Combining
\eqref{eq:global-dirichlet-liminf} with the energy identity in
\eqref{eq:rescaled-equation}, we obtain $D_{\mathrm E}(v)\le A^{-1}$.

The same lemma gives local $L^1$ convergence of the fourth-power
densities. Since the charts cover each bounded region, a partition of unity gives weak
convergence of $v_j^4d\nu_j$ to $v^4c\,dx$ there. Euclidean distance spheres have zero
measure, so
the integrals over $B_{h_j}(p_j,R)$ converge to those over $B_R$.
With $R=1$, it follows from \eqref{eq:rescaled-local-mass} that the limit is nonzero. By exhausting $\mathbb R^4$ and using $\int v_j^4d\nu_j=1$, we obtain
\begin{equation}\label{eq:limit-extremizer}
 0\le v\le1,\qquad 0<m:=\int_{\mathbb R^4}v^4c\,dx\le1,
 \qquad -\gamma\Delta_{\mathrm E}v=A^{-1}v^3.
\end{equation}

\paragraph{Step 7.}
It remains to justify testing the limiting equation by $v$ on
the noncompact limit. Choose $0\le\chi_R\le1$, supported in
$B_{2R}$, equal to one on $B_R$, and with
$|d\chi_R|\le C/R$. H\"older's inequality on the Euclidean annulus gives
\[
 e_R:=\gamma\int v^2|d\chi_R|^2c\,dx
 \le CR^{-2}\Vol(B_{2R}\setminus B_R)^{1/2}
       \left(\int_{B_{2R}\setminus B_R}v^4c\,dx\right)^{1/2}
 \longrightarrow0.
\]
Testing \eqref{eq:limit-extremizer} with $\chi_R^2v$, we obtain
\[
 \gamma\int\chi_R^2|dv|^2c\,dx
 +2\gamma\int\chi_Rv\langle d\chi_R,dv\rangle c\,dx
 =A^{-1}\int\chi_R^2v^4c\,dx.
\]
The absolute value of the middle term is at most
$2D_{\mathrm E}(v)^{1/2}e_R^{1/2}\to0$. Dominated convergence
in the other terms therefore gives $D_{\mathrm E}(v)=m/A$.
Moreover $\chi_Rv\to v$ in $L^4$ and
$d(\chi_Rv)\to dv$ in $L^2$, the latter by finite Dirichlet
energy and $e_R\to0$. Applying the homogeneous Euclidean inequality to $\chi_Rv$ and passing to the limit, we obtain
\[
 \sqrt m\le\CE D_{\mathrm E}(v)=\frac{\CE}{A}m.
\]
Since $m>0$, this forces $\sqrt m\ge A/\CE>1$, contrary to
$m\le1$. This excludes the assumed sequence and proves the
existence of $\delta_A>0$.
\end{proof}

\begin{remark}\label{rem:tian-wang-sobolev}
The proof of Theorem~\ref{thm:barrier} follows the argument of
Tian and Wang~\cite[proof of Proposition 3.1]{TianWang}.
They use a minimizer of the logarithmic Sobolev functional on a ball,
whereas we use an extremizer for the critical Sobolev inequality on $X$.
In both proofs, the limiting function contradicts the corresponding
sharp Euclidean inequality. Bakry, Coulhon, Ledoux, and
Saloff-Coste~\cite[Theorem 10.2]{BCLS} show that, for the same
Dirichlet energy, the critical Sobolev inequality and the dimensional
logarithmic Sobolev inequality imply each other up to a change of constants.
\end{remark}

\label{sec:flow-application}
Theorem~\ref{thm:barrier} excludes the value $A$ along any path whose
energy stays below $\delta_A$. To carry the strict bound from the
reference metric to the other endpoint, we need continuity of the
optimal Sobolev constant along the path.

\begin{lemma}\label{lem:sobolev-continuity}
The quantity $\CS(\omega_u)$ is continuous in the $C^2$ topology on
positive potentials in a fixed K\"ahler class.
\end{lemma}
\begin{proof}
If $u_j\to u$ in $C^2$ with $\omega_u>0$, then the metrics, inverse
metrics, and volume densities converge uniformly. Since the limiting metric is
positive on the compact manifold $X$, these comparisons hold up to multiplicative factors
tending to one. Thus the
Sobolev spaces agree as sets, and there is a sequence $\epsilon_j\downarrow0$ such
that the numerator and denominator in \eqref{eq:sobolev-definition}
each lie between $1-\epsilon_j$ and $1+\epsilon_j$ times their values
at $u$, uniformly in the test function for all large $j$. Taking the supremum over all
nonzero test functions, we obtain
\[
 \frac{1-\epsilon_j}{1+\epsilon_j}\CS(\omega_u)
 \le\CS(\omega_{u_j})
 \le\frac{1+\epsilon_j}{1-\epsilon_j}\CS(\omega_u).
\]
\end{proof}

\begin{proof}[Proof of Theorem~\ref{thm:sobolev-small}]
Fix $A>\CS(\omega_*)$ and let $\delta_A>0$ be the number supplied
by Theorem~\ref{thm:barrier}. Set $\varepsilon_A=\delta_A$.
For a smooth $u$ with $\Ca(u)<\varepsilon_A$, normalize $\AM(u)=0$.
Theorem~\ref{thm:connected} gives a path of smooth potentials
$\gamma:[0,1]\to\HH_0$, continuous in $C^\infty$, with
\[
 \gamma(0)=0,\quad\gamma(1)=u,\qquad
 \Ca(\gamma(s))<\delta_A\quad(0\le s\le1).
\]
The function $s\mapsto\CS(\omega_{\gamma(s)})$ is continuous by
Lemma~\ref{lem:sobolev-continuity}. It starts below $A$. If its value at the endpoint
were at least $A$, it would take the value $A$ at some
$s_0$, contradicting Theorem~\ref{thm:barrier} and the strict energy
inequality along the path. Hence $\CS(\omega_u)<A$.
\end{proof}

By the extension theorem of Li, Wang, and Zheng~\cite[Theorem 1.2]{LWZSobolev}, a smooth Calabi flow on a compact K\"ahler surface extends past
any finite time if its Sobolev constants remain uniformly bounded.
In their Definition 2.1, the authors use the same probability measure as here and the
real gradient in the Dirichlet term. Since
$|\nabla f|^2=2|\partial f|^2$, their optimal constant
$\CS^{\mathrm{real}}$ satisfies
\[
 \tfrac12\CS(\omega_u)\le\CS^{\mathrm{real}}(\omega_u)
                    \le\CS(\omega_u).
\]
Thus our uniform bound verifies precisely the extension
theorem's Sobolev hypothesis.

\begin{proof}[Proof of Corollary~\ref{cor:calabi-global}]
Choose $\varepsilon_A$ as in Theorem~\ref{thm:sobolev-small}.
Since $\omega_{u_0}$ is smooth and strictly positive on the compact manifold $X$, the
initial norm and ellipticity bounds required by Chen and He~\cite[Theorem
3.2]{ChenHeStability} hold. Let $[0,T_{\max})$ be
the maximal smooth interval of the resulting flow. It preserves the K\"ahler
class, and the Calabi-energy dissipation identity is
\[
 \frac{d}{dt}\Ca(u(t))
 =-2\int_X|\bar\partial\nabla^{1,0}S_{u(t)}|_{u(t)}^2
                  d\mu_{u(t)}\le0.
\]
See Chen and He~\cite[Section 2.2]{ChenHeStability} and
Li, Wang, and Zheng~\cite[equation (2.6)]{LWZSobolev}. Consequently, we have
\[
 \Ca(u(t))\le\Ca(u_0)<\varepsilon_A\qquad(t<T_{\max}).
\]
Theorem~\ref{thm:sobolev-small} therefore gives
$\CS(\omega_{u(t)})<A$ throughout this interval. If
$T_{\max}<\infty$, the flow extends past $T_{\max}$ by the Sobolev extension theorem, a contradiction. Thus $T_{\max}=\infty$.
\end{proof}

\vskip10pt
Haozhao Li, Institute of Geometry and Physics, and School of Mathematical
Sciences, University of Science and Technology of China, No. 96 Jinzhai
Road, Hefei, Anhui Province, 230026, China.
Email: \href{mailto:hzli@ustc.edu.cn}{hzli@ustc.edu.cn}\par\medskip

Bing Wang, Institute of Geometry and Physics, and School of Mathematical
Sciences, University of Science and Technology of China, Hefei 230026,
China, and Hefei National Laboratory, Hefei 230088, China.
Email: \href{mailto:topspin@ustc.edu.cn}{topspin@ustc.edu.cn}\par\medskip


\begin{thebibliography}{99}
\raggedright

\bibitem{AgmonDouglisNirenberg}
S.~Agmon, A.~Douglis, and L.~Nirenberg,
\emph{Estimates near the boundary for solutions of elliptic partial
differential equations satisfying general boundary conditions. I},
Comm. Pure Appl. Math. \textbf{12} (1959), 623--727.

\bibitem{BCLS}
D.~Bakry, T.~Coulhon, M.~Ledoux, and L.~Saloff-Coste,
\emph{Sobolev inequalities in disguise},
Indiana Univ. Math. J. \textbf{44} (1995), no.~4, 1033--1074.

\bibitem{BHPV}
W.~P.~Barth, K.~Hulek, C.~A.~M.~Peters, and A.~Van de Ven,
\emph{Compact Complex Surfaces}, second enlarged edition,
Ergebnisse der Mathematik und ihrer Grenzgebiete (3),
\textbf{4}, Springer-Verlag, Berlin, 2004.

\bibitem{BN}
I.~D.~Berg and I.~G.~Nikolaev,
\emph{Quasilinearization and curvature of Aleksandrov spaces},
Geom. Dedicata \textbf{133} (2008), 195--218.

\bibitem{BB}
R.~J.~Berman and B.~Berndtsson,
\emph{Convexity of the K-energy on the space of K\"ahler metrics
and uniqueness of extremal metrics},
J. Amer. Math. Soc. \textbf{30} (2017), 1165--1196.

\bibitem{BDL}
R.~J.~Berman, T.~Darvas, and C.~H.~Lu,
\emph{Regularity of weak minimizers of the K-energy and applications
to properness and K-stability},
Ann. Sci. \'Ec. Norm. Sup\'er. (4) \textbf{53} (2020), 267--289.

\bibitem{BDLflow}
R.~J.~Berman, T.~Darvas, and C.~H.~Lu,
\emph{Convexity of the extended K-energy and the large time behaviour
of the weak Calabi flow},
Geom. Topol. \textbf{21} (2017), 2945--2988.

\bibitem{BlockiGeodesics}
Z.~B{\l}ocki,
\emph{On geodesics in the space of K\"ahler metrics},
in \emph{Advances in Geometric Analysis},
Advanced Lectures in Mathematics \textbf{21},
International Press, Somerville, MA, 2012, 3--20.

\bibitem{BHJ}
S.~Boucksom, T.~Hisamoto, and M.~Jonsson,
\emph{Uniform K-stability, Duistermaat--Heckman measures and
singularities of pairs},
Ann. Inst. Fourier \textbf{67} (2017), 743--841.

\bibitem{BJGlobal}
S.~Boucksom and M.~Jonsson,
\emph{Global pluripotential theory over a trivially valued field},
Ann. Fac. Sci. Toulouse Math. (6) \textbf{31} (2022), 647--836.

\bibitem{BJDiv}
S.~Boucksom and M.~Jonsson,
\emph{A non-Archimedean approach to K-stability, II:
Divisorial stability and openness},
J. Reine Angew. Math. \textbf{805} (2023), 1--53.

\bibitem{BrezisLieb}
H.~Br\'ezis and E.~Lieb,
\emph{A relation between pointwise convergence of functions and
convergence of functionals},
Proc. Amer. Math. Soc. \textbf{88} (1983), 486--490.

\bibitem{CalabiExtremal}
E.~Calabi,
\emph{Extremal K\"ahler metrics},
in \emph{Seminar on Differential Geometry} (S.-T.~Yau, ed.),
Annals of Mathematics Studies \textbf{102},
Princeton University Press, Princeton, NJ, 1982, 259--290.

\bibitem{CC}
X.~X.~Chen and J.~Cheng,
\emph{On the constant scalar curvature K\"ahler metrics
(I)---A priori estimates},
J. Amer. Math. Soc. \textbf{34} (2021), 909--936.

\bibitem{Chen}
X.~X.~Chen,
\emph{The space of K\"ahler metrics},
J. Differential Geom. \textbf{56} (2000), 189--234.

\bibitem{ChenHe}
X.~X.~Chen and W.~He,
\emph{The Calabi flow on K\"ahler surfaces with bounded Sobolev constant (I)},
Math. Ann. \textbf{354} (2012), 227--261.

\bibitem{ChenHeStability}
X.~X.~Chen and W.~He,
\emph{On the Calabi flow},
Amer. J. Math. \textbf{130} (2008), 539--570.

\bibitem{CTW}
J.~Chu, V.~Tosatti, and B.~Weinkove,
\emph{On the \(C^{1,1}\) regularity of geodesics in the space of
K\"ahler metrics},
Ann. PDE \textbf{3} (2017), article 15.

\bibitem{Darvas}
T.~Darvas,
\emph{The Mabuchi geometry of finite energy classes},
Adv. Math. \textbf{285} (2015), 182--219.

\bibitem{DarvasCompletion}
T.~Darvas,
\emph{The Mabuchi completion of the space of K\"ahler potentials},
Amer. J. Math. \textbf{139} (2017), 1275--1313.

\bibitem{Deimling}
K.~Deimling,
\emph{Nonlinear Functional Analysis},
Springer-Verlag, Berlin, 1985.

\bibitem{DRKahler}
R.~Dervan and J.~Ross,
\emph{K-stability for K\"ahler manifolds},
Math. Res. Lett. \textbf{24} (2017), 689--739.

\bibitem{DGL}
E.~Di Nezza, V.~Guedj, and C.~H.~Lu,
\emph{Finite entropy vs finite energy},
Comment. Math. Helv. \textbf{96} (2021), 389--419.

\bibitem{DonaldsonCalabi}
S.~K.~Donaldson,
\emph{Lower bounds on the Calabi functional},
J. Differential Geom. \textbf{70} (2005), 453--472.

\bibitem{ELMNPBase}
L.~Ein, R.~Lazarsfeld, M.~Musta\c t\u a, M.~Nakamaye, and M.~Popa,
\emph{Asymptotic invariants of base loci},
Ann. Inst. Fourier \textbf{56} (2006), 1701--1734.

\bibitem{ELMNPVol}
L.~Ein, R.~Lazarsfeld, M.~Musta\c t\u a, M.~Nakamaye, and M.~Popa,
\emph{Restricted volumes and base loci of linear series},
Amer. J. Math. \textbf{131} (2009), 607--651.

\bibitem{GilbargTrudinger}
D.~Gilbarg and N.~S.~Trudinger,
\emph{Elliptic Partial Differential Equations of Second Order},
Classics in Mathematics, Springer-Verlag, Berlin, 2001.

\bibitem{HisamotoCalabi}
T.~Hisamoto,
\emph{Lower bound of the Calabi functional and non-Archimedean K-energy},
\href{https://arxiv.org/abs/2611.11938v1}{arXiv:2611.11938v1} (2026).

\bibitem{Kolodziej}
S.~Ko\l odziej,
\emph{The complex Monge--Amp\`ere equation},
Acta Math. \textbf{180} (1998), 69--117.

\bibitem{Lazarsfeld}
R.~Lazarsfeld,
\emph{Positivity in Algebraic Geometry I--II},
Ergebnisse der Mathematik und ihrer Grenzgebiete (3),
\textbf{48--49}, Springer-Verlag, Berlin, 2004.

\bibitem{LL}
M.-C.~Lee and T.-K.~Lee,
\emph{Local mollification of metrics with small curvature concentration},
Compos. Math. \textbf{162} (2026), 2675--2697.

\bibitem{LiRays}
C.~Li,
\emph{Geodesic rays and stability in the cscK problem},
Ann. Sci. \'Ec. Norm. Sup\'er. (4) \textbf{55} (2022), 1529--1574.

\bibitem{LiSemiFano}
C.~Li,
\emph{Yau--Tian--Donaldson correspondence for K-semistable
Fano manifolds},
J. Reine Angew. Math. \textbf{733} (2017), 55--85.

\bibitem{ControlledCompanion}
H.~Li, B.~Wang, and K.~Zheng,
\emph{Calabi flow in the controlled K\"ahler cone},
manuscript (2026).

\bibitem{LWZSobolev}
H.~Li, B.~Wang, and K.~Zheng,
\emph{Long-time existence of the Calabi flow on K\"ahler surfaces
with bounded Sobolev constants},
\href{https://arxiv.org/abs/2610.05217v1}{arXiv:2610.05217v1} (2026).

\bibitem{PolarizedCompanion}
H.~Li, B.~Wang, and K.~Zheng,
\emph{Calabi flow on polarized K\"ahler surfaces},
\href{https://arxiv.org/abs/2611.01968}{arXiv:2611.01968} (2026).

\bibitem{WeakCompanion}
H.~Li, B.~Wang, and K.~Zheng,
\emph{Convergence of the Calabi flow on K\"ahler surfaces
with bounded Sobolev constants},
\href{https://arxiv.org/abs/2611.01882}{arXiv:2611.01882} (2026).

\bibitem{LuSeyyedali}
Z.~Lu and R.~Seyyedali,
\emph{Remarks on a result of Chen--Cheng},
Complex Manifolds \textbf{11} (2024), article 20240004.

\bibitem{Lunardi}
A.~Lunardi,
\emph{Analytic Semigroups and Optimal Regularity in Parabolic Problems},
Progress in Nonlinear Differential Equations and Their Applications
\textbf{16}, Birkh\"auser, Basel, 1995.

\bibitem{MabuchiKenergy}
T.~Mabuchi,
\emph{K-energy maps integrating Futaki invariants},
Tohoku Math. J. (2) \textbf{38} (1986), 575--593.

\bibitem{Martens}
A.~Martens,
\emph{Sharpening a gap theorem: nonnegative Ricci and small curvature
concentration},
Calc. Var. Partial Differential Equations \textbf{64} (2025), article 79.

\bibitem{MP}
P.~Mesquita-Piccione,
\emph{A non-Archimedean theory of complex spaces and the cscK problem},
Adv. Math. \textbf{481} (2025), article 110543.

\bibitem{MPN}
P.~Mesquita-Piccione and D.~Witt Nystr\"om,
\emph{A transcendental non-Archimedean Calabi--Yau Theorem with
applications to the cscK problem},
\href{https://arxiv.org/abs/2509.09442v1}{arXiv:2509.09442v1} (2025).

\bibitem{Perelman}
G.~Perelman,
\emph{The entropy formula for the Ricci flow and its geometric applications},
\href{https://arxiv.org/abs/math/0211159}{arXiv:math/0211159} (2002).

\bibitem{Shi}
W.-X.~Shi,
\emph{Deforming the metric on complete Riemannian manifolds},
J. Differential Geom. \textbf{30} (1989), 223--301.

\bibitem{SDSemi}
Z.~Sj\"ostr\"om Dyrefelt,
\emph{K-semistability of cscK manifolds with transcendental
cohomology class},
J. Geom. Anal. \textbf{28} (2018), 2927--2960.

\bibitem{Streets}
J.~Streets,
\emph{Long time existence of minimizing movement solutions of Calabi flow},
Adv. Math. \textbf{259} (2014), 688--729.

\bibitem{Talenti}
G.~Talenti,
\emph{Best constant in Sobolev inequality},
Ann. Mat. Pura Appl. (4) \textbf{110} (1976), 353--372.

\bibitem{TianWang}
G.~Tian and B.~Wang,
\emph{On the structure of almost Einstein manifolds},
J. Amer. Math. Soc. \textbf{28} (2015), no.~4, 1169--1209.

\bibitem{TosattiWeinkoveSmallEnergy}
V.~Tosatti and B.~Weinkove,
\emph{The Calabi flow with small initial energy},
Math. Res. Lett. \textbf{14} (2007), no.~6, 1033--1039.

\bibitem{Trusiani}
A.~Trusiani, with an appendix by S.~Boucksom,
\emph{A solution to the Yau--Tian--Donaldson Conjecture through
Special Fujita Approximations},
\href{https://arxiv.org/abs/2605.30063v3}{arXiv:2605.30063v3} (2026).

\bibitem{Xia}
M.~Xia,
\emph{On sharp lower bounds for Calabi-type functionals and
destabilizing properties of gradient flows},
Anal. PDE \textbf{14} (2021), 1951--1976.

\bibitem{Yau}
S.-T.~Yau,
\emph{On the Ricci curvature of a compact K\"ahler manifold and
the complex Monge--Amp\`ere equation, I},
Comm. Pure Appl. Math. \textbf{31} (1978), 339--411.

\end{thebibliography}
\end{document}